\magnification=\magstep1

\input amstex
\documentstyle{amsppt}
\vsize 45pc                  

\input diagrams.tex

\ifx\labelsloaded\relax\endinput\else\let\labelsloaded\relax\fi 
\catcode`\@=11 
\newwrite\@auxout 
\newif\if@fwdref \@fwdreffalse 
\newcount\Ch@nro \Ch@nro=0 
\newcount\Sec@nro \Sec@nro=0 
\newcount\Tag@nro \Tag@nro=0  
\newcount\Thm@nro \Thm@nro=0 
\newcount\Fig@nro \Fig@nro=0 
\newcount\Table@nro \Table@nro=0 
\newif\ift@gsandthms \t@gsandthmsfalse 
\def\theCh{\the\Ch@nro} 
\def\theSec{\the\Sec@nro} 
\def\theTag{\the\Tag@nro} 
\def\theThm{\the\Thm@nro} 
\def\theFig{\the\Fig@nro} 
\def\theTable{\the\Table@nro} 
\def\TagsAndTheorems{\t@gsandthmstrue} 
\def\stepbyone#1{\global\expandafter 
  \advance\csname#1@nro\endcsname by1\relax} 
\def\stepcounter#1#2{\global\expandafter 
  \advance\csname#1@nro\endcsname by #2\relax} 
\def\setcounter#1#2{\global\expandafter\csname#1@nro\endcsname =#2\relax} 
\def\Ch@{\setcounter{Sec}{0}\setcounter{Tag}{0}\setcounter{Thm}{0}%
   \xdef\@currentlabel{\ifnum\Ch@nro=0 \else\theCh\fi}\@currentlabel} 
\def\Ch{\stepbyone{Ch}\Ch@} 
\def\Sec@{\setcounter{Tag}{0}\setcounter{Thm}{0}%
  \xdef\@currentlabel{\ifnum\Ch@nro=0 \else\theCh.\fi 
     \ifnum\Sec@nro=0 \else\theSec\fi}%
  \@currentlabel} 
\def\Sec{\stepbyone{Sec}\Sec@} 
\def\LeftTagForm{(} 
\def\RightTagForm{)} 
\def\Tag@{\xdef\@currentlabel{\LeftTagForm 
    \ifnum\Ch@nro=0 \else\theCh.\fi\ifnum\Sec@nro=0 \else\theSec.\fi 
    \theTag\RightTagForm}%
  \@currentlabel} 
\def\Tag{\stepbyone{Tag}%
    \ift@gsandthms\else\stepbyone{Thm}\fi\Tag@} 
\def\Thm@{\xdef\@currentlabel{\ifnum\Ch@nro=0 \else\theCh.\fi 
    \ifnum\Sec@nro=0 \else\theSec.\fi\theThm}%
  \@currentlabel} 
\def\Thm{\stepbyone{Thm}%
    \ift@gsandthms\else\stepbyone{Tag}\fi\Thm@} 
\def\Fig@{\xdef\@currentlabel{\theFig}\@currentlabel} 
\def\Fig{\gdef\CorrectCounter{\stepcounter{Fig}{-1}}%
  \stepbyone{Fig}\Fig@} 
\def\Table@{\xdef\@currentlabel{\theTable}\@currentlabel} 
\def\Table{\gdef\CorrectCounter{\stepcounter{Table}{-1}}%
  \stepbyone{Table}\Table@} 
\def\thepage{\number\pageno} 
\def\refto#1{\@ifundefined{r@#1}{{\bf ??}\global\@fwdreftrue\@warning 
   {Reference `#1' on page \thepage \space 
    undefined}}{\edef\@tempa{\@nameuse{r@#1}}\expandafter 
    \@car\@tempa \@nil\null}} 
\def\reftopage#1{\@ifundefined{r@#1}{{\bf ??}\global\@fwdreftrue\@warning 
   {Reference `#1' on page \thepage \space 
    undefined}}{\edef\@tempa{\@nameuse{r@#1}}\expandafter 
    \@cdr\@tempa\@nil\null}} 
\newif\if@filesw \@fileswtrue 
\def\label#1{\@bsphack\if@filesw {\let\thepage\relax 
   \xdef\@gtempa{\write\@auxout{\string 
      \newlabel{#1}{{\@currentlabel}{\thepage}}}}}\@gtempa 
   \fi\@esphack} 
\def\rosterlabel#1{\@bsphack\xdef\@currentlabel{%
   \therosteritem{\number\rostercount@}}%
   \if@filesw {\let\thepage\relax 
   \xdef\@gtempa{\write\@auxout{\string 
      \newlabel{#1}{{\@currentlabel}{\thepage}}}}}\@gtempa 
   \fi\@esphack} 
\def\newlabel#1#2{\@ifundefined{r@#1}{}{\@warning{Label `#1' multiply 
   defined}}\global\@namedef{r@#1}{#2}} 
\def\@currentlabel{} 
\newif\if@indexsw 
\newwrite\@idxout 
\def\makeindex{\@indexswtrue\immediate\openout\@idxout=\jobname.idx 
   \immediate\write\@idxout{\string\relax}} 
\def\index#1{\@bsphack\if@indexsw {\let\thepage\relax 
   \xdef\@gtempa{\write\@idxout{\string 
      \indexentry{#1}{\thepage}}}}\@gtempa 
   \fi\@esphack} 
\def\indexentry#1#2{\noindent {#1}{\unskip\nobreak\hfil\penalty50 
   \hskip2em\hbox{}\nobreak\hfil#2%
   \parfillskip=0pt \finalhyphendemerits=0 \par}} 
\def\make@ref#1{\@ifundefined{b@#1}{{??}\global\@fwdreftrue\@warning 
   {Citation `#1' on page \thepage \space 
    undefined}}{\@nameuse{b@#1}}} 
\def\citeto#1{\relaxnext@ 
 \def\nextiii@##1,##2\end@{\immediate\write\@auxout{\string\citation{##1}}%
   [{\bf\make@ref{##1}},##2]}%
 \in@,{#1}\ifin@\def\next{\nextiii@#1\end@}\else 
   \def\next{\immediate\write\@auxout{\string\citation{#1}}%
     [{\bf\make@ref{#1}}]}\fi 
     \next} 
\def\noOf#1{\global\advance\c@num by1 
   \if@filesw \immediate\write\@auxout 
   {\string\bibcite{#1}{\the\c@num}}\fi\no\make@ref{#1}} 
\def\keyOf#1#2{\if@filesw 
   {\def\protect##1{\string ##1\space}\immediate 
   \write\@auxout{\string\bibcite{#1}{#2}}}\fi 
   \key{$\lbrack${\bf #2}$\rbrack$}} 
 
\let\@@end=\end 
\let\@@supereject=\supereject 
\def\supereject{\@@supereject} 
\def\end{\if@fwdref 
  \typeout{Undefined forward references/citations. Rerun.}\fi 
  \@@end} 
\newwrite\@unused 
\def\typeout#1{{\let\protect\string\immediate\write\@unused{#1}}} 
\def\@warning#1{\typeout{Warning: #1.}} 
\def\@namedef#1{\expandafter\def\csname #1\endcsname} 
\def\@nameuse#1{\csname #1\endcsname} 
\def\@car#1#2\@nil{#1}               
\def\@cdr#1#2\@nil{#2} 
\long\def\@ifundefined#1#2#3{\expandafter\ifx\csname 
  #1\endcsname\relax#2\else#3\fi} 
\def\@ifnextchar#1#2#3{\let\@tempe #1\def\@tempa{#2}\def\@tempb{#3}\futurelet 
    \@tempc\@ifnch} 
\def\@ifnch{\ifx \@tempc \@sptoken \let\@tempd\@xifnch 
      \else \ifx \@tempc \@tempe\let\@tempd\@tempa\else\let\@tempd\@tempb\fi 
      \fi \@tempd} 
\def\?{\let\@sptoken= } \?  
\def\?{\@xifnch} \expandafter\def\? {\futurelet\@tempc\@ifnch} 
\newdimen\@savsk 
\newcount\@savsf 
\def\@bsphack{\@savsk\lastskip 
    \ifhmode\@savsf\spacefactor\fi} 
\def\@esphack{\relax\ifhmode\spacefactor\@savsf 
     {}\ifdim \@savsk >\z@ \ignorespaces\fi\fi} 
 
\newcount\c@num \c@num=0 
\def\bibcite#1#2{\@ifundefined{r@#1}{}{\@warning{Referencelabel `#1' multiply 
   defined}}\global\@namedef{b@#1}{#2}} 
\def\@gobble#1{} 
\let\citation\@gobble 
\def\CorrectCounter{} 
\def\ins@#1{\relaxnext@ 
 \smallcaptionwidth@\captionwidth@\gdef\thespace@{#1}%
 \def\next@{\ifx\next\space@\def\next. {\futurelet\next\nextii@}\else 
  \def\next.{\futurelet\next\nextii@}\fi\next.}%
 \def\nextii@{\ifx\next\caption\def\next\caption{\futurelet\next\nextiii@}%
  \else\let\next\nextiv@\fi\next}%
 \def\nextiv@{\vnonvmode@ 
  {\ifmid@\let\next\midinsert\else\let\next\topinsert\fi 
  \next\vbox to\thespace@{}\endinsert} 
  {\ifmid@\nonvmodeerr@\midspace\else\nonvmodeerr@\topspace\fi}}%
 \def\nextiii@{\ifx\next\captionwidth\let\next\nextv@ 
  \else\let\next\nextvi@\fi\next}%
 \def\nextv@\captionwidth##1##2{\smallcaptionwidth@##1\relax\nextvi@{##2}}%
 \def\nextvi@##1{\def\thecaption@{##1}%
  \def\next@{\ifx\next\space@\def\next. {\futurelet\next\nextvii@}\else 
   \def\next.{\futurelet\next\nextvii@}\fi\next.}%
  \futurelet\next\next@}%
 \def\nextvii@{\vnonvmode@ 
  {\ifmid@\let\next\midinsert\else 
  \let\next\topinsert\fi\next\vbox to\thespace@{}\nobreak\smallskip 
  \setbox\z@\hbox{\noindent\ignorespaces\thecaption@\unskip}%
  \ifdim\wd\z@>\smallcaptionwidth@\centerline{\vbox{\hsize\smallcaptionwidth@ 
  \gdef\CorrectCounter{}\unskip}}%
  \else\centerline{\box\z@}\fi\endinsert} 
  {\ifmid@\nonvmodeerr@\midspace 
  \else\nonvmodeerr@\topspace\fi}}%
 \futurelet\next\next@} 
\openin1 \jobname.aux \ifeof1 \typeout{No file \jobname.aux.}\closein1 
    \else\closein1 \relax\input \jobname.aux \fi 
\immediate\openout\@auxout=\jobname.aux 
\immediate\write\@auxout{\string\relax} 
\catcode`\@=\active 


\NoBlackBoxes
\TagsAsMath
\define\PP{\Bbb P}
\define\CC{\Bbb C}

\define\QQ{\Bbb Q}
\define\RR{\Bbb  R}
\define\ZZ{\Bbb Z}

\define\A{\Bbb{A}}
\define\spec{\text{Spec }}
\define\proj{\text{Proj }}
\define\enddemoo{$\square$ \enddemo}
\define\LL{\Cal L}
\define\MM{\Cal M}

\define\rmap{\dashrightarrow}

\NoRunningHeads
\document
\topmatter
\title   Rational and Non-rational Algebraic Varieties: \\
Lectures of J\'anos Koll\'ar 	
\endtitle
\rightheadtext{Rational Varieties}
\author  
By Karen E. ~Smith \\
with an Appendix by Joel Rosenberg
\endauthor

\address{J. Koll\'ar, Dept. of Mathematics, University of Utah, Salt Lake City UT 84112 }
\endaddress 
\email {kollar\@math.utah.edu}
\endemail
\address
{K.E. Smith, Dept. of Mathematics,
Massachusetts Institute of Technology,
Cambridge, MA 02139}
\endaddress
\email {kesmith\@math.mit.edu}
\endemail 

\thanks {These notes are based on a course given by 
J\'anos Koll\'ar at the  European Mathematical Society Summer School 
   in Algebraic Geometry at Eger, Hungary,
July 29- August 9, 1996.
These notes are in draft form only, and the author takes responsibility for
all errors. 
Please let me know if you find any!
Special thanks to Paul Taylor's Diagram program and Joel Rosenberg 
for help using it.}
\endthanks
\date  July 14, 1997
\enddate
\endtopmatter

\setcounter {Ch}{0}
\heading   Introduction \endheading

Rational varieties are among the simplest possible algebraic varieties.
Their study  is as old as algebraic geometry itself, 
yet it remains a remarkably difficult area of research today. 
  These notes offer an introduction to the study of rational varieties.
We begin with 
the  beautiful classical geometric 
approach to  finding examples of  rational varieties,
and end with some  subtle algebraic
arguments that   have  recently  established
 non-rationality of  varieties that otherwise 
share many of their traits.

In lecture one, 
rationality  and unirationality
are defined, and illustrated with a series of examples.
We also introduce some
easily computable invariants, the  plurigenera,
that vanish for all rational varieties. Using  the plurigenera,
one develops a sense how rare rational varieties are.
For example,  we will  immediately see  that no smooth projective hypersurface
whose degree exceeds its embedding dimension 
can be rational: no plane cubic curve is rational, no space quartic surface
is  rational,
and so on.

The second and third lectures focus on the rationality question for
smooth 
cubic surfaces over arbitrary fields, an
issue  thoroughly explored by 
 B. ~Segre in the forties.
This 
is a pretty story, depending 
  subtly  on the field of definition. It was already understood 
one hundred years ago that every  cubic surface
is rational over the complex numbers; however, the situation is 
quite complicated over
the rational numbers. 
In the second lecture, we construct examples of non-rational smooth 
cubics over $\QQ$ by considering 
 the orbits of 
 the Galois group of $\bar \QQ/\QQ$ on the
 twenty seven lines on the cubic surface.
The construction makes use of a beautiful theorem of Segre:
no smooth cubic surface of Picard number one can be rational. 
  In the third lecture, 
we prove  this theorem of Segre. Using the same techniques, we 
also prove the following stronger theorem of 
Manin:
 two smooth  cubic surfaces of Picard number one are birationally equivalent
if and only if they are projectively equivalent.

The cubic surface case  
is especially  interesting in light of the fact that 
there exist no smooth rational cubic threefolds.
 This longstanding--- and sometimes controversial---  question
 was  settled  in the early seventies by
Clemens and Griffiths; see 
Section \refto{history}.
 Still  to this day, no one
knows whether or not there exists a smooth  rational quartic fourfold.

When the obvious numerical obstructions
 to rationality vanish,
it can be difficult to determine whether or not a variety is rational.
The purpose of the final two lectures is  to show that indeed,
there are abundant examples of non-rational, and even non-ruled,
 varieties of every dimension  that would  appear to be rational
from the point of view of na{\"\i}ve numerical invariants.
The examples of non-ruled {\it smooth Fano varieties\/}
(definition \refto{fanodef}) are 
 defined over fields of characteristic zero, 
but the only known proof that they are non-ruled
uses the technique of reduction to prime characteristic.
The characteristic $p$ argument, presented in Lecture 4, relies
heavily on the special nature of differentiation in prime characteristic.
In Lecture 5, the characteristic zero theorems are deduced from the
characteristic $p$ results.

These techniques are applied in 
an appendix
 written by Joel Rosenberg to construct 
new 
explicit examples of varieties defined over $\QQ$ that are Fano but not 
ruled.

The notes  have been written with the goal of making them 
accessible to 
students 
with a basic training in algebraic geometry at the level of
 \cite
{H}. The first lecture is to be relatively
easy, with 
subsequent lectures requiring more of the reader.
The first three lectures are in the realm of classical algebraic geometry,
while a scheme theoretic approach is indispensible for the last  lecture.

Finally, included are 
detailed 
solutions to nearly all the exercises scattered throughout the notes.
The exercises formed an important part of the summer course, as
they form an important part of the written version. Many of the
solution ideas were worked out together with the summer school participants,
especially S\'andor Kov\'acs,
 G\'abor Megyesi, and   Endre Szab\'o.
Finally, special thanks are due to 
K\'aroly B\"or\"oczky Jr.  for organizing the summer school,
and to the European Mathematical Society for funding it.

\subhead{Notation}\endsubhead
 
Varieties defined over {\it non-algebraically closed\/} fields
occupy a central position in our study. 
The notation for the ground field is suppressed when the field
is clear from the context or irrelevent, but
occasionally  we write $X/k$ to emphasize that the variety $X$ is defined
over the ground field $k$.  
The algebraic closure of $k$ is denoted $\bar k$.

Varieties are 
 assumed reduced and 
irreducible, except where explicitly stated otherwise.
Because we are concerned with birational properties,  
there is  no loss of generality in assuming all varieties
 to be quasi-projective.
In any case, our main interest is in smooth projective varieties.

Morphisms 
and rational maps  between varieties are always assumed
to be defined over the ground field, except where explictly stated otherwise.
Likewise, linear systems on a  variety  $X/k$ are assumed defined 
over $k$.

 Morphisms are denoted by solid arrows $@>>>$ and rational
maps by dotted arrows $\rmap$.
 The ``image'' of a rational map is 
the closure of the image of the   morphism obtained by restricting 
the rational map to some non-empty open set where it
is defined;  in the same way, we define the image of
a subvariety under a rational map, provided that the map is defined
at its generic point. 
In the case of birational maps, the image of a subvariety is also called the
birational transform.  Likewise, the graph of a rational
map $X \overset{\phi}\to\rmap Y $
is the closure in $X \times Y$  of the graph of the 
morphism obtained by restricting  $\phi$ to a non-empty open set where it
is defined.

If  $L$ is any  field containing $k$,
the symbol $X(L)$ denotes  set of 
the $L$-rational points, or $L$-points for short,  of $X$.

The final lecture  deals
 with  schemes  quasi-projective over an arbitrary affine base scheme
$S$. In this case, all  morphisms and all rational maps
are assumed defined over $S$. Likewise, all products are defined over $S$.

\setcounter {Ch}{0}
\heading{\Ch.   First Lecture }
\label{lec1}
 \endheading

We begin with the precise definitions of rationality and unirationality.

\proclaim{\Sec. Definition} \label{ratdef}
A  variety $X$ is rational if there exists a birational map
$
{\PP^n \overset{\phi}\to\rmap   X,} 
$  that is, if  $X$ is 
birationally equivalent to projective space.
\endproclaim

According to the  conventions agreed upon in the introduction, 
 implicit in the above definition is an unnamed
ground field $k$ and the birational map $\phi$ is assumed to be defined
over $k$. Occasionally, we will  emphasize the ground field
by  saying that  ``$X$ is rational over $k$.'' 
Of course, if $X$ is defined over $k$ and $L$ is any 
field extension of $L$, then $X$ may also be considered to be defined
over $L$.
If $X$ is rational over $k$, then $X$ is also rational over $L$.

\proclaim{\Sec. Definition} \label{unidef}
A  variety $X$ is  unirational if there exists a 
generically finite dominant rational map 
$  \PP^n  \overset{\phi}\to\rmap X$.
\endproclaim

Roughly speaking, a variety is {\it unirational\/} if a dense
open subset is 
parameterized  by projective space, and
 {\it rational\/}
 if such a  parametrization is one-to-one.

\smallskip
The  disconcerting 
use of the prefix ``uni"  in referring to a map which is 
finite-to-one  makes more sense when viewed in historical context.
Rational varieties were once called ``birational," in reference to the
rational maps between them and projective space in each direction.
``Unirationality" thus refers to the  map from $\PP^n$ to the variety,
defined in one direction only.

\medskip

These lectures  treat the following general question:
{{{\it{Which varieties are rational or unirational?}}}}

\medskip

It is important to realize that the rationality or unirationality 
of a variety may depend on what we take to be the field of definition.
For example, a variety $X$ defined over $\QQ$ may be considered as
a variety defined over $\RR$. It is possible that there
is  a  birational map 
given by polynomials with {\it real\/} coefficients  from projective
 space to $X$,
 but there is no such  birational map
from given by polynomials with {\it rational\/}
coefficients.
Our first example nicely illustrates this point.

\smallskip

{\bf{\Sec. Quadrics.\/}} \label{conic}
A smooth quadric hypersurface in projective space
 is rational  over $k$
  if and only if it  has at least one $k$-point.
In particular, every smooth quadric over an algebraically closed field
is rational.

\demo{Proof}
If  a variety $X$ is unirational over $k$, then it has a $k$-point.
This obvious when $k$ is  infinite:
the map 
$\PP^n \overset{\phi}\to\rmap X$ is defined on some  Zariski open subset of
$\PP^n$ and because
 $\PP^n$ has plenty of $k$-points in every open set,
 the image of any one of them 
under $\phi$ will be a $k$-point of $X$.
This is not so obvious when $k$ is finite, because $\PP^n$ has open sets
with no $k$-points. The finite field case follows from Nishimura's Lemma
(Exercise 1 below).

Conversely, let $X$ be a smooth quadric in $\PP^{n+1}$, 
defined  over $k$, and  with a  $k$-point
$P \in X(k)$. 
 Let $\phi $ be the
 projection from $P$ onto any $n$-plane not containing $P$.
  Choosing coordinates so that
  $P = (0: 0: \cdots :0: 1)$ , we have
$$
 \PP^{n+1} \overset{\phi}\to\rmap \PP^n
$$
$$
(x_0 : \cdots : x_{n+1}) \mapsto (x_0: \cdots : x_n).
$$
Restricting $\pi$ to $X$,  we expect a generically 
one-to-one map of $X - \{P\}$ to
 $\PP^n$;
the exceptions occur when an entire line on $X$ is collapsed to a point by
 the projection. 
But most points
 can not lie on lines that are collapsed by this map, as this 
would force $P$ to be a singular point of $X$.
\enddemoo

\smallskip
{\bf{Exercise 1:}}
Prove
 Nishimura's Lemma:   If $Y$ is smooth, $Y'$ is projective, 
 and there is a rational map 
$Y \rmap Y'$, then if $Y$ has a  $k$-point, so does  $Y'$.  
Also, find a counterexample 
when $Y$ is not smooth.

 \bigskip

Every  smooth cubic surface  in 
 $ \PP^3$ defined over an algebraically closed field
is rational, since it is isomorphic to  the blowup of $\PP^2$ at six points
\cite{H, p 395}.
The rationality question for 
 cubic surfaces over a non-algebraically closed field
is more subtle. This will be our main focus in lectures 
two and three. For now, we discuss a few simple examples to indicate
their richness.

\subhead{\Sec. A non-rational cubic surface}\label{nonrat}\endsubhead
Let  $X$ be a cubic surface in real  projective three space,
 defined by an equation  (in affine coordinates) 
 $x^2 +  y^2 = f_3(z)$, where $f_3$ has three distinct real roots. 
Then $X$  is not rational  over $\RR$.

\demo{Proof}
Consider the graph of  $f_3$. The set of points in the source
 $\RR$ where $f_3$ takes positive values has two disjoint components. 
The equation 
 $x^2 + y^2  = f_3(z) $   has  real solutions  (in fact, a circle's worth)
 if and only if  $f_3 \geq 0$, so  we see that as a real manifold $X(\RR)$
 has two distinct components.  
But if 
$X$ is birationally equivalent to $\PP^2$
over $\RR$, then because
 $\PP^2(\RR)$ is connected as a real manifold,  so would be  $X(\RR)$.
\enddemoo

However, it turns out that the cubic surface $X$ defined 
above  is unirational over $\RR$. We leave the reader the pleasure of
finding a  map $\PP^2_{\RR} \rmap X$ that is two-to-one onto one 
of the manifold
 components and misses the other 
 component entirely.  Indeed, the preimage
of a point in the missed manifold component  can be interpreted as
 a pair of 
 complex conjugate points 
in the complex manifold $\PP^2(\CC)$.

\medskip
\subhead{\Sec. Singular cubics}\label{singrat}\endsubhead
An irreducible  cubic hypersurface  in projective space
(that is not a cone over a cubic hypersurface of lower dimension)
is rational over $k$ if 
it   has  a  {\it singular\/} $k$-point.

\demo{Proof}  
Let $X$ be the cubic hypersurface in $\PP^{n+1}$, defined over the ground field
$k$.
Project from the singular 
point  $P \in X(k)$ onto a general
hyperplane defined over $k$.  Since
$P$ has multiplicity two on $X$,  any line through 
$P$ 
has a unique third point of intersection with $X$. Its projection onto the
hyperplane  gives the one-to-one map from $X$ to $\PP^n$.
Of course, this makes sense only when the line through $P$ 
does not lie on $X$. This is where we use 
the assumption that $X$ is not a 
cone: because $X$ is not a cone, the generic line through 
$P$ does not lie on $X$. 
\enddemoo

\subhead{\Sec. Rationality of cubic hypersurfaces}\label{linspac}\endsubhead
If a smooth cubic hypersurface of even dimension
 contains
two  disjoint
linear spaces, each of half the dimension,
 then the cubic hypersurface is 
 rational.  In particular, a smooth cubic surface 
is rational over $k$  if  it contains
two  skew lines  defined over $k$ (of the twenty seven
lines on the surface defined over $\bar k$). 

\demo{Proof}
Let $X \subset \PP^{2n+1}$ be the 
cubic hypersurface, and let $L_1$ and $L_2$ be
 the two linear spaces on $X$.
Consider the map 
$$
L_1 \times L_2 \overset{\phi}\to\rmap  X
$$
$$
(P, Q) \mapsto {\text{third intersection point }} X \cap \overline{PQ}.
$$
This defines a birational map from $L_1 \times L_2$  to $X$. 
The map is well defined because each line intersects $X$ in exactly three
points (counting multiplicities).  
This map is birational: if the pre-image of $x \in X$ includes
 two distinct pairs $(P_1, Q_1)$
and $(P_2, Q_2)$ on $L_1 \times L_2$, then the  projections of 
the linear spaces $L_1$ and $L_2$ from $x$ onto a general hyperplane  would 
intersect each other in more than one 
point, which  is impossible
(see  \refto{mapdisc} for a more general discussion).
Because 
$$
 {\PP^{2n}} 
\rmap {\PP^n \times \PP^n} \rmap 
L_1 \times L_2 \rmap X 
$$
are birational equivalences, we conclude that $X$ is rational. Note that
all maps above are defined over the ground field $k$. 
\enddemoo

{\bf{Exercise 2:} }
\roster
\item
Find examples of smooth 
 cubic  hypersurfaces in $\PP^{2n+1}$ containing two disjoint $n$-planes.
\item What is the dimension of the variety of all such cubics?
\item
Why have we  not considered linear spaces of non-equal dimension?
\endroster

\medskip
\subhead{\Sec. Discussion}\label{mapdisc}\endsubhead
More generally,  
given  any two  subvarieties, $U$ and $V,$
 of a degree three hypersurface $X$,  
one is tempted to form a similar map:
$$
\phi:  U  \times  V \rmap  X
$$
$$
(u, v)  \mapsto {\text{third intersection point }} X \cap \overline{uv}.
$$
If $U$ and $V$ are disjoint, this  map is a  morphism except at pairs of
points $(u, v)$ spanning a line on $X$; in  general, 
 it is not defined on $U \cap V$.

 The  map $\phi$ can not be dominant unless  $\dim U + \dim V \geq \dim X$
 and
 it can not be generically finite unless $\dim U + \dim V= \dim X$.
When $\phi$ is finite,    how does one compute its degree ?

To determine the  pre-image of a general point $x \in X$, consider the
projection $\pi_x$ from $x \in X$ to a general hyperplane.
The set $\pi_x(U) \cap \pi_x(V)$ consists of all points $\pi_x(u)= \pi_x(v)$,
  with $u, v, $ and $x$ collinear. In this case, assuming
  that $u \neq v$, the points
  $(u, v) \in U \times V$ are the pre-images of $x$ under $\phi$. 
So, if $U \cap V = \emptyset$,  we expect that the 
degree  of 
$\phi$ is the cardinality of  $\pi_x(U) \cap \pi_x(V) $. More generally, 
we must subtract something for the intersection
points of $U$ and $V$.

In Example 
 \refto{linspac}, we applied this idea with $U$ and $V$ linear subspaces
and deduced that cubic hypersurfaces are rational if they contain 
two disjoint linear  subvarieties of half the dimension.
More generally, the 
idea is useful  for detecting {\it unirationality\/} of some cubics,
 as the next example shows.

\subhead{\Sec. Unirationality of Cubic Surfaces}\label{unirational}
\endsubhead It is easy to see that 
 a  smooth 
cubic surface in $\PP^3$
containing  two  non-coplanar rational curves is 
unirational. (We remind the reader that  implicit in this statement
is that both curves are defined  over the ground field $k$,
and that the  surface is unirational over $k$.) 

Indeed, let $C_1$ and $C_2$ be rational 
curves on the surface 
$X$, and define the map $C_1 \times C_2 \overset{\phi}\to\rmap X$
as above. 
Because $C_1$ and $C_2$ do not lie in the same plane, their join
(meaning the locus of points lying on lines joining points on $C_1$ to 
points on $C_2$) must be all of $\PP^3$. This ensures that the map $\phi$ is
dominant, and hence finite. Because $C_1$ and $C_2$ are rational (over $k$),
we conclude that $X$ is unirational (over $k$).

\subsubhead{\Thm}\label{uni}\endsubsubhead
Using \refto{unirational}, we can easily deduce that 
 a sufficiently general  cubic  surface 
containing two $k$-points will be unirational.
Indeed, two rational curves $C_1/k$ and $C_2/k$ can be found by 
intersecting $X$ with the tangent plane at each of the two $k$-points. 
Assuming both $C_1$ and $C_2$ are irreducible over $\bar k$
 (as usually happens),
  each is a plane cubic curve with a singular $k$-point, and 
hence rational 
 by Example \refto{singrat}.
Furthermore, $C_1$ and $C_2$ are not coplanar; otherwise, their union is a
plane section of the cubic surface, and hence a plane cubic,  so it
could have at most one singular point.  Using the map 
$C_1 \times C_2 \overset{\phi}\to\rmap X$ defined in \refto{mapdisc},
we conclude that $X$ is unirational. Finally, note that $\phi$ is degree
six--- the  plane projections of 
$C_1$ and $C_2$  are both cubic plane curves, so they 
intersect in nine points,
 but three of these nine points come from 
 the intersection
points of $C_1$ and $C_2$.

 Even in the degenerate case where $C_1$ or $C_2$ is reducible over $\bar k$,
the argument often  goes through unchanged. 
If either  $C_1$ or $C_2$ is a union of a line and an irreducible
 quadric, then each of these components is defined  and rational over $k$,
and the above argument goes through using these rational curves---
 assuming they are not coplanar. The only way in which they can be
coplanar is when $C_1 = C_2$ is a plane section of $X$ consisting of
three lines intersecting in three distinct points. In this case,
the line through the  original two $k$-rational points is
 a $k$-rational line on $X$, and so at least 
when $k$ is infinite  it contains  many $k$-rational 
points on $X$.  Applying the argument to a generic pair of these points,
we again see that $X$ is  unirational.

There is one   honest exceptional case  where this argument breaks down:
the tangent section curve $C_1$ or $C_2$ could  be a
 union of three lines with none of these lines  defined over $k$.
The only way this can happen is when  the three lines meet in 
a single $k$-rational point, which actually happens in 
some interesting examples.
A point $P$ on a cubic surface such that $T_PX \cap X$
 is a union of three lines is called an Eckardt point.
A cubic surface can have at most finitely many Eckardt points, since they
occur exactly when three of the twenty seven lines on the 
surface intersect in a single point. Indeed, cubic
 surfaces with Eckardt points are rather special among all cubics;
see \cite{Ec}.

\subsubhead{\Thm}\endsubsubhead
The above argument shows that 
a smooth cubic surface containing a single $k$-rational point (that is not an
Eckardt point) is unirational, at least over an infinite field.
 The point is that that tangent plane to 
this point intersects the surface in a 
singular cubic, giving rise to a rational curve on $X$.
This curve contains plenty of
$k$-points, and so we can apply  \refto{uni}. This argument is due to 
B. Segre \cite{S43}. In fact, Segre later 
 showed that  if a
smooth cubic surface (over an infinite field $k$) contains a $k$-point, 
then it contains   infinitely many $k$-points
 \cite{S51}.  It follows from the argument described above that 
any smooth cubic surface  containing a $k$-point is unirational.

\subsubhead{\Thm}\endsubsubhead
An interesting variation on the 
map discussed in \refto{mapdisc} is when we allow 
$U = V$.
For example,  suppose that $X$  is a smooth cubic 
four-fold in $\PP^5$
containing a smooth surface $S$.

Consider the map 
$$
\frac{S \times S}{\sim}  \overset{\phi}\to\rmap  X
$$
$$
(P, Q) \mapsto  {\text{third point of intersection }}
 X\cap \overline{PQ}.
$$
Here, $\frac{S \times S }{\sim}$ is the symmetric product of $S$, 
the quotient variety of $S \times S$ 
by the action of the two-element group interchanging the factors.
If $S$ is unirational over $k$, then so is $S \times S$, and hence so is the
 image $\frac{S \times S}{\sim}$ under the generically two-to-one quotient map.

 Consider a general point $x \in X$,
say not on $S$. When is $x$ in the image of $\phi$? 
Consider the family of lines $\{\overline{sx}\}_{s \in S}$. 
The point $x$ is in the image of $\phi$ precisely when 
at least one of these lines intersects $S$ in a point other than $s$.
In particular, the projection from $x$,
 $S @>{\pi_x}>>  S' \subset \PP^4$ can not
be one-to-one. 
 Indeed, $x$ has a unique
pre-image under $\phi$ precisely when 
the projection
$\pi_x$ collapses  exactly
 two points of $S$ to a single point.
On the other hand, if $S @>{\pi_x}>> S'$ is not of degree one, 
then $x$ will have infinitely many preimages under $\phi$.

Thus $\phi$ is finite and dominant if and only if the generic projection 
of $S$ from a point $x \in X$ is one-to-one except on a finite set.
The next exercise provides one case where this condition can be verified.

\smallskip

{\bf{Exercise 3.} } 
Find examples of smooth surfaces in $\PP^5$ such that 
the generic projection to $\PP^4$ has exactly one singular point.
 Use this to give some more examples of unirational cubic four-folds.
(Hint: Consider  the linear system of plane cubics
through four points.)

\bigskip

\subhead{\Sec. Numerical Constraints} \endsubhead
Rationality and unirationality force strong numerical constraints 
on a variety. Let $\Omega_{X}$ be the sheaf of differential forms (K\"ahler
differentials) of $X$ over $k$.

\proclaim{\Sec. Theorem} \label{diffvan}
If a smooth projective variety $X$  is rational, then 
 $H^0(X, (\Omega_X)^{\otimes m})  $ is zero  for all $ m \geq 1$.
The same holds for unirational  $X$, provided the ground field 
has characteristic zero. 
\endproclaim

\demo{Proof}
Suppose we have a generically finite, dominant map 
$\PP^n \overset \phi\to \rmap X$.
 Let $U \subset \PP^n$ be an open set over which $\phi$ is defined;
 its complement may be assumed to have codimension
at least two.

Non-zero differential 
 forms on $X$ pull back to 
non-zero differential forms on $U$,  that is, we have an inclusion 
$ (\phi^*\Omega_X^1)^{\otimes m} \hookrightarrow (\Omega^1_U)^{\otimes m}$.
This is obvious when $\phi$ is birational,
and easy to check when $\phi$ is finite (assuming $k$ has characteristic zero).
Because the complement of $U$ has codimension at least two,  
the differential forms
 on $U$ extend  uniquely to forms on $\PP^n$, 
 so that 
$$
H^0(U, (\Omega^1_U)^{\otimes m}) = 
H^0(\PP^n, (\Omega^1_{\PP^n})^{\otimes m}).
$$
Therefore $H^0(X, (\Omega^1_X)^{\otimes m}) \subset
 H^0(\PP^n, (\Omega_{\PP^n}^1)^{\otimes m}),$
  and the problem is reduced to
 proving the vanishing for $\PP^n$, left as  Exercise 4.
\enddemoo

{\bf{Exercise 4.} }
 Complete the proof  by showing that
  $H^0(\PP^n, (\Omega_{\PP^n}^1)^{\otimes m}) $ is zero.

\subsubhead{\Thm. Remark}
\endsubsubhead
In prime characteristic, unirationality of $X$ does not necessarily 
force the vanishing of the invariants $H^0(\Omega_X^{\otimes m})$.
Indeed, the pull back map for differential forms can be the zero map, 
so  the argument above  fails.
For example, 
consider the Frobenius map
$F$ on $\A^n$ sending  
 $(\lambda_1, \dots, \lambda_n) \mapsto
(\lambda_1^p, \dots, \lambda_n^p)$, where $p>0$ is the characteristic of
 the ground field. The induced map of differential forms
$F^*\Omega \rightarrow \Omega$  sends every differential
 $dx$ to $d(x^p) = px^{p-1}dx = 0$.

In characteristic $p$, we
are led to the more sensible 
 notion of {\it separable unirationality.\/}
A variety $X$ is  separably unirational if 
there is a dominant generically finite map  $\PP^n \rmap X$
such that the induced inclusion of function fields is {\it separable.\/} 
In other words, separably unirational is equivalent to ``unirational by a
 generically \'etale map." 
 By the definition of \'etale, the pull back map of
differential forms is 
injective (in fact, an isomorphism) on a dense open set,
 so the proof of \refto{diffvan}
 shows that $H^0(X, \Omega_X^{\otimes m}) = 0$ for a  smooth projective
 separably
unirational variety $X$ of arbitrary characteristic.

\proclaim{\Sec. Corollary} \label{canvan}
If $X$ is  a smooth projective variety that is separably
unirational (for example,  rational),
then  $H^0(X, \Cal O(mK_X)) = 0$ for all $m \geq 1$.
\endproclaim 

Here $K_X$ denotes the canonical class of $X$, that is, the divisor
class representing the canonical sheaf $\omega_X = \wedge^n\Omega_X$ of
highest differential forms.
The dimension of the $k$-vector space $H^0(X, \Cal O(mK_X))$
is called the $m^{th}$ plurigenus of $X$, and is denoted $h^0(mK_X)$.
The plurigenera are  easily computable obstructions to rationality.

\bigskip
It is conjectured that for a smooth projective variety, the plurigenera
are the {\it only\/} obstructions to a related property of varieties
called {\it uniruledness;\/}  see Definition \refto{ruled}. 
A variety for which all plurigenera
vanish is said to have {\it negative Kodaira dimension\/}  (or
Kodaira dimension  $-\infty$). 
An understanding of these varieties is an essential feature
of Mori's Minimal Model Program for  the birational classification
of algebraic varieties; see   \cite{KaMM}. Miyaoka proved that 
 uniruledness is equivalent to negative Kodaira dimension
for smooth three-folds of characteristic zero \cite{Mi}; see also \cite{S-B}.
The conjecture remains open in higher dimension.
  
\medskip
A similar conjecture, attributed to Mumford, predicts  that the vanishing
of $H^0(X, \Omega_X^{\otimes m})$ for all
$m \geq 0$ 
is the only obstruction to a   smooth projective variety $X$ 
being {\it rationally connected, \/} at least in characteristic zero.
A variety is said to be rationally connected if every pair of points
can be joined by a rational curve.
 See \cite{K96, p202} for a
discussion of this conjecture, which has recently 
been proved in dimension three by 
  Koll\'ar, Miyaoka and Mori \cite{KoMM}.

\bigskip

{\bf{Exercise 5.} }
  \roster 
\item  Prove  Corollary \refto{canvan}.
\item
Show that the plurigenera of a 
 smooth hypersurface of degree $d$ in $\PP^n$
 do not vanish when $d > n$.
  Conclude that no smooth hypersurface whose degree exceeds its
embedding dimension is separably unirational. 
\endroster
\medskip

Of course, 
 in characteristic zero, separable  
unirationality  is the same as unirationality.  In prime characteristic, 
unirational but not separably unirational varieties exist, and in
fact, there are unirational hypersurfaces of arbitrary degree:

\medskip
{ \bf {Exercise  6.}} 
\roster
\item 
Show that a purely inseparable cover of a unirational variety
over a perfect field is unirational.
\item Show that there exist unirational hypersurfaces of arbitrarily large
 degree
relative to their embedding dimension.
\endroster
\medskip 

The next four exercises provide more examples of rational and non-rational 
varieties. 

\medskip
{\bf{Exercise 7.}}
\roster
\item
Let $k = \CC(t)$, and let $X$ be a degree $d$ hypersurface in 
 $\PP^n$ defined over $k$.
Prove that if $d \leq n$, then  $X$ has at least one $k$ point.
Find an example with 
exactly one $k$-point. For $d > n$, find a hypersurface with no $k$-points.
Explain why such a hypersurface is non-rational.
\item
 Same as (1), but with  $k = \CC(t, s)$ and $d^2 \leq n$.
\endroster

\medskip
{\bf{Exercise 8.} } 
\roster
\item
Let $X_{a, 2} \subset \PP^1 \times \PP^n$  be a smooth hypersurface
of bi-degree $(a, 2)$.    For  $n \geq 2$, 
show that $X_{a, 2}$ is rational over $\CC$.
\item 
Let $X_{a, 2} \subset \PP^2 \times \PP^n$ 
be a smooth hypersurface of bi-degree $(a, 2)$.
For $n \geq 4$, show that $X_{a, 2}$ is rational over $\CC$.
\endroster

\medskip

{\bf{Exercise 9.}}
\roster
\item
Prove that the variety of $m \times n$ matrices of rank at most $t$ is
rational over any field. Find its non-smooth   locus.
\item 
Consider an $n \times n$ array of general linear forms on $\PP^n$.
Prove that the 
hypersurface defined by the determinant of this array is a rational variety.
When is this variety smooth?
\item Prove that every smooth cubic surface over an algebraically closed
field  is determinantal.
\endroster

\medskip

{\bf{Exercise 10.}}
Let $X$ be a smooth hypersurface in $\PP^n$ of degree $d \leq n $. 
Assuming the field is algebraically closed,
find a rational curve passing through every point of $X$.

\medskip
It is an open question whether or not such a smooth hypersurface has a 
rational
 surface through every point.

\setcounter {Ch}{1}
\heading  \Ch. Second Lecture  
\endheading
\label{lec2}
\setcounter{Sec}{0}

In  lecture one,  we saw that
smooth 
 hypersurfaces of large degree relative to their dimension are never
rational.
 On the other extreme, a linear hypersurface is 
obviously rational, and we observed that a quadric hypersurface is
rational over $k$, provided it has a  $k$-point. For cubics, however, 
we saw both rational and non-rational  examples, indicating that the 
rationality question for cubics is more subtle.

The purpose of this lecture is to consider the rationality question for
smooth cubic surfaces in detail. 
Because every smooth cubic surface over an algebraically closed field
is isomorphic to a blowup of the 
plane at six points,
{\footnote{For this,
 and other basic properties of cubic surfaces, the reader is
 referred to \cite{H, V 4},  or the more  elementary
account in \cite{Ger}. The discussion in \cite{R} is also quite fun
 and informative, although lacking in proofs.}
 every such surface is rational.
The interesting issue  is 
 the rationality or non-rationality of cubic surfaces defined 
over non-algebraically closed fields.

\bigskip
Rationality for cubic surfaces has
interesting applications to  Diophantine equations.  
For example, consider a cubic  $f(x, y, z) = 0$ 
defined over $\ZZ$
 or some other number ring. Let $S/\QQ$ be the 
 cubic surface it defines.
If $S$ is unirational, 
then the rational   map $ \A^2  \overset{\phi}\to\rmap S$ can be used to find
rational, and hence integer, solutions of the  equation:
 for each  
$(s, t) \in \A^2(\QQ)$,
the image 
$$
\phi(s, t) = (\phi_1(s,t), \phi_2(s, t), \phi_3(s, t))
$$
is a $\QQ$-solution.
 If  $S/\QQ$ is actually rational, then
$\phi$ is invertible on an open subset and we have an essentially complete
parameterization of the solutions: 
except on a 
locus of finitely many curves, every  solution of $f_3 =0$ is 
 uniquely  described 
 by the parameterization $(\phi_1(s, t), \phi_2(s, t), \phi_3(s,t))$.
{\footnote{Do not be misled to believe that if a map $\phi$  
defined over $\QQ$ is dominant, 
essentially all $\QQ$-solutions are parameterized by $\phi$. Indeed,
it is possible to miss most of them. 
For example, the squaring map
 $\A^1_{\QQ} \rmap \A^1_{\QQ}$ sending $x$ to $x^2$  is 
 dominant, but  its image, while Zariski dense, is 
fairly sparse, consisting  only of the perfect squares.}}
Segre's  1943 paper \cite{S43} contains  applications of this
 idea to problems posed by Mordell regarding
the  representation of  rational numbers by ternary cubic
forms. A nice survey of recent  developments in this direction is offered
by Colliot-Th\'el\`ene  in \cite{C86}.

\bigskip
{\bf{Exercise 11.}} (Swinnerton-Dyer \cite{S-D})
Consider the cubic surface defined by 
$$
t(x^2 + y^2) = (4z - 7t)(z^2 - 2t^2).
$$
 Prove that
\roster
\item The real points of this surface, considered as a real two-manifold,
 consist of 
two connected components.
\item  On one manifold component, $\QQ$-points are dense.
\item On the other manifold component, there are no $\QQ$-points.
\endroster
Hint:  Let $a$ 
 and $b$ be  integers,  and let $\prod p_i^{m_i}$ be a prime factorization of 
$a^2 + b^2$.  Then if $p_i = 3 {\text{ modulo }} 4$, then $m_i$ is even.

\bigskip
Rationality for cubic surfaces is quite subtle. 
Our  goal is  to clarify the 
situation by proving the  following theorem of Beniamino Segre.

\proclaim{\Sec. Segre's Theorem} \label{segre} 
If the Picard number of a smooth cubic surface is one, then the
surface is not rational.
\endproclaim

By definition, 
 the Picard number $\rho_k$ of a normal projective variety
over $k$ is the rank of its
 N\'eron-Severi group, the group of divisors up to numerical 
equivalence. Of course, this depends on the ground field.
 The Picard number of a smooth cubic surface $S$ over
an algebraically closed field  is seven: thinking of $S$ as the
blowup of $\PP^2$ at six points, the Picard group (and hence the
 N\'eron-Severi group) is freely
generated by the six exceptional lines
and the pull back of the hyperplane class.
On the other hand, 
the cubic surface $S$ may be defined over some non-algebraically closed field
 $k$,  even if the individual points we blow up are not defined over $k$.
In this case, the Picard number of $S/k$ may be less than seven.

Both the hypothesis that 
$X$ is smooth and the hypothesis that its Picard number is one 
are essential in Segre's theorem.  See Remark \refto{essential}.

\bigskip

Before proving Theorem \refto{segre}, we establish a 
criterion for detecting when a cubic surface
 has Picard number one. Recall that 
every cubic surface  over an algebraically
closed field contains exactly twenty-seven distinct 
lines.

\proclaim{\Sec. Theorem} \label{gal}
Let $S/k$ be 
a smooth cubic surface  in $\PP^3$ and consider the
action of the Galois group of $\bar k /k$ on the twenty seven 
lines of $S/\bar k$.  
 The following are equivalent.
\roster
\item
The Picard number  $\rho_k(S)$ is  one.
\item
The sum of the lines in each  orbit is linearly equivalent to 
a multiple of the hyperplane class on $S$.
\item
No  orbit consists of disjoint lines on $S$.
\endroster
\endproclaim

\subsubhead{\Thm}\endsubsubhead\label{orbits}
The proof of Theorem \refto{gal} makes use of the following general 
principle. If $k \subset L$ is a Galois extension of fields,
and $X$ is quasi-projective variety defined  over $k$, 
then  the Galois group  $G$ of $L/k$ acts on the $L$-points of $X$. 
An $L$-point of $X$ is a $k$-point if and only if it is fixed by this
action of $G$. Likewise,  the group $G$ permutes around the 
subschemes of $X$ defined over $L$, and such a subscheme is defined over
$k$ if and only if it is fixed by this action. Indeed, the subschemes of
$X$ defined over $k$ can be interpreted as orbits of the action of $G$ on 
the $L$-subschemes of  $X$. All these facts follow easily from 
 the simple observation that $G$ acts on 
a polynomial ring $L[X_1, \dots, X_n]$ in such a way that the
fixed subring is the polynomial ring $k[X_1, \dots, X_n]$, and the induced
map of schemes $$
\spec {L[X_1, \dots, X_n]} @>>> \spec {k[X_1, \dots, X_n]}
$$
is a quotient map: 
lying over each prime ideal in  $k[X_1, \dots, X_n]$ are 
prime ideals of $L[X_1, \dots, X_N]$ making up an orbit under the action 
of $G$. 
See  \cite{Na, 10.3}, or \cite{Se; p 108}.

\demo{Proof of Theorem \refto{gal}}
The twenty seven lines $\{L_i\}_{i=1,\dots, 27}$
on $S/\bar k$ span the N\'eron-Severi group of $S/\bar k$.
In fact, these lines generate the ``cone of curves'' for $S/\bar k$,
meaning that every effective curve on $S/\bar k$ 
is numerically equivalent to 
a non-negative integer combination the $L_i$.
See,  for example, 
\cite{H p405}.

Now, 
without  loss of generality, the ground field $k$  may be assumed  perfect.
Indeed, if  $k^{\infty}$ denotes the perfect closure
 $\cup_{e} k^{1/p^e}$ of $k$, then the 
  automorphism groups 
${\text{Aut }}\bar k/k$ and 
${\text{Aut }}\bar k/k^{\infty}$ are  identical, 
because $k^{\infty}$ is the precisely the fixed  
field of $G = {\text{Aut }}\bar k/k$. 
The cubic surface  $S$ is smooth
whether regarded over $k$ or $k^{\infty}$, and 
the orbits of the action of $G$ on the twenty-seven lines are independent
of this choice. Furthermore, any  divisor  on $S$ 
defined over  $k$  is {\it a priori\/}  defined over $k^{\infty}$, while
any divisor   defined over $k^{\infty}$ is defined over some purely inseparable
extension  $k^{1/p^e}$ of $k$, so that some ${p^e}^{th}$ multiple 
 is defined over $k$. 
This implies that the N\'eron Severi group has
 the same
rank over $k$ or $k^{\infty}$.  In particular, conditions (1), (2), and (3)
are equivalent over $k$ if and only if they are equivalent over $k^{\infty}$.

Assuming $k$ is perfect, we can
 understand $S$ over the non-algebraically closed field $k$ by 
considering the action of the Galois group $G$ of $\bar k/k$
 on the $\bar k$ points of $S$. 
The extension $k \hookrightarrow \bar k$ is Galois, so 
 divisors of $S$ defined
over $k$ can be identified with the $G$-invariant divisors of $S/\bar k$,
which in turn can be identified with the 
$G$-orbits of divisors defined over $\bar k$. 
In particular, the Picard group of $S/k$ is the  subgroup 
of the  Picard group of $S/{\bar k}$ invariant under the 
natural action of $G$. 
 If $\{L_{i_1}, \dots, L_{i_t}\}$ is an
orbit of $G$  on the twenty-seven lines, then $L_{i_1} + \dots + L_{i_t}$ is
an effective curve on $S/k$. We leave it as an exercise 
to check that these orbit sums span the Picard group 
(and hence N\'eron-Severi group) of $S/k$, and in fact, they 
 generate the cone of curves for $S/k$.

Thus,
the Picard
number is  one if and only if these orbit sums
are all multiples of each other. Equivalently, 
$\rho_k = 1$ 
if and only if every orbit sum is
 a  multiple of the
 hyperplane class. This establishes the equivalence of
(1) and (2).

To see that (1) implies (3), suppose that $\{L_1, \dots, L_t\}$ is
an orbit consisting of non-intersecting lines.
If the 
Picard
number is one, then any other orbit sum $L_{t+1} + \dots + L_{r}$
is a positive $\QQ$-multiple of $L_1 + \dots + L_t$. 
Then 
$$
(L_1 + \dots + L_t)^2 = \frac{m}{n}(L_1 + 
\dots  + L_t)\cdot(L_{t+1} + \dots + L_{r})
\,\,\,\,\,\,\,{\text{   with }} \frac{m}{n} \geq 0.
$$
This produces a contradiction: the left-hand-side is
 $\sum_{i=1}^t L_i^2 = -t$,
while the right-hand-side is non-negative because only distinct effective
divisors are intersected. This contradiction proves that (1) implies (3).

Finally, to prove (3) implies (1), we invoke the following easy
fact whose proof is left as an exercise:
 An effective curve with positive self-intersection lies in the interior
of the cone of curves.
Assuming that the Picard number of $S/k$ is at least two, 
the cone of curves has dimension at least two, and therefore has a
non-zero non-interior point. This forces 
one of the orbit sum generators
for the cone of curves, 
say $ L_1 + \dots  + L_t$, to have non-positive self-intersection.

Consider an 
 orbit $\{L_1, \dots, L_t\}$, and 
assume that $L_1$ intersects $L_2$. Note that $L_1$ can intersect
no other line $L_i$ in this orbit, because 
$$
\align
(L_1& + \dots + L_t)^2 = \sum_{i=1}^t L_i^2 + 2(\sum_{i<j} L_i \cdot L_j)\\ 
&= (-t) + 2({\text{number of pairs of intersecting lines in the orbit}})
\endalign
$$
is non-positive, and $G$ acts transitively on $\{L_1, \dots, L_t\}$.
Thus the orbit $\{L_1, \dots, L_t\}$ is partitioned into pairs of intersecting 
lines, with all pairs disjoint from each other, and $G$ acts transitively
on the set of  these  pairs.

Now consider the degenerate conic $L_1 + L_2$ on $S/\bar k$.
The plane $H$ spanned by $L_1$ and $L_2$ intersects $S$ in a third
line $L$. 
Consider the linear system on $S/\bar k$ of plane sections 
containing $L$. Throwing away the fixed component $L$, this gives a
base-point-free pencil of conics on $S/\bar k$. Obviously $L_1 + L_2$ 
belongs to this linear system, but we claim that in fact, every 
pair of intersecting lines in the orbit 
$\{L_1, \dots, L_t\}$ determines a conic in this linear system.
To see this, note that the linear system determines a morphism 
$S/k \rightarrow \PP^1$ whose fibers are the conics of the pencil.
 For $i \neq 1, 2$, each $L_i$ is disjoint
from  $L_1 + L_2$, so it must be contained in a fiber of this morphism.
From this we conclude that for every 
pair  $\{L_i, L_j\}$  of intersecting
lines in the orbit
$\{L_1, \dots, L_t\}$,
we have a linear equivalence $L + L_1 + L_2 = L + L_i + L_j$. 
Because  $G$ permutes around the intersecting pairs  $\{L_i, L_j\}$,
the line $L$  must be fixed
by $G$. This is true even if there are no intersecting pairs in the orbit
besides $L_1$ and $L_2$, for then $G$ fixes $L_1+L_2$, so that $G$ fixes $L$.
In either case, $L$ makes up its own orbit, contrary to assumption (3). 
\enddemoo

\medskip
{\bf {Exercise 12.}}
Complete the proof of
 Theorem \refto{gal} by proving the following two lemmas.
\roster
\item
The sums of the lines in each 
 orbit of the action of Galois group of $\bar k/k$ 
 on the twenty seven lines
on a cubic surface $S$ over $k$ 
generate the cone of curves.
\item 
On  a non-singular projective surface, an effective curve  with 
positive self intersection lies in the interior of the cone of curves.

\endroster

\medskip
Specific non-rational cubic surfaces over $\QQ$ are constructed in the
next exercise.

\medskip
{\bf{Exercise 13.} }
\roster
\item
Find all lines on a smooth cubic surface given by an equation of the form
  $u^3 = f_3(x, y)$ in affine
 coordinates. In particular,  find the lines on the Fermat hypersurface
 $a_1 x_1^3 + a_2x_2^3 + a_3x_3^3  = a_0$.
 Do the same for $u^2 = f_3(x, y)$. 
\item
Show that if $a$ is a rational number that is not a perfect cube,
then the rational surface defined by $x_1^3 + x_2^3 + x_3^3 = a $ 
has Picard number one over $\QQ$. Conclude that such surfaces are not rational.
\endroster

In fact, Segre showed that a surface over $\QQ$ defined by 
the  equation 
$
a_0 x_0^3 + a_1 x_1^3 + a_2x_2^3 + a_3x_3^3 = 0
$
 has
Picard number one if and only if,  for all permutations  $\sigma$ of four
 letters,
the rational number
$$
\frac{a_{\sigma(0)}a_{\sigma(1)}}{a_{\sigma(2)}a_{\sigma(3)}}
$$
is not a cube \cite{S43}. 
The proof of  Exercise 13  (on page 50) 
easily generalizes to  yield this stronger 
result.

\bigskip
\subhead{\Sec. Maps to Projective Space}\label{maps}\endsubhead
To prove Segre's theorem, 
we  need a good understanding  of what is involved in
mapping a smooth surface to the plane, so we digress to 
 consider this general problem.

A  rational map  $S \overset{\phi_{\gamma}}\to\rmap 
 \PP^2$ is given by a two-dimensional fixed-component-free linear
system of curves on a surface $S$,  say $\gamma$,
 contained in some complete linear system $|C|$.
If the map  were everywhere defined, we could compute the self 
intersection multiplicity of $\gamma$, by which we mean 
the intersection multiplicity of two general members of
$\gamma$, from 
the self intersection of its image. In particular, when 
 $\phi_{\gamma}$ is a {\it morphism\/}, 
 $\gamma^2 $ is its degree, assuming the map is finite.

However,  the map may not be
 defined everywhere.
The points where it is not defined
are precisely the base points of $\gamma$, call them
$P_1, \cdots, P_r$  and let their multiplicities  be $m_1, \cdots, m_r$.
 In general, the expected contribution of  the 
base point  $P_i$ to the  self intersection multiplicity of $\gamma$
 is
$m_i^2$. However, this  is valid only when two general
curves in  $\gamma$ have distinct tangents at $P_i$.
The number will be even higher if the  curves 
share tangents at  $P_i$: this is the case where
$\gamma$ has  
 ``infinitely near" base points.

To make this  precise, let $S' @>{\pi}>> S$ be the blowup of $S$ at a base 
point $P$ of  $\gamma$.
The birational transform of $\gamma$ is the linear system $\gamma'$ on $S'$
whose generic member is the birational transform of the
generic member of $\gamma$.
The base points of $\gamma'$ that lie in the exceptional fiber
 are called {\it base points of
$\gamma$ infinitely near $P$.\/}
They represent the tangent directions at $P$ that are shared by all members
of $\gamma$.

Let $C$ be a general member of $\gamma$ and let $C'$ be its birational
transform on the blowup $S'$. It is easy to verify that 
$
\pi^* C = mE +  C'
$
where $m$ is the multiplicity of $\gamma$ at the base point $P$ and $E$
is the exceptional fiber of the blowup $\pi$.
 The linear system $\gamma'$  on $S'$
is obtained by pulling back the linear system $\gamma$ and throwing away the
fixed component $mE$.
It particular, to the extent to which $S$ and $S'$ 
are ``the same,''  $\gamma$ and $\gamma'$ determine  ``the same''
 map to $\PP^2$: 
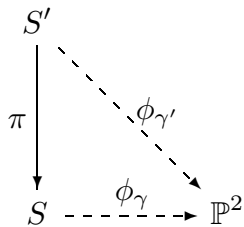
\newarrow{Nothing}{}{}{}{}{}
\diagram
 S' &   &    \\
\dTo^{\pi}   & \rdDashto[PS] \rdNothing^{\phi_{\gamma'}}  & \\
S      & \rDashto^{\phi_\gamma} & \PP^2 \\
\enddiagram

\smallskip
{ \bf{Exercise 14.} } 
With notation as above,  verify that    
$$
{C'}^2 = C^2 - m^2 \,\,\,\,
{\text{    and     }} \,\,\,\,\, C'\cdot K_{S'} = C\cdot  K_S + m,
$$
 where $C$ is a general member of $\gamma$ and $C'$ is 
a general member of $\gamma'$. We write this also as
 ${\gamma'}^2 = \gamma^2 - m^2 $ 
and  $\gamma'\cdot K_{S'} =  \gamma\cdot  K_S + m$.
\medskip

The process of blowing up base points
can be iterated until the multiplicities of all base points drop to zero.
In this way, we arrive at a smooth surface $\bar S$, and
a  base point free linear system $\bar \gamma$ 
defining the  ``same map" 
 (ie,  composed with the blowing up map) to  projective space.
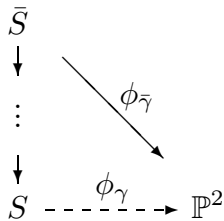
\diagram[height=2ex]
\bar S &    &   \\
\dTo   & \rdTo(2,6)^{\phi_{\bar\gamma}} &  \\
 {}    &    &   \\
\dDots &    &   \\
 {}    &    &   \\
\dTo   &    &   \\
S      & \rDashto^{\phi_\gamma} & \PP^2 \\
\enddiagram
This process is called
{\it resolving the 0\/} of the rational map $\phi_{\gamma}$.

\newarrow{Dashonto}{}{dash}{}{dash}{>>}

Beginning with a {\it birational\/} map 
$S \overset{\phi_{\gamma}}\to\rmap \PP^2,$
the 
linear system $\bar \gamma$ on $\bar S$ 
defines a {\it  birational  morphism\/} to 
$\PP^2$. In this case,  $\bar C^2 = 1$ and
$\bar C \cdot K_{\bar S} = -3$, where $\bar C$ is a general member of
$\bar\gamma$. Repeated applications of  Exercise 14 
express this in terms of 
the divisor $C$ on $S$:
$$
C^2 - \sum m_i^2 = 1 \,\,\,{\text{ and    }}  \,\,\,K_S \cdot C + \sum m_i = -3
\tag{\Thm.}\label{num}
$$
where the sums are taken over all base points, including the infinitely
near ones, and the $m_i$
are their multiplicities.

This leads to the following theorem.
\proclaim {\Sec. Theorem}\label{cubicsurface}
Let  $S$ be a smooth projective surface over $k$.
Then $S$ is rational over $k$ 
if and only if 
  $S$ admits a fixed-component-free
 two-dimensional linear system $\gamma$ defined over $k$ 
satisfying
 $$\gamma^2 - \sum m^2_i = 1$$ and 
$$K_S \cdot \gamma + \sum m_i = -3,$$ where the $m_i$ are the multiplicities
of all base points of $\gamma$, including the infinitely near ones.
\endproclaim

\demo{Proof}
Assume that $S$ is rational over $k$, and let 
$S \overset{\phi}\to\rmap \PP^2$ be a birational map.
 Let $\gamma$ be the fixed-component-free linear system obtained by 
pulling back the linear system of hyperplanes on $\PP^2$.
The dimension of $\gamma$ is two and the desired 
numerical conditions have been computed
already  in \refto{num}.

\medskip
Conversely, given a linear system $\gamma$ satisfying the given 
numerical conditions, 
it determines a rational map $S \overset{\phi_{\gamma}}\to\rmap \PP^2$
 defined over $k$.
We need only verify that this map is actually birational.
Because the map is {\it a priori\/} defined over $k$,
to check that it is birational  we are free to assume 
that $k$ is algebraically closed since 
whether or not the map  is dominant and degree one 
is unaffected by replacing $k$ by its algebraic closure.

Blow up the base
 points of $\gamma$,
including the infinitely near ones, 
to obtain  a  morphism $\bar S \overset{\bar\phi}\to\rightarrow  \PP^2$
resolving the indeterminacies of $\phi$.
The dimension of  $\bar \gamma$  is 2, and the numerical conditions 
$\bar {\gamma}^2 = 1$
and $\bar \gamma \cdot K_{\bar S} = -3$ hold. Because $S$ and $\bar S$ are 
birationally equivalent, it is sufficient to show that the morphism 
$\bar S \overset{\phi_{\bar{\gamma}}}\to\rightarrow \PP^2$ is a birational 
equivalence.

To  check that the map
 $\bar S \overset \phi_{{\bar {\gamma}}} \to \rightarrow \PP^2$ is
 surjective, assume, on the contrary,
that its image is a plane curve, $B$.
The fibers of $\bar S \rightarrow B$ would then be the elements  of
the linear system $\bar \gamma$.
 But thinking of the elements $\bar C \in \bar \gamma$
as the fibers of this map, we would have $\bar \gamma^2 = 0$.
 This contradicts the fact
that  $\bar\gamma^2 = 1$, so that  $\phi_{\bar{\gamma}}$ is  surjective.

Finally, the morphism determined by $\bar\gamma$ is generically
one-to-one because its degree is determined by the formula
$\bar\gamma^2 =  \,$(deg $\phi_{\bar\gamma}) H^2$.  Because
$\bar\gamma^2= 1$, we conclude that the map  must be one-to-one.
\enddemoo

\subhead{\Sec. Cautionary Example}\label{caution}\endsubhead
Let $F(X) \in \QQ[X]$.  As $\lambda $ and $\mu$ vary through $\CC$,
 the linear system  $|\lambda Y + \mu F(X)|$ is a one
dimensional linear system on $\A^2$ defined over $\QQ$.
The zeros of $F(X)$ determine the base points, since $(x, y)$ is a base point 
if and only if  $(x, y) = (0, \alpha) $ where $\alpha $ is a root of $F$.
These base points may not be defined over $\QQ$, although the
linear system is defined over $\QQ$. The 
map to projective space determined by this linear system
is
defined over $\QQ$ even when its base points are not.

\setcounter {Ch}{2}
\heading \Ch. \ Third Lecture 
\endheading
\label{lec3}
\setcounter {Sec}{0}
        
In this lecture we will prove Segre's Theorem. 
Essentially the same argument, with minor modifications to be made 
afterwards,   will prove the following  stronger 
theorem of Manin \cite{M66}.

\proclaim{\Sec. Theorem}\label{manin}
 Two smooth cubic surfaces defined over a perfect field,
 each of Picard number one, 
are birationally equivalent if and only if they are
projectively equivalent.
\endproclaim

{\it Caution:\/} Manin's theorem does not assert that every birational
equivalence is a projective equivalence. It guarantees 
 only that  if two surfaces
are birationally equivalent, then there exists a projective equivalence 
between them.

\bigskip 

\subsubhead{\Sec. Remark}\label{essential}\endsubsubhead
Manin's theorem 
holds under weaker hypotheses: the ground field need not be perfect and
smoothness can be weakened to non-singularity.
Because the proof in this case requires some distracting technicalities, 
we do not bother with it here; see \cite{K97, 5.3.3}.

However, the hypothesis of smoothness can not be weakened to include
singular varieties. For instance, consider a plane
conic defined over $k$,
 together with  six  points on it
 conjugate, but not individually 
defined, over $k$.
 By blowing up the six points and then 
contracting the conic, we achieve a singular
 cubic surface with  Picard number one.
All such surfaces are birationally equivalent to each other, but
two such are projectively equivalent if and only if the corresponding 
six-tuples of  points are
projectively equivalent in $\PP^2$.

Similarly, the hypothesis that the  Picard number is one is essential. 
For example,
there exist 
smooth cubic surfaces  of Picard number two containing 
exactly one line defined over $\QQ$. Because the line contains plenty
of $\QQ$-points, the surface is  unirational  by \refto{uni}, 
so it
contains a Zariski dense sets of $\QQ$-points.
Contract the $\QQ$-rational line and then blow up a general
 $\QQ$-point. The resulting
surface is a 
smooth cubic surface of Picard number two,
birationally equivalent over $\QQ$ to the original surface.
 However,  two such surfaces are not 
isomorphic in general. 
 Any isomorphism between them would amount to an automorphism of
the original cubic surface  interchanging the two distinguished points, but
because  the automorphism group of the cubic surface is finite,
this is not possible in general.

\bigskip
The proof of Segre's (and Manin's) theorem  begins 
with the  general observation 
that the Picard group of a smooth cubic surface 
 of Picard  number one is generated 
by the class of a hyperplane section. 
Indeed, 
the 
Picard group of $S$ is torsion-free,
 because 
Pic$(S) \subset $ Pic $(S_{\bar k}) \cong \ZZ^{\oplus 7}$.
Also, the 
 hyperplane class
 $H$ is not divisible, for otherwise  $H = mD$ for some $D$,
and  this would force $H^2 = 3 = m^2D^2$, which forces  $m = 1$.

Segre's theorem 
asserts that no cubic surface 
with Picard number one can be rational.
If this were false, there would be a 
 birational map $S \overset\phi\to\rmap \PP^2_k$ defined
over $k$.
The pull back of the hyperplane linear system on $\PP^2$
 is base point free on  the  dense open set of $S$ where $\phi$ is defined, and
the closure of this linear system in  
$S$ is the linear system $\gamma$
 defining the rational map $\phi$.
Because the Picard group is generated by the  hyperplane class $H$, 
the linear system  $\gamma$ must be  contained in the complete linear system
$|dH| $ for some $d$.
Therefore, the proof of  Segre's theorem will be complete upon 
proving  the following
theorem.

\proclaim{\Sec. Theorem} \label{inter}
If $S \subset \PP^3_k$ is a smooth cubic surface, then there is no 
fixed-component-free linear system
 on $S$ 
contained in 
$ |dH|$ that defines a birational map
 to the projective plane.
\endproclaim

Although the statement of this theorem is less appealing than Segre's Theorem,
we have, in effect, reduced the proof of Segre's theorem
 to the case where the
ground field is algebraically closed: if such a  linear system is defined 
over $k$, then it is defined over the algebraic closure of $k$.
Note that a na\"\i ve reduction to the algebraically closed field
case is not possible, 
as the Picard number is never  one over an 
algebraically closed field.

\demo{Proof of \refto{inter}}
Suppose that such a linear system, $\gamma$, exists and defines the birational
map 
$S \overset\phi_{\gamma}\to \rmap \PP^2_k$.
Without loss of generality, we assume $k$ is algebraically closed,
as explained above.

Let $P_1, \cdots, P_r$ be the base points of $\gamma$, including the infinitely
near ones, and let $m_1, \cdots, m_r$ be their multiplicities. From the
computation \refto{num},
 we have
$$
\align
&\sum m_i^2 = \gamma^2 - 1 = 3d^2 -1  \\
&\sum m_i=  -K_S\cdot \gamma - 3 = 3d-3. \tag{\Thm}\label{useful}
\endalign
$$
If all $m_i$  are less than or equal to $d$, then 
$$
3d^2 -1 = \sum m_i^2 \leq d\sum m_i = 3d^2 - 3d < 3d^2 - 1.
$$ 
This contradiction ensures that at least one $m_i$  is 
greater than $d$. 

Let $P$ be a base point of $\gamma$ of multiplicity
greater than $d$. There is no loss of generality in assuming that
$P \in S$, that is,
that $P$ is an actual base point, not 
an infinitely near one.
  This is because  the multiplicity of a  base point 
 is greater than or equal to the multiplicity of any base point
 infinitely near it. Indeed, the multiplicity  of $P$ is 
at least the sum of the 
multiplicities of all the base points infinitely near to $P$ to first order.
(We leave as an exercise the following fact:
 the multiplicity of $P \in S$ as a base point of the linear system $\gamma$
 is greater than or equal to the sum of the multiplicities of all base points
of $\gamma'$ which lie over $P$, where $\gamma'$ is the birational transform
of the linear system $\gamma$ under the blowing up map of 
$S$ at $P$.)

Furthermore,  the high multiplicity base point $P$ 
can not lie on any line
on $S$. Indeed, since $\gamma \subset |dH|$, we must have that 
$L\cdot C \leq d$ for all lines $L$ on $S$ and all $C \in \gamma$.
Computing $C\cdot L$ as the sum over all points 
 (with  multiplicities)  in $C \cap L$,
we see that $C$ can not have a multiple point of order more than $d$ on $L$.

\medskip
The proof of Theorem \refto{inter} proceeds by induction on 
 $d$. The inductive step is
 accomplished by  finding a birational self map of $S$ that takes
$\gamma$ to a linear system contained in the linear system $|d'H|$
with $d' < d$.

\medskip

How can we find a birational self-map of the cubic surface?
First, recall the following involution of 
a plane cubic curve $E$: fixing a point $P$
on $E$, define the map $\tau$  which  sends $Q \in E$ to the
third point of intersection of $E$ with $\overline{PQ}$. 
The map $\tau$ 
extends to an involution defined everywhere on $E$ by 
sending
the point $P$ to the intersection of $E$ with the tangent line through 
$P$.

We attempt to construct a similar involution of the cubic surface
$S \subset \PP^3$.
Define a self-map $\tau$ of $S$ as follows: fix a point $P $ on $S$ and for
each 
$Q \in S$, let $\tau(Q)$ be the third point of intersection of
$S$ with the line $\overline{PQ}$.
This defines a rational map
$$
S \overset\tau\to\rmap S
$$
such that $\tau^2 = id$. 
If we assume that $S$ contains no lines through $P$,
 then $\tau$ is defined  everywhere on $S$, except at $P$.
However, unlike the situation of the plane cubic, there is a whole
plane of tangent lines to  $S$ at $P$, so 
there is no way to extend $\tau$ to a morphism  at $P$.
 Indeed, $\tau$ 
contracts the  entire curve $D =  T_PS \cap S$ to the point $P$ on $S$.

As usual, the best way to sort out  different tangent 
directions at a point is to blow up.
Let $S' \overset\pi\to\rightarrow S$ be the blowup of $S$ at $P$, let $E$
be the exceptional fiber, and let  $D'$ be the birational transform  of the
curve $D =  T_PS \cap S.$ 

By definition, the blowup of $\PP^3$ at $P$ 
consists of those points 
$\{(x, L)\}$ in $ \PP^3 \times \PP^2$, where
$\PP^2 = \PP(T_P(\PP^3))$, 
such that $x \in L$.
The blowup of $S$ at $P$ is identified with  the birational transform of $S$.
 The blowing up map $\pi$ is the projection of $S'$ onto the first
 factor $S \subset \PP^3$. Let $q$ denote the projection of $S'$ onto
 the second factor $\PP^2$:

\newarrow{Into} C--->
\diagram
S' &   \rInto & S \times \PP^2  \\
\dTo^{\pi}   & \rdTo(2,2)^{q} &  \\
S      & \rDashto & \PP^2 \\
\enddiagram

\medskip
{\bf Exercise 15.} Show that $q$ is two-to-one and ramified along a
smooth  curve of degree 4. Find the equation of this branch locus.

\medskip
The rational map $\tau$ of $S$ extends to a  morphism $\tilde \tau$
of $S'$.
The following facts about $\tilde \tau$ are easily verified:
\roster
\item 
$\tilde\tau $ is the unique non-trivial Galois automorphism of the 
degree two cover
$S'$ of $\PP^2$.
\item $\tilde\tau(E) = D'$ and $\tilde\tau(D') = E.$
\item  $|\pi^*H - E | = |q^*L|$, where $L$ is a line in $\PP^2$;
that is $\pi^*\Cal O_S(1)(-E)  \cong q^*\Cal O_{\PP^2}(1)$.
\endroster

\medskip
To complete the proof of Theorem \refto{inter}, let $P$ 
be a base point of the linear system $\gamma$, and suppose that the 
multiplicity  of $P$ is $m > d$. 
 Let $\gamma'  = \pi^*\gamma - mE$ be the birational transform of
$\gamma$ via the blowup $S' @>{\pi}>> S$
 at $P$,  as discussed in \refto{maps}.
 Because
$\gamma \subset |dH|$, we have
$$
\gamma' + (m-d)E  = 
\pi^* \gamma  - dE  \subset  |\pi^*(dH) - dE| = |d(\pi^*H - E)|
=|q^* (dL)|.
 $$ 
 Applying the automorphism $\tilde\tau$ to $S'$, the elements of 
 $\gamma' + (m-d)E$ are  taken to another linear system inside 
$|q^*(dL)|$
  because  
$\tilde\tau \in Gal(S'/\PP^2)$ preserves any linear system pulled back from
 $\PP^2$.
 Therefore
 $$
\tilde \tau(\gamma' + (m-d)E) = 
\tilde\tau(\gamma') + (m-d)D' \subset |q^*(dL)| = |d(\pi^* H - E)|
 \subset |\pi^*(dH)|.
 $$
 Pushing back down to $S$, we have
 $$
 \tau(\gamma)  + (m-d)D   \subset |dH|.
 $$
 Because $D$ is a hyperplane section of $S$, we conclude that
 $$
 \tau(\gamma) \subset |(d - (m-d))H|.
 $$
 Because $m > d$, the linear system $\tau(\gamma)$ is contained in 
 $|d'H|$, with $d' < d$. 
By induction, we conclude eventually that $\gamma \subset |H|$,
that  is, that $d=1$.
 But now considering the two useful formulas in \refto{useful},
we see that $\sum m_i = 3d -3 = 0$ and $\sum m^2_i = 3d^2 -1 = 2$, where
the $m_i$ are the multiplicities of the base points of $\gamma$. 
The first equation forces all the $m_i$ to be zero, while the second 
forces two of the $m_i$ to be exactly one. This  contradiction
 completes the proof of Segre's theorem.
 \enddemoo
  
The proof Segre's theorem is easily altered to 
produce the following proof of Manin's theorem.
\demo{Proof of Theorem \refto{manin}}
Assume that $S \overset\phi\to\rmap S'$ is a birational equivalence.
Let $\gamma$ be the linear system on $S$ obtained by pulling back
the hyperplane system on $S'$ 
 via $\phi$.
Thus $\phi = \phi_{\gamma}$, and $\gamma$ is a linear system of dimension 3.

Because the Picard number of
 $S$ is one, we can assume, as before, that 
$\gamma \subset |dH|$ for some $d$, where $H$ is a hyperplane
section of $S$. Let $P_1, \dots, P_r$ be the 
base points of $\gamma$, including the infinitely near ones, and suppose their
multiplicities are $m_1, \dots, m_r$.
 Because $H^2 = 3 $ and 
$K_{S'} \cdot H = -3$, we again compute using Exercise 14
that 
$$
\sum m_i = 3d -3 \,\,\,{\text{   and    }} \,\, \sum m_i^2 = 3d^2 - 3.
$$

If all $m_i \leq d$, then  
 $3d^2 - 3 = \sum m^2_i \leq d(\sum m_i) = 3d^2 - 3d$.
This is possible only if $d = 1$, in which case $\phi$ is an 
automorphism of $\PP^3$.

Therefore, we may assume without loss of generality that 
some $m_i > d$, and the corresponding $P_i$ may be assumed to be on $S$.
The same  trick that was used to accomplish the
inductive step in the proof of Segre's theorem works here too.
The only problem is that the involution
$\tau$
will not be defined over $k$ unless the base point $P$ is defined over $k$.
A priori, the argument shows only that 
if $S$ and $S'$ are birationally equivalent over $k$, then
they are projectively equivalent over $\bar k$.

To see that $S$ and $S'$ are actually projectively equivalent over $k$, we
need to construct an involution $\tau$ defined over $k$. 
Because the Galois group 
of $\bar k$ over $k$ acts on the $P_i$ preserving multiplicities $m_i$,
 it follows from the fact that 
$\sum m_i = 3d - 3$ that 
at most two of the base points $P_i$ can have multiplicity greater than $d$.
If exactly one, say  $P_1$,
 has multiplicity greater than $d$, then the Galois group
fixes this base point. Because $k$ is perfect, this implies that $P_1$ is
defined over $k$, so the involution $\tau$ is defined over $k$ and the
proof of Manin's Theorem is complete. 
If exactly two base points, say $P_1$ and $P_2$, have multiplicity exceeding
$d$, then the Galois group must fix their union,  and so $P_1 \cup P_2$
is defined over $k$.  As before, $P_1$ may be assumed to be on $S$,
but it is possible that $P_2$ is infinitely near $P_1$.
The existence of an involution 
defined over $k$  needed to  complete the proof 
is left as an exercise.
\enddemoo

{\bf Exercise 16.}
Let $S$ be a cubic surface defined over $k$ and
let $P_1$ and $P_2$ be distinct 
$\bar k$-points of $S$ such that $P_1 \cup P_2$ is
defined over $k$.  Assuming that 
the line through $P_1$ and $P_2$ does not lie on $S$,
construct an involution $\tau$ of $S$ defined over $k$.
Similarly, interpret and prove a version of this statement 
 in the case where $P_2$ is an point on the blowup of
$S$ at $P_1$.

Use this involution to complete the proof of Manin's theorem.

\bigskip

{\bf 
{\Sec. Some Historical Remarks and Subsequent Developments.\/}}\label{history}
The geometry of cubic surfaces and the configuration of the twenty seven
lines on them occupied a tremendous amount of attention of the nineteenth
century algebraic geometers. 
However, many of the  beautiful geometric arguments presented here in our
discussion of rationality of  cubic
surfaces go back to Segre in the middle of this century,
 in the series of papers listed in the bibliography.
Segre was motivated, at least at first, by arithmetic questions of Mordell
on representing integers by ternary forms. The arithmetic applications
are still important today; see \cite{C86}, \cite{C92}, \cite{K96a}.

\medskip
Manin's theorem--- that any two birationally equivalent smooth cubic
surfaces of Picard number one
 are actually projectively equivalent---
 has a stronger  analog for three-folds, at least over an 
 algebraically
closed field of characteristic zero.
Indeed, in 1971,  Iskovskikh  and Manin
proved  that any  birational equivalence 
 between smooth hypersurfaces of degree
four in $\PP^4$ must be a projective equivalence
\cite{IM}.  This implies that no smooth  quartic threefold is 
rational: the birational automorphism group of the quartic threefold
is the same as its group of projective automorphisms; 
since the latter is finite, the threefold can not be birationally 
equivalent to $\PP^2$, which has an infinite automorphism group.
 Iskovskikh and Manin acknowledge their indebtedness to Fano,
who, despite some serious errors,  had laid the foundations
of this investigation, using ideas of Max Noether.

This theorem of  Iskovskikh and Manin 
 settled  the longstanding 
L\"uroth problem: 
\smallskip
{\centerline{\it Is every unirational variety actually rational?\/}}
\smallskip

For curves, L\"uroth's theorem states that every unirational curve is
rational; 
 Castlenuovo and Enriques later showed that every unirational surface is 
rational (with some help from 
 Zariski in  prime characteristic, where ``unirational'' should be 
read ``separably unirational''). 
Fano and others had long expected a negative answer to the higher
 dimensional L\"uroth problem, but the question 
 remained open for many years, despite some 
 erroneous counterexamples proposed by prominent mathematicians. 
The theorem of Iskovskikh and Manin on the non-existence of rational 
quartic three-folds settled--- negatively--- the L\"uroth problem, because
 undisputed examples of 
  unirational quartic three-folds had already been
 constructed by Segre \cite{S60}.

Around the same time, Clemens and Griffiths also resolved the L\"uroth problem,
by showing that 
there exist no
smooth rational {\it cubic\/} three-folds 
 \cite{CG}. (It is not hard to see  that 
{\it every\/} smooth cubic threefold is 
unirational, an elementary fact Clemens and Griffiths attribute to Max Noether;
see \cite{CG, p 352}.)
Clemens and Griffiths approach was 
entirely different, based on a study of the Intermediate Jacobian
of the cubic threefold.

Also in the early seventies, 
Artin and Mumford gave a simple proof  of
the existence of non-rational unirational varieties in all dimensions
three or higher,  using the observation
that the torsion subgroup of the third integral
singular cohomology group of a non-singular projective variety over
$\CC$ is a birational
invariant \cite{AM}.

\medskip

We now have a reasonably complete understanding of rationality 
for smooth three-folds; 
 see the papers of
Sarkisov, Iskovskikh, Bardelli 
 and Beauville  listed in the bibliography.  
Meanwhile, Colliot-Th\'el\`ene and Ojanguren 
developed  further  the
  examples of Artin and Mumford  in all dimensions \cite{CO}.

For four-folds, considerably less is known. 
Pukhlikov  generalized the methods of \cite{IM} to show that
 a birational equivalence between  two 
smooth quintics in $\PP^5$  is actually a projective
equivalence. Again the corollary follows:
 there 
exist no smooth rational quintic four-folds \cite{P}.
The short paper of Tregub \cite{Tr} gives some nice
 examples  of rational cubic hypersurfaces
in $\PP^5$. 
But there is not a single known
 example of a smooth rational quartic hypersurface of 
dimension four or higher; nor is there any proof that one can not exist.

\bigskip
Although we now know that not every 
 unirational  variety is  rational, it is natural to wonder about 
other classes of varieties 
that are ``close to rational.'' 
Recent work of Koll\'ar has established the
existence of abundant examples of non-rational varieties with 
various other nice properties, such as {\it rational connectivity\/}
 \cite{K95}, \cite{K97}. A smooth projective variety is rationally 
connected if every two points are joined by a rational curve.

We have seen that a  smooth hypersurface whose
 degree is no greater  than its embedding dimension is covered by 
rational curves. In fact, such a hypersurface
has the  even stronger property of rational connectivity.
This is a special case of a theorem of  Koll\'ar, Miyaoka and Mori 
 stating that  
every smooth Fano variety  of characteristic zero
is rationally connected \cite{KoMM}.
Koll\'ar's book \cite{K96} presents an in-depth study of rational curves
on algebraic varieties.

\bigskip

\setcounter {Ch}  {3}
\heading  \Ch.\ Fourth  Lecture 
\endheading
\label {lec4}

It is easy to find smooth projective varieties that are not rational:
a high degree hypersurface in $\PP^n$ is not even separably 
unirational, 
as we observed in 
the first lecture.
It is much harder to find examples of non-rational varieties 
when the obvious obstructions to rationality, such as the plurigenera,
 vanish.

\proclaim {\Sec. Definition} \label{fanodef}
A {\it Fano variety\/} is
a smooth projective variety 
whose anti-canonical bundle is ample.
\endproclaim 

By anti-canonical bundle we mean the 
 the dual of the sheaf
of differential $n$-forms,  where $n$ is the dimension.
More generally, it 
is not really necessary to assume the variety is smooth in the
definition of Fano, provided one can make sense of
the sheaf of differential $n$-forms 
as a line bundle.

The plurigenera 
vanish  for any Fano variety, so 
it is natural to wonder whether all Fano varieties might be rational.
In the final two lectures, we
 construct specific examples of non-rational Fano varieties.
In this lecture,  we focus mostly on the prime characteristic case.
We will construct hypersurfaces
with
ample anti-canonical divisor that are not even {\it ruled\/} in
 characteristic  $p$. 
The argument makes  
heavy use of the special nature of derivations in characteristic $p$, and
is not valid in characteristic zero. Unfortunately, our hypersurfaces are not
smooth.
In Lecture 5, we will show how to deduce the existence
of {\it smooth\/} projective
Fano varieties that are not ruled. In particular, we will
 produce explicit and abundant  examples
of non-rational smooth  Fano hypersurfaces.
This complements  the rather special examples previously known;
 see \refto{history}.

\bigskip
It will be necessary to consider geometric conditions weaker than rationality,
so we recall the following definitions.

\proclaim{\Sec. Definition}\label{ruled}
A  variety $X$ is  {\it ruled\/} if there exists a
variety $Y$ and a birational map 
$Y \times \PP^1 \rmap X$.

A   variety $X$ is  {\it uniruled\/} if there exists
a variety $Y$ and a generically finite dominant map
$Y \times \PP^1 \rmap X$. 
The variety $X$ is {\it separably uniruled\/}
 if this birational map is generically 
\'etale, that is, if the function field extension is separable. 
\endproclaim

Loosely speaking, a variety is {\it uniruled\/} if it is covered by 
 rational curves.
Of course, every rational variety $X$ is ruled, since 
$\PP^n$ is birationally equivalent to
$\PP^{n-1}\times \PP^1$. 
Furthermore, every ruled variety is separably uniruled, and
every unirational variety $X$ is uniruled. In characteristic zero, 
separably uniruled is equivalent to uniruled.

In this lecture, we prove the following theorem.

\proclaim{\Sec. Theorem}\label{charp}
For any prime $p > 0$,
there exists a projective variety  of characteristic $p$ 
that is not separably uniruled, but whose anti-canonical sheaf is
an ample invertible sheaf. In fact, 
 there are examples given by  hypersurfaces with local
equation of the form $y^p = f(x_1, \dots, x_n)$.
\endproclaim

We will prove even more: for every sufficiently general $f$ 
 whose  degree is
in a certain range (described in \refto{summary}), the
affine hypersurface defined by $y^p -f(x_1, \dots, x_n)$ has
a natural compactification that is  not separably uniruled 
but 
 whose anti-canonical 
sheaf is ample invertible.

It is relatively easy to find singular non-rational varieties
with ample anti-canonical sheaf over any field:
 the cone over an elliptic curve is an immediate
example.  It is harder to produce non-ruled examples. 
 The point here
 is not so much the mere existence of such varieties,
but  rather the explicit construction of them in 
\refto{construct}.
Essentially the same construction 
yields  {\it smooth\/} Fano varieties that are
not ruled  in characteristic zero, 
but we will need the characteristic $p$ results to deduce this.

\medskip
Non-ruledness 
will be proved by appealing to the next proposition.
Recall that a 
 line bundle $\MM$ on a variety 
 $X$ is big if some multiple defines a birational map 
$X \overset\phi_{|\MM^m|}\to\rmap X' \subset \PP^n$ 
 to its image in projective space.
Intuitively, big line bundles are ``birationally ample.''

\proclaim{\Sec. Proposition}\label{rulecrit}
If $X$ is a smooth projective variety admitting a big line bundle
$\MM \subset \bigwedge^i \Omega_X$, then $X$ is not
separably uniruled.
\endproclaim

\medskip
{\bf Remark.} 
	If $X$ is a smooth projective variety of characteristic zero,
 the Bogomolov-Sommese vanishing theorem  guarantees that 
the only sheaf of differential forms 
that can contain {\it any \/} big line bundle is the line bundle  
$\bigwedge^{\dim X}\Omega_X$ itself;  see \cite{EV, p 58}.
In particular, in characteristic zero, 
Proposition \refto{rulecrit} degenerates  to  the following statement:
no smooth projective  variety of general type is uniruled. 

\demo{Proof of \refto{rulecrit}}
Suppose that $X$ is separably uniruled, and let 
$Y \times \PP^1 \overset{\phi}\to\rmap X$ 
be a generically \'etale 
 map. Without loss of generality, $Y$ may be replaced by a smooth
affine
open subset on which $\phi$ is a morphism.
  The pull back map on differential forms
$$
\phi^*\Omega_X \overset{\phi^*}\to\rightarrow \Omega_{Y \times \PP^1}
$$  
is an isomorphism on a dense open set because $\phi$ is generically \'etale.
In fact, all we need is that $\phi^*$ is injective, 
so that  for the subsheaf $\Cal M$ of $\bigwedge^i\Omega_X$,
we have the inclusion
$\phi^* \MM \subset \bigwedge^i\Omega_{Y \times \PP^1}$.

Replacing $Y$ by an even smaller affine open subset if necessary,
we may assume its sheaf of differential forms is free.
In particular, 
$$
\Omega_{Y \times \PP^1} \cong \Omega_Y \oplus \Omega_{\PP^1} \cong
\Cal O_Y^{\oplus n-1} \oplus \Cal O_{\PP^1}(-2),
$$ where
we have abused
 notation slightly by omitting the symbols for pull back of $\Omega_Y$
and $\Omega_{\PP^1}$ to the product $Y \times \PP^1$.
We have
$$
\phi^*\Cal M \subset
 \bigwedge^i (\Omega_{Y \times \PP^1})
\cong
\bigwedge^i(\Cal O_Y^{\oplus n-1} \oplus \Cal O_{\PP^1}(-2)).
$$
If $\Cal M$ is big, so is its pull back under any generically finite map, so
in particular, the powers of 
$\phi^*\Cal M$ should have enough  global sections to separate points
on an open set of $Y \times \PP^1$.
 But this is impossible, since the  symmetric powers of
 $$
\bigwedge^i(\Cal O_Y^{\oplus n-1} \oplus \Cal O_{\PP^1}(-2))
$$
obviously have very few global sections. 
\enddemoo

\subsubhead
{\Thm. Remark} 
\endsubsubhead
The  inclusion 
$\phi^* \MM \subset \Omega_{Y \times \PP^1}$ always holds in
 characteristic zero. But in characteristic $p$ it can fail when the map 
is inseparable.
In fact, the varieties  to which we apply Proposition \refto{rulecrit}
to verify non-separably uniruledness {\it are\/} uniruled. See
also  Exercise 6.

\bigskip
\subhead{\Sec. The Construction}\label{construct}\endsubhead
Consider an affine  hypersurface in $\A^n$ defined by the equation 
$y^p = f_{mp}(x_1, \dots x_n)$,
 where $f_{mp}$ is a (non-homogeneous) polynomial
of degree $mp$. 

\medskip
{\bf{Exercise 17:}} 
Prove that 
the non-smooth points of 
the affine hypersurface  defined by 
$y^p - f$  are described  in terms of the critical points of $f$ as
follows.
 In the case where the ground field has
characteristic $p$,  
the non-smooth points  are in one-to-one correspondence with the
critical points of $f$, with 
$(y, x_1, \dots, x_n) = (\lambda,  \lambda_1, \dots, \lambda_n)$
a non-smooth point if and only if $(\lambda_1, \dots, \lambda_n)$
is a critical point  of $f$.
In the case where the ground field has characteristic not equal to $p$,
$(\lambda,  \lambda_1, \dots, \lambda_n)$
is a non-smooth point if and only if $(\lambda_1, \dots, \lambda_n)$
is a critical point of $f$ of critical value zero.
Show that for sufficiently general choice of $f$, the 
hypersurface has only isolated non-smooth points in
the characteristic $p$ case, 
and is everywhere smooth in characteristic zero case.
 
\medskip
The non-rational Fano variety will be a hypersurface defined by an equation
of the  form $y^p - f_{mp}$.
We need a compactification
{\footnote{More accurately, we need  a completion of our variety,
 not a compactification. The Zariski topology is always compact in any case. 
However, over the complex numbers, a variety is complete
 precisely when it is compact as a complex manifold.
For this reason, it is customary to refer to  ``compactification'' 
of a variety.}}
 of
this affine hypersurface. The obvious one, namely its closure in $\PP^{n+1}$, 
is insufficient for our purposes; the next exercise indicates why.

\smallskip
{\bf{ Exercise 18:}}
Prove  that the projective closure in $\PP^{n+1}$ of the affine hypersurface
defined by
$y^p - f_{mp}$ is never smooth in any 
 characteristic (whenever $m \geq 2$).
\smallskip

\subsubhead{\Thm}\label{proj}\endsubsubhead
Instead, we  compactify the hypersurface by taking its closure in a
weighted projective space. Let $F_{mp}(X_0, \dots, X_n)$ be a homogeneous
polynomial of degree $mp$ in the variables $X_0, \dots, X_n$.
Let $\bold P$ be the weighted projective space with coordinates 
$$Y, X_0, X_1, \dots, X_n,$$
where $Y$ is degree $m$ and the $X_i$ are all degree
one; that is $\bold P = \proj k[Y, X_0, \dots, X_n]$.
A common notation for this weighted projective space
is $\PP(m, 1, 1, \dots, 1)$, but we will shorten this to just
$\bold P$.

Consider the closed subscheme $Z$ in $\bold P$ defined by 
$Y^p - F_{mp}$. In the affine chart  of $\bold  P$ where $X_i$ does not
vanish, it is an affine hypersurface of the type considered above.
Because these charts cover $Z$, 
the projective scheme $Z$ is
 smooth 
for generic choice of $F_{mp}$ 
(at least when the characteristic is not $p$). 
Our goal is to show that we can choose the integers $m, n$ and $p$ such that
$Z$ is Fano but not separably uniruled.

\subsubhead{\Thm}\endsubsubhead
\label{alt}
The variety $Z$ has an 
alternate description as a cyclic cover of
projective space.
 Let $\PP = \PP^n$ be the projective space with homogeneous
coordinates $X_0: \dots : X_n$, and let $V_i$ be the open affine set
where $X_i \neq 0$. The affine coordinates on $V_i$ are $\frac{X_j}{X_i}$.

Consider the line bundle $\Cal O_{\PP}(m)$ on $\PP$. 
Fixing local generators $s_i$ on  $V_i$, we have  patching data for 
$\Cal O_{\PP}(m)$ on $V_i \cap V_j$
$$
s_i = \left(\frac{X_i}{X_j}\right)^ms_j.
$$
Let $U$ be the variety formed by the union of 
 the open sets $U_i = V_i \times \A^1$,
patched together by the relations $y_i = (\frac{X_i}{X_j})^{-m}y_j$,
where $y_i$ is the local coordinate for the copy of $\A^1$ in $U_i$.
The natural  projection 
$U \overset\pi\to\rightarrow \PP^n$ defines
an $\A^1$-bundle over $\PP$. This is the total space of the line bundle whose
sheaf of sections is $\Cal O_{\PP}(m)$; it is neither affine nor projective.

Now let $F_{mp}(X_0, \dots, X_N)$ be a homogeneous polynomial of degree $mp$,
 and let $Z$ be the subvariety of $U$ defined locally by the equations
$y_i^p - \frac{F_{mp}}{X_i^{mp}}$ in the open subset $U_i$.
In each patch, the variety $Z$  has exactly the form of the affine 
hypersurface $y^p = f_{mp}$.  Thus,
for a generic choice of the homogeneous polynomial $F_{mp}$, 
the variety $Z$ is smooth, at least assuming
$k$ is not of characteristic $p$.
This construction obviously produces 
a variety isomorphic to the $Z$
constructed  in \refto{proj} as a 
subscheme of the weighted projective space $\bold P$.

\subsubhead{\Thm}\label{-Z}\endsubsubhead 
There is a natural isomorphism
$\Cal O_U(-Z) = \pi^* \Cal O_{\PP}(-mp).$ \,
Indeed, the patching data for 
$\Cal O_U(-Z)$, the defining ideal for $Z$ as a closed subvariety
of $U$, has the same transition 
functions    as $\pi^*\Cal O_{\PP}(-mp)$: a local generator for either sheaf
on the affine neighborhood
 $U_j$ is transformed into a local generator on $U_i$
by multiplication by  
  $(\frac{X_i}{X_j})^{-mp}$.

\medskip

The map $Z @>{\pi_Z}>> \PP$, obtained by restricting the natural 
projection $U @>{\pi}>> \PP$,  is a finite  surjective 
map, of degree $p$. This is easily
checked locally: each point of $\PP$ has precisely $p$ preimages.
Of course, when $k$ has
characteristic $p$, the preimage of a point in $\PP$
 is a single 
point of multiplicity $p$, which is to say, $Z$ is purely inseparable over 
$\PP$.

\subsubhead{\Thm. The Fano range}\label{-F}\endsubsubhead 
For appropriate choices of the integers 
 $m, n$ and $p$, the variety $Z$ will be Fano. Specifically, 
 the anti-canonical sheaf of $Z$ is an ample
invertible sheaf
whenever $m, n$ and $p$ satisfy
$$ 
mp - m < n+1,
$$
 whether or not $Z$ is smooth.
There are several  different ways to see this.
We explain the least succinct way first, because the computation will
be of the most use later.

\medskip
{\it Method 1.}
Compute $K_Z$ using the adjunction formula for $Z \subset U$.
First compute $K_U$ using the
exact sequence 
$$
0 \rightarrow \pi^*\Omega_{\PP} \rightarrow \Omega_U \rightarrow
\pi^*\Cal O_{\PP}(-m) \rightarrow 0.
\tag{\Thm.}\label{seq1}
$$
The exactness of \refto{seq1} is easily verified
by the local  computation 
$$
dy_i = d\left[\left(\frac{X_i}{X_j}\right)^{-m}y_j\right] =
\left(\frac{X_i}{X_j}\right)^{-m}dy_j + 
d\left[\left(\frac{X_i}{X_j}\right)^{-m}\right]y_j,
$$
observing that $d\left[\left(\frac{X_i}{X_j}\right)^{-m}\right]$
 is pulled back from 
$\Omega_{\PP}$ and that  $(\frac{X_i}{X_j})^{-m}dy_j$ maps
to a local generator $\pi^*\Cal O_{\PP}(-m)$.
Therefore, 
$$
\omega_U = \bigwedge^{n+1} \Omega_U = 
\left(\bigwedge^n\pi^*\Omega_{\PP}\right) \otimes \pi^*\Cal O_{\PP}(-m) = 
 \pi^*\Cal O_{\PP}(-n-1-m).
$$
By adjunction, therefore, 
$$
\omega_Z = {(\omega_U \otimes \Cal O_U(Z))|}_{Z} = \pi_Z^*\Cal O_{\PP}(-n-1-m)
\otimes \pi_Z^*\Cal O_{\PP}(mp) = \pi_Z^*\Cal O_{\PP}(mp - n - 1- m),
$$
where $\pi_Z$ is the restriction of $\pi$ to $Z$.
Because $Z @>{\pi_Z}>> \PP$ is a finite map, the pull back of an ample
line bundle on $\PP$ is ample.
This proves that 
$\omega_Z^{-1}$ is ample whenever the numerical condition $mp - m < n+1$ 
is satisfied.

\smallskip
{\it Method 2.\/}
Compute 
$K_Z$ by the adjunction formula for  $Z$ in the
 weighted projective space $\bold P$.
This yields  $K_Z = (K_{\bold P} + Z)|_Z, $ so that 
$\omega_Z = \Cal O_{Z}(-m-n-1 + mp)$. This sheaf is invertible, because 
$Z$ is locally a hypersurface in a smooth variety,
 and hence Gorenstein. It follows that $\omega_Z^{-1}$ is ample 
if and only if $mp - m < n+1$.

\smallskip
{\it Method 3.\/} Compute $K_Z$ using the Hurwitz
formula for the finite map of $Z$ to $\PP^n$. This formula predicts that 
the canonical class for $Z$ is the pullback of the canonical class of $\PP^n$
plus the ramification divisor. This can not be applied 
when the characteristic is $p$, however, 
because $Z @>>> \PP$  is  inseparable, ie, the map is everywhere ramified.

\subsubhead{\Thm. A special subbundle of differential forms}\label{special}
\endsubsubhead 
Our goal
is to 
 apply Proposition \refto{ruled}
to conclude that the variety $Z$ is not ruled. In order to do so,  
we need to find a 
big line sub-bundle of a sheaf of differential forms on $Z$.
This will be an exterior power of a special subsheaf $Q$
of $\Omega_Z$, which we now construct.

Consider the familiar   exact sequence 
$$
{\Cal O_U(-Z)}|_Z \overset d\to\rightarrow
{\Omega_U|}_Z \rightarrow \Omega_Z \rightarrow 0.
\tag{\Thm}\label{seq3}
$$
(This is  the ``conormal" or the ``second exact sequence;"
see \cite{H, II 8.12}.)
Let us scrutinize the map $d$. 
In the chart $U_0 = \{X_0 \neq 0\}$, set $x_i = \frac {X_i}{X_0}$, 
$y = \frac{Y}{X_0},$  and
  $\frac{F_{mp}}{X_{0}^{mp}} = f(x_1, \dots, x_n)$.
The map 
$$
{\Cal O_U(-Z)}|_{Z} \overset d\to \rightarrow {\Omega_U}|_Z
$$
sends the local generator $y^p - f(x_1, \dots, x_n)$
to 
$$
d(y^p - f(x_1, \dots, x_n)) = -\frac{\partial f}{\partial x_1}dx_1 - \dots
-\frac{\partial f}{\partial x_n}dx_n  + p y^{p-1}dy.
$$
Something very interesting happens in characteristic $p$:
the image of $d$ is contained in the subsheaf of ${\Omega_U}|_Z$ generated 
by the differentials $dx_1, \dots, dx_n$, that is, 
$$
d(\Cal O_U(-Z)|_Z) \subset
 \pi_Z^*\Omega_{\PP}.
$$
Making use of the $\Cal O_Z$-module isomorphism
 $\Cal O_U(-Z)|_Z \cong \pi_Z\Cal O_{\PP}(-mp)$ derived in \refto{-Z},
we can define a
 $\Cal O_Z$-module map {\it in characteristic $p$ only\/}
$$
\pi_Z^*\Cal O_{\PP}(-mp) \overset d\to\rightarrow \pi^*_Z\Omega_{\PP}
\tag{\Thm}\label{d}
$$
sending a local generator $f$ to
 $df = \sum 
\frac{\partial f}{\partial x_i}dx_i$ and extending $\Cal O_Z$-linearly.
The use of the symbol $d$ to denote this map is somewhat misleading,
 since the map is not a derivation, but is $\Cal O_Z$-linear.
Miraculously, this  is a well defined  $\Cal O_Z$-module
map in characteristic $p$, because the transition functions for 
$\Cal O_{\PP}(-mp)$ are $p^{th}$ powers, and are therefore killed by $d$.

Let $Q$ be the cokernel of the $\Cal O_Z$-module 
 map \refto{d}. There is
an exact sequence
of
$\Cal O_Z$-modules
$$
\pi_Z^*\Omega_{\PP} \rightarrow {\Omega_{U}}|_Z \rightarrow
\pi_Z^*\Cal O_{\PP}(-m) \rightarrow 0.
\tag{\Thm.}
\label{seq2}
$$
obtained by 
restricting 
 sequence \refto{seq1}
 to $Z$. Combining this with the 
exact sequence \refto{seq3}, we 
get an exact sequence of $\Cal O_Z$-modules:
$$
0 \rightarrow Q \rightarrow \Omega_Z  \rightarrow \pi_Z^*\Cal O(-m) 
\rightarrow 0.
\tag{\Thm}\label{seq4}
$$
An exterior power of  $Q$ will  give us the desired big subbundle of a 
 sheaf of differential forms (at least after desingularizing). 
It is important to realize that the 
assumption that $k$ has characteristic $p$ is essential:
 nothing like this is possible in characteristic zero.

 The next exercise is not essential
for our computation, but it should help clarify what is
going on in the construction of $Q$.

\smallskip
{\bf Exercise 19:}
Let $X$ be an arbitrary variety over a field $k$.
A connection on an invertible sheaf  $\LL$ of $\Cal O_X$-modules 
is a  $k$-linear map
$$
\LL @>{\nabla}>> \LL \otimes \Omega_X
$$
satisfying $\nabla(fs) = f\nabla(s) + s \otimes df$ for local sections 
$s \in \LL$
 and $f \in \Cal O_X$.
\roster
\item 
Explain how a connection can be interpreted as a rule for differentiating
sections of line bundles.
\item Show that if $k$ has prime characteristic $p$, then 
any line bundle that is a $p^{th}$ power admits a  connection.
\item Observe  that the
 map  $d$ above in \refto{d} can be constructed from
the composition 
$$
\Cal O @>{{\text{mult by }s}}>> \LL^p @>{\nabla}>>
\LL^p \otimes \Omega_X
$$
where $s$ is a global section of $\LL^p$,
using the identifications $H^0(\Cal L^p \otimes \Omega_X) = 
{\text{Hom}}(\Cal O_X, 
\Cal L^p \otimes \Omega_X) = {\text{Hom}}(\Cal L^{-p}, \Omega_X).$ 
\endroster

\bigskip
\subsubhead{\Thm.  Bigness of the special subbundle of differential forms}
\endsubsubhead
With an eye towards applying Proposition \refto{rulecrit},
we hope to find integers $m, n$ and $p$ so that 
 $\bigwedge^{n-1}Q$ will be a big invertible sheaf.
Assume, for a moment,
 that
$Z$ is smooth.
The  sequence \refto{seq4} would imply that $Q$ is locally free
of rank $n-1$, and that 
$\bigwedge^{n-1}Q \hookrightarrow \bigwedge^{n-1}\Omega_Z$.
We could easily determine the range of values for $m, n$ and $p$,
for which
 $\bigwedge^{n-1}Q$  is big.
Indeed,
 $$
\omega_Z = \bigwedge^n\Omega_Z = 
\bigwedge^{n-1}Q\otimes \pi_Z^*\Cal O_{\PP}(-m),
$$
so that 
$$
 \bigwedge^{n-1}Q = \omega_Z\otimes \pi_Z^*\Cal O_{\PP}(m)
= 
\pi^*\Cal O_{\PP}(mp - n - 1),
$$
using the
isomorphism 
$\omega_Z = \pi^*\Cal O_{\PP}(mp-m-n-1)$  verified in \refto{-F}. 
From this we could conclude that $\bigwedge^{n-1}Q$ is ample (and hence big) 
whenever  $n+1 < mp$.

\smallskip
Because there 
 are plenty of choices of integers $m, p$ and $n$ for which the
constraints 
$$
mp -m < n+1 < mp  
$$
hold, we would expect to 
be able to
 use Proposition \refto{rulecrit} to 
find plenty of non-rational Fano varieties. 
Unfortunately, however, this 
argument  fails because 
the variety $Z$ is {\it virtually never\/} smooth:
this is the  price to be paid for the characteristic $p$ trickery that allowed
us to construct $Q$. Indeed, $\wedge^{n-1}Q$ is not even an invertible 
sheaf in general. 

\medskip
It is easy to alter  $\wedge^{n-1}Q$ so as to get  big invertible sheaf.
On the smooth locus of $Z$, 
the  sheaf $\wedge^{n-1}Q$ 
is naturally isomorphic to 
 $\pi^*\Cal O_{\PP}(mp - n- 1)$.
Now, because $Z$ is a hypersurface in a smooth variety, it is normal,
so that  any invertible sheaf 
defined on  the complement of a codimension two closed subscheme 
  extends {\it uniquely\/} to  a reflexive sheaf 
of $\Cal O_Z$-modules.
Since $\bigwedge^{n-1} Q$ agrees with $\pi_Z^*\Cal O_{\PP}(mp-n-1)$ on the
smooth locus, 
its ``reflexive hull'' 
 $$(\bigwedge^{n-1}Q)^{**} = 
\Cal{H}om_{\Cal O_Z}(\Cal{H}om_{\Cal O_Z}
(\bigwedge^{n-1}Q, \Cal O_Z), \Cal O_Z)
$$
 is an
 invertible sheaf of $\Cal O_Z$-modules  
isomorphic to $\pi_Z^*\Cal O_{\PP}(mp-n-1)$.
This sheaf  is  ample when 
$mp > n+1$, and is  a subsheaf of $(\bigwedge^{n-1}\Omega_Z)^{**}$. 
On  the smooth locus of $Z$, 
this  sheaf restricts to an invertible subsheaf of
 differential $n-1$ forms on $Z$.

\medskip
\subsubhead{\Thm. Desingularizing $Z$}\label{desing}\endsubsubhead
Proposition \refto{rulecrit} is only valid for smooth varieties,
so we 
must resolve  the singularities of $Z$.
Bigness is preserved under birational pull back, 
so pulling back $(\bigwedge^{n-1}Q)^{**}$
 to a desingularization, it is still a big
invertible sheaf. 
But we must then  check that this pull back 
is a subsheaf of some sheaf of differential forms.
We accomplish this  by choosing  the polynomial $F_{mp}$ 
so as to make an
 explicit resolution
 straightforward.

Recall that the non-smooth points of $Z$ are are given precisely by
the critical points of the polynomials $f_{mp}(x_1, \dots, x_n)$ \,\,
(Exercise 17). 
It is easy to desingularize $Z$
  under the the following non-degeneracy assumption:

\medskip
{\bf Assumption \Thm:\/}\label{nondeg} 
 {\it{Each dehomogenization of $F_{mp}$   is a polynomial
with only non-degenerate
critical points.\/}}

This means that for each $i = 0, \dots, n$, the polynomial 
$$
f = F_{mp}(X_0, \dots, X_{i-1}, 1, X_{i+1}, \dots, X_n)
$$
in $n$ variables  
has 
only non-degenerate critical points.  As usual, 
a critical point of a polynomial $f$  is a point where 
all the partial derivatives of $f$ vanish, and the critical  point  is
non-degenerate if the 
determinant of the Hessian matrix of second derivatives does not vanish there.
Here, ``point'' means point defined over the algebraic closure of the ground 
field. Such $F_{mp}$ exist over any infinite field (with some exceptions in 
characteristic two); see \refto{exist}. 

\bigskip
Assuming now that $F_{mp}$ has the non-degeneracy condition
described above, we  now complete the proof of Theorem \refto{charp}
by desingularizing $Z$ and 
verifying that  $(\bigwedge^{n-1}Q)^{**}$ pulls
back to a subsheaf of regular differential forms.

\medskip
The advantage of non-degenerate critical points is
that, 
 after possibly enlarging the ground field,
the affine equation of the hypersurface $Z$ can be assumed of the 
the form 
$$
y^p =  c 
+ x_1x_2 + x_3x_4 + \dots + x_{n-1}x_n + \,\, f_3
\,\,\,\,{\text{  if $n$ is even, or}}
$$
$$
y^p = c + x_1x_2 + x_3x_4 + \dots + x_{n-2}x_{n-1} + x_{n}^2 +
 f_3
  \,\,\,\,\,{\text{ if $n$ is odd and $p > 2$.}}
$$
where the $y, x_i$'s are local coordinates at a non-smooth point of $Z$,
$c$ is a constant, 
and $f_3 = f_3(x_1, \dots, x_n)$ 
 is a polynomial
of order three 
or more
 in the $x_i$ (see Exercise 21).
 Fortunately, desingularizing such a hypersurface is easy.

\medskip
{\bf Exercise 20:} 
 Show that if $f$ has only  isolated non-degenerate
critical points, then the affine
hypersurface defined by 
$y^p - f$  becomes smooth
 upon 
  blowing up each 
non-smooth point (over the algebraic closure of the ground field),
regardless
of the characteristic of the ground field.

\subsubhead{\Thm. Verification that $\Cal M$ is a subsheaf of differential
forms}\label{ver}\endsubsubhead
Having shown that $Z$ can be smoothed by blowing up points, 
let $Z' @>{q}>>  Z$ be this desingularization
of $Z$.
Assuming that the ground field is characteristic $p$, consider the sheaf 
 $$
\Cal M = q^*(\bigwedge^{n-1} Q)^{**} = q^*\pi^*\Cal O_{\PP}(mp - n-1)
$$
as in \refto{special}. 
We know that $\Cal M$ is big, and we wish to show that it
is contained in 
$\bigwedge^{n-1}\Omega_{Z'}$.  This is just a matter
of computing local generators for $\Cal M$ and comparing them to 
local generators for $\bigwedge^{n-1} \Omega_{Z'}$.

We return  to the somewhat mysterious definition of $Q$.
Recall that $Q$ is the cokernel 
of the  very special  $\Cal O_Z$-module map 
$$
\pi_Z^*\Cal O_{\PP}(-mp) \overset{d}\to\rightarrow \pi_Z^*\Omega_{\PP},
$$ 
defined in \refto{special}.  Think of
 $d$  as the pull back of 
a map of $\Cal O_{\PP}$-modules
$
\Cal O_{\PP}(-mp) \overset{d'}\to\rightarrow \Omega_{\PP},
$
 sending the local generator $f_{mp}$ to 
$$df_{mp} =
 \sum_{i =1}^n{\frac{\partial f_{mp}}{\partial x_i}} dx_i.
$$
(We reiterate that this map is 
deceptively subtle: its existence is a very special
consequence of the fact that the ground field is characteristic $p>0$;
Cf \refto{special}.)

 Taking the  $(n-1)^{st}$ exterior power of the cokernel $Q$  of
$d = \pi_Z^*d'$, 
 we have  convenient local generators 
$$
\eta_i 
= (-1)^i  \frac{dx_1 \wedge \dots \wedge
 \widehat{dx_i} \wedge \dots \wedge dx_n}{{\partial f_{mp}/\partial x_i}}
$$
for $\bigwedge^{n-1}Q$ 
on the open set where 
${\partial f_{mp}/\partial x_i}$ is non-zero.
Note that   $\eta_i = \eta_j$ whenever both are defined.
The locus where  no $\eta_i$ is defined is precisely the non-smooth
locus of $f_{mp}$. Since this set has codimension at least two,
this sheaf extends uniquely to 
a sheaf   $(\bigwedge^{n-1}Q)^{**}$  on 
all of $Z$. The extension can be defined as a subsheaf of the 
constant sheaf of rational differential forms on $Z$  generated by the
$\eta_i$.
By definition of $\Cal M$,
these pull back to  local generators of $\Cal M$ on the
desingularization $Z'$.

To check that $\Cal M \subset \bigwedge^{n-1}\Omega_{Z'}$, we 
 only need check what happens along the exceptional
fibers of $Z' @>>> Z$, since we already  know the inclusion holds
  on the smooth locus of $Z$. 
This is  a straightforward computation;
we work it out in one case below.

Let $y, x_1, \dots, x_{n}$ 
be local coordinates  for $Z$ near a non-smooth point,
and let 
$y, x_1', \dots, x_{n}'$ denote  local coordinates on the blowup $Z'$
of  the ideal $(y, x_1, \dots, x_{n})$, with
 $x_i = yx_i'$.
Computing the pull back of, say,  $\eta_{n}$ when $n$ is even and $p$ is 2,
 we have
$$
q^*\eta_n =  \frac{d(yx_1') \wedge  \dots \wedge d(yx_{n-1}')}
{\partial (y^2 + x_1x_2 + \dots +x_{n-1}x_{n}+ g)/ \partial x_n}
$$
where $g$ has order 3 or more in $(y, x_1, \dots, x_{n})$.
Performing the differentiation, we see that the denominator is 
$
x_{n-1} + h,
$  
 where $h$ is order 2 or more in 
 $(y, x_1, \dots, x_{n})$, which we write as 
 $y(x'_{n-1} + yh')$ in local coordinates  on the blowup,
with  $h'$ in $(y, x_1', \dots, x_n')$. 
Thus
$$
\align
q^*\eta_n
=&
 \frac{y^{n-1}(dx_1' \wedge \dots  \wedge dx_{n-1}') + \sum_{j=1}^{n-1}y^{n-2}
(dx_1' \wedge \dots \wedge dy \wedge \dots \wedge dx'_{n-1})}
{y(x'_{n-1} + yh')} \\
=& y^{n-3}\left(\frac{y(dx_1' \wedge \dots  \wedge dx_{n-1}') 
+ \sum_{j=1}^{n-1}
(dx_1' \wedge \dots \wedge dy \wedge \dots \wedge dx'_{n-1})}
{(x'_{n-1} + yh')}\right),
\endalign
$$
where the $j^{th}$ term in the sum has $dy$ in the $j^{th}$ position.

To check that the local generator $q^*\eta_n$ has no pole along the
 exceptional fibers,  
we can compute
 compute in any open set 
which intersects the  exceptional divisor $E$.
In the neighborhood considered above,
 the exceptional divisor $E$
is defined by $y$, so the 
 generator $q^*\eta_n$  vanishes along $E$ to at least 
order $n-3$. So  $\Cal M$ 
has no poles along $E$
 whenever $n \geq 3$.
Because the computation along each exceptional divisor is essentially
the same,
we conclude that $\Cal M$ is a subsheaf of $\bigwedge^{n-1}\Omega_{Z'}$,
whenever $n \geq 3$.

\subsubhead{\Thm. $Z$ is not ruled}\endsubsubhead
Now we are in position to apply Proposition \refto{rulecrit}.
The blowup variety $Z'$ is smooth and carries the
 big invertible sheaf $\Cal M$ 
that is a subsheaf of a sheaf of differential forms on $Z'$.
Proposition \refto{rulecrit} 
implies that $Z'$ cannot be separably uniruled.
 Because $Z'$ is  birationally equivalent to $Z$,
it follows also that $Z$ is not 
separably uniruled.  In particular, $Z$ is not ruled.

\medskip
This essentially completes 
the proof of Theorem \refto{charp}.
Any subvariety $Z$
of $\bold{ P}$  defined by an equation of the form   $Y^p - F_{mp}$,
where $F_{mp}$ satisfies
Assumption \refto{nondeg} (ensuring the singular points of $Z$ are 
non-degenerate)  and where the numerical constraints $pm - m < n+1 < pm$ hold,
is an example of a non-ruled projective variety 
whose  anti-canonical
sheaf is ample invertible. 
We now only need to observe that such $F_{mp}$ satisfying this assumption
{\it do\/} exist. 

\medskip

{\bf Exercise 21:} 
A critical point $P$ of a polynomial 
$f$ is   non-degenerate if  the determinant of the Hessian matrix 
$(\frac{\partial^2{f}}{{\partial x_i \partial x_j}})$ does not vanish at $P$,
or equivalently, if 
$\{\frac{\partial f}{\partial x_i}\}_{i = 1}^{n}$
generate the maximal ideal of $P$. 
Prove the following  Morse Lemma  for  polynomials
over an infinite field $k$.
\roster
\item 
 If the  characteristic  of $k$ is greater than two,
then a sufficiently general polynomial function of degree $d$ in $n$
 variables over $k$ has only non-degenerate critical points. 
\item 
If $k$ has characteristic two,  then every critical point of a
 polynomial in an odd number of variables is degenerate, where as the 
general polynomial function of an even number of variables has only
non-degenerate critical points. 
\endroster

\subsubhead{\Thm}\endsubsubhead\label{exist}
When the ground field is infinite, the Morse Lemma ensures that
{\it every\/}
sufficiently general choice of $F_{mp}$ satisfies Assumption \refto{nondeg}.
Thus, over an infinite field, there are many examples of non-ruled $Z$. 
It is not obvious that  such $F_{mp}$ exist over finite fields.
Fortunately, 
explicit  examples, due to Joel Rosenberg, show that 
 polynomials $F_{mp}$ satisfying Assumption \refto{nondeg} 
and  also satisfying the numerical constraints $mp - m < n+1 < pm$
exist  over {\it any\/} finite field.
A specific example, over any field of characteristic $p$, is 
 the hypersurface in the weighted projective space $\bold P$
defined by 
$$
Y^p - \sum_{i=0}^{{n}}X_i^{mp-1}X_{i+1} 
\tag{\Thm}\label{rosen}
$$
where the subscripts are taken modulo $n + 1$.
Here,  $n$ is any integer greater than two
satisfying $pm -m < n+ 1 < pm$. See 
 Appendix I.
\medskip

\subsubhead{\Thm. Summary}\label{summary}\endsubsubhead
We have established the following.
 Fixing positive integers $p, m$ and $n$
satisfying
$$
pm - m  < n+1 < pm \,\,\,\,\,\,\,\, 
$$
with $p$ prime and $n$ at least three, let $F_{mp}$ be a homogeneous
polynomial of degree $mp$ in $n+1$ variables.
If $Z$ denotes the hypersurface defined by $Y^p - F_{mp}$ in weighted
 projective space $\bold P$, the following hold for sufficiently
 general  choices of $F_{mp}$:
\roster
\item 
The anti-canonical sheaf of $Z$ is ample invertible;
\item 
When the characteristic is not $p$, the variety  $Z$ is smooth;
\item 
When the characteristic is $p$, the variety  $Z$ is not ruled  (even after
arbitrary base extension). 
\endroster
The precise meaning of ``sufficiently general'' here is
that $F_{mp}$ should define a projective hypersurface with only isolated 
non-smooth points, all of which are non-degenerate;
that is, $F_{mp}$ must satisfy Assumption \refto{nondeg}.

 Unfortunately, however, the varieties  $Z$ are  non-smooth in characteristic
$p$, precisely the case  where we have proven them non-ruled. In 
  Lecture 5, we will see that the same polynomial
defines
 a non-ruled smooth Fano variety  in characteristic zero.
In particular, these are examples of non-rational smooth Fano varieties.

\setcounter {Ch}  {4}
\heading  \Ch.\  Final Lecture 
\endheading
\label {lec5}

The goal of this lecture is to establish the 
existence of   a multitude of smooth 
non-rational
Fano varieties  in every dimension.
In fact, we construct specific families of {\it non-ruled\/}
 smooth Fano varieties. The main point is the method of the proof.
Here we content ourselves with  simple applications, producing
only  some examples of non-rational Fano varieties. With minor
modifications the method produces many more such examples, see
\cite{K95} or \cite{K97, V 5} for details.

For a homogeneous polynomial $F_{mp}$ of degree $mp$ let $Z=Z(F)$ 
denote the hypersurface 
$Y^p - F_{mp}$ in the weighted  projective space $\bold P = \PP(m,
1, 1, \dots, 1)$. This was
   constructed in \refto{construct},  and its key properties are
summarized in \refto{summary}.

\proclaim{\Sec. Theorem}\label{char0}
Fix an arbitrary ground field $k$ of characteristic zero. 
Fix  integers $m \geq 1,  n \geq 3$ and a prime  $p$ satisfying 
$$(p-1)m < n+1 < pm.$$ 
Then there exists a homogeneous polynomial $F_{mp} \in k[X_0, \dots, X_n]$
such that the corresponding variety $Z(F)$ 
is a smooth projective Fano variety that  is not ruled over $\bar k$. 
\endproclaim

\subhead{\Sec. Remarks}\endsubhead

\subsubhead{\Thm}\endsubsubhead
We will prove here only a weaker form Theorem \refto{char0} in full detail:
we will show that $Z(F)$ is not {\it rational.\/}   

\subsubhead{\Thm}\endsubsubhead
 It is quite likely that if $(p-1)m < n+1 < pm$ then
every smooth variety $Z(F)$ as above is nonrational. Unfortunately
this is not known.

\subsubhead{\Thm}\endsubsubhead
For most values of $m$, $n$, and $p$, one can use Theorem
\refto{char0} to  write down explicit
examples of smooth Fano varieties over $\QQ$ that are not ruled.
For instance, the examples of Rosenberg (see Appendix I)
show that if $n\not\equiv{-1}\mod p$ and $mp\ge3$
$$
\{Y^p-\sum_{i=0}^{n}X_i^{mp-1}X_{i+1}=0\}  \subset \bold P 
$$
is a smooth Fano variety that is not ruled.
See also \cite{K96, V 5.16.3} for other examples deduced
by a related theorem.

\subsubhead{\Thm}\endsubsubhead
 Let $V_{mp}$ be the vector space of  homogeneous polynomials
$F_{mp}\in k[X_0, \dots, X_n]$. By
  \cite{K97 IV.1.8.3}, there are countably many subvarieties
$W_i\subset V_{mp}$ such that if $F\not\in \cup_iW_i$ then $Z(F)$
is not ruled. Thus if
$k$ is uncountable, then finding just one non-ruled example
guarantees that most of the varieties $Z(F)$ are not ruled, hence
also not rational.

In any case, 
if $k$ is not algebraic over $\QQ$, then the Morse Lemma
guarantees that a generic choice
of $F_{mp}$ produces a non-ruled example; see the explanation 
\refto{generic}.

\bigskip
 Lecture 5 is devoted to the proof of Theorem \refto{char0}.

\medskip
The idea of the proof is simple.
We already know that 
$Z$ is a smooth Fano variety in characteristic zero \refto{-F}. 
The proof that $Z$ is non-ruled uses reduction to characteristic $p$.
 To get a rough idea how this is done,
suppose first that we wish to construct an example defined over $\QQ$.
Consider 
$F_{mp} \in  \ZZ[X_0, \dots, X_n]$, with $m, n$ and $p$ 
as in Theorem \refto{char0}.
  Consider the scheme 
$$
\proj \frac{\ZZ[Y, X_0, \dots, X_n]}{Y^p - F_{mp}} =  
Z_{\ZZ} \subset \bold P_{\ZZ},
$$
 a closed subscheme of 
the weighted projective
space  $\bold P_{\ZZ} = \proj \ZZ[Y, X_0, \dots, X_n]$.
There is a morphism 
$$
Z_{\ZZ} @>>> \spec \ZZ,
$$
whose special fiber over $(p)$ is  the  singular 
Fano variety as constructed in Lecture 4.
Furthermore,  we can arrange this so the special fiber is non-ruled: 
we need only choose 
$F_{mp} $ so that its reduction modulo $p$ 
 has only non-degenerate
critical points in each affine patch. 
Now the idea is to apply the 
 following theorem of Matsusaka \cite{Mats} about the
behavior of ruledness in families.

\proclaim{\Sec. Theorem\/}\label{specialfiber}
Let $V$ be a 
 discrete valuation ring 
 that is a localization of a finitely generated algebra
 over a field or over the integers.
{\footnote{In fact, the base $V$ can be any {\it excellent\/}
 discrete valuation ring, although we do not need to 
apply the theorem in such generality. For the
definition of excellent, see \cite{M, p260}.}}  
 Let 
  $K$  (respectively, $k$) denote its  
 quotient field  (respectively,  residue field).
Let $Z_S$ be a normal irreducible projective scheme over $S = \spec V$.
 If the generic
fiber of the natural projection $Z_S @>>> S$ is ruled over $K$,
 then each  irreducible component of the 
special fiber is ruled over $k$.
\endproclaim

{\bf Cautionary Remark.\/}
Theorem \refto{specialfiber} underscores the reason we are led to
consider non-ruled
 varieties in our quest for non-rational 
ones: ruledness is better behaved in families than 
rationality. 
We can not conclude a special member  of a family is rational
when we know the generic member is rational.
 For example,  a family
of  degree three 
hypersurfaces in $\PP^3$  has a smooth cubic surface as its 
generic member, but it can have singular members that are cones over
elliptic curves.
The generic member is rational,
 whereas the special member
of this family is only ruled. 

\medskip
{\bf \Thm. Remark.\/}\label{weak}
To keep things elementary, we prove here only the following weak form
 of Matsusaka's theorem: If the general fiber above is {\it rational,\/}
then the components of the special fiber are ruled.
This will be sufficient to conclude the existence of non-rational
 $Z(F)$. For the full
proof of Theorem \refto{specialfiber},
 the reader is referred to \cite{K96, p 184}.

\medskip
Before proving Theorem \refto{specialfiber}, we show how to use it
to deduce the existence of non-rational Fano varieties in characteristic zero.

\demo{Deduction of  Theorem \refto{char0} from Theorem \refto{specialfiber}}
Let $Z$ be the hypersurface in $\bold{P}$  defined over $k$.
Choose a finitely generated $\ZZ$-algebra 
$A $ contained in $k$ over
which  $Z$ is defined, that is, containing all the coefficients
of the defining equation for $Z$.
Having chosen $A$, let $S = \spec A$, and observe that there
is an $S$-scheme $Z_S$ 
such that 
the scheme $\spec k \times_{S}  Z_S$ is naturally isomorphic
to $Z$.

 Assuming $Z$ were ruled,
fix a  birational map 
$$ W \times_k \PP^1_k @>\phi>>  Z.$$ 
Now choose $A$ large enough so that 
it  contains all the elements of $k$ necessary to describe 
the $k$-scheme $W$ and the map $\phi$.
This gives rise to an $S$-scheme $W_S$ and an $S$-scheme map 
$$
W_S \times_S \PP^1_S  \overset{\phi}\to\rmap  Z_S
$$
which is a birational equivalence.

The base scheme $S$ may be replaced by the spectrum of
a   discrete valuation ring  as follows.
Normalizing if necessary, there is no loss of generality in assuming
the base ring
 $A$ is  normal. Now localizing  $A$ at a minimal prime of $(pA)$, we achieve
a discrete valuation ring, say  $V$.
By base change, we replace $S = \spec A$ with the two-point scheme 
$S = \spec V$, and assume the map $Z_S @>>> S$ is a family over spectrum of
the discrete valuation ring $V$.

The generic fiber of $Z_S @>>> S$ would clearly be ruled.
However, the special fiber is a variety of prime  characteristic $p$,
defined in $\bold{P}$ 
by an equation  of the form $Y^p - {F_{mp}}$.
We proved in Theorem \refto{charp}
 that such a variety is {\it not\/}
ruled in general.
Indeed, if   $F_{mp}$ is any homogeneous polynomial
having only non-degenerate critical points locally (ie, satisfying
  Assumption \refto{nondeg})
 over  the residue field $V/uV$, then the special fiber $Z_p$ is not ruled.
So lifting $F_{mp}$ to a polynomial defined over $V$, we have a scheme
$Z_S$ in $\bold P_{S}$ whose fiber over the generic point of $S$ is 
a smooth  Fano   variety $Z \subset \bold {P}$, 
defined 
over the fraction
 field of $V$, which can not be ruled by Theorem \refto{specialfiber}.
This contradiction completes the proof that $Z_k$ is not ruled.
\enddemoo

\medskip
{\it \Thm. Remark.\/}\label{generic} 
If $k$  
 is non-algebraic
over $\QQ$, we expect some of the coefficients of $F_{mp}$ to be
non-algebraic over $\QQ$. This forces $V/uV$ to be an infinite
field in general. Now the Morse Lemma implies that  sufficiently general 
 choices of $F_{mp}$ modulo  $uV$ 
satisfy Assumption \refto{nondeg}
(see Exercise 21). This means that 
all sufficiently general choices of $F_{mp}$ defined over 
 $k$ give rise to varieties $Z$ that are smooth Fano non-ruled varieties.

\medskip
We now complete the proof of the existence of non-rational Fano varieties
by proving the weak form of Theorem \refto{specialfiber}.
With some extra work, 
 Theorem \refto{specialfiber} can be proved in full,
but we refer to  \cite{K96, p. 184} for the details.

\demo{Proof Theorem \refto{specialfiber} in the weak form \refto{weak}}
Let $\xi$ be the generic point of $S$ and assume the fiber $Z_{\xi}$ 
of $Z_S @>>> S$
over $\xi$ is rational. 
Since $S$ is birationally equivalent to ${\xi}$, 
a birational map
$$
\PP^n_S \times_S \xi =   \PP^n_{\xi} \rmap Z_{\xi} = Z_S \times_S \xi
$$
defines a birational map
$$ \PP^n_S \overset\phi\to\rmap Z_S
$$
over $S$. 
(If the generic fiber of $Z_S @>>> S$ is only ruled, we can
 find a $S$-scheme $\PP^1_S \times_S W$ mapping birationally onto $Z_S$
over $S$.  But  to make the following proof work,  this scheme must be 
 both regular and proper over $S$.)

\medskip
Let $\Gamma_S$ be the  normalization of the (closure of the)
graph (in $\PP^n_S \times_S Z_S$)  of the birational map $\phi$. 
The normalization is a finite, birational morphism.
{\footnote{In general, the normalization map
of an arbitrary scheme can fail to  be finite, although it is finite
for schemes of essentially of  finite type over $\ZZ$, or more generally
for any excellent scheme.}}  
So composing with the natural projections, we have proper birational
morphisms 
$\Gamma_S @>{\pi_1}>>  Z_S$ and
 $\Gamma_S @>{\pi_2}>> \PP^n_S$.

Consider the special fiber 
 $Z_p$ of $ Z_S @>>> S$.
Being defined by a single equation, namely the
pullback of the uniformizing parameter $u$, all components of $Z_p$ have
codimension one. 
Likewise, 
the components of the special fiber
$\Gamma_p$ \, of \, $\Gamma_S @>>> S$ 
are 
 all  codimension one.
Now, 
because $\Gamma_S @>{\pi_1}>>   Z_S$ is a proper birational map
of normal schemes, it is an isomorphism in codimension one (on the base).
 This means that for each irreducible component of $Z_p$, there 
corresponds a 
 unique  irreducible 
component $\Gamma_0$ of $\Gamma_p$ mapping birationally to it.
Therefore, 
it suffices to show that the  reduced 
irreducible divisor $\Gamma_0$   is ruled. 

Consider the restriction of  $\pi_2$ to $\Gamma_0$.
This gives a   morphism $\Gamma_0 @>>> \PP^n_p$,
where 
$$
\PP^n_p = \PP^n_S \times_S \spec(V/uV)
$$
 is the special fiber of 
$\PP^n_S @>>> S = \spec V$.
If  $\Gamma_0 @>>> \PP^n_p$ is
 birational, the proof of Theorem \refto{specialfiber} is complete:
then $\Gamma_0$, and hence the birationally equivalent variety  $Z_p$,
 would be birationally equivalent to 
$\PP^n_p$.

On the other hand, the map 
$\Gamma_S @>{\pi_2}>> \PP^n_S$ is a proper birational morphism,
so if its restriction to the divisor $\Gamma_0$ is not birational, then 
$\Gamma_0 $ must be an exceptional divisor for this map.  But exceptional
divisors of proper birational maps to regular schemes are {\it always\/}
ruled, as the following theorem of Abhyankar shows
\cite{Ab, p336}.  
\enddemo

\proclaim{\Sec. Theorem}\label{ram}
Let $Y @>{\pi}>> X$ be a proper birational morphism of irreducible schemes,
with $Y$ normal and $X$ regular.
{\footnote{Again, some very mild reasonability  condition on $X$
  is required: it is sufficient if
 $X$ is  of finite type over a localization of a 
finitely generated algebra over a field or $\ZZ$. Indeed, $X$ can be 
any excellent scheme.}}
Then every exceptional divisor of $\pi$ is ruled over its image.
That is, if $E$  is an integral subscheme of $Y$
 of codimension 
one, whose image $E'$
 has codimension greater than one in $X$, then 
$E$ is birationally equivalent over $E'$ 
to a scheme $W \times_{E'} \PP^1_{E'}$. 
\endproclaim

\medskip
Abhyankar's proof uses valuation theory. Rather than reproduce his proof here
in full generality, we provide a  nice geometric proof
in the case where  $Y$ is of finite type over $X$, which is certainly
sufficient for our purposes.

\medskip
We first point out that  when
 $Y$ and $X$ are algebraic varieties
defined over a field $k$ of characteristic zero,
Theorem \refto{ram} follows easily from Hironaka's theorem on 
resolution birational maps.
{\footnote{Although Abhyankar's 1956 proof came first.}}
In this case, 
the morphism  $\pi$ factors through a sequence of
blow-ups along smooth centers:
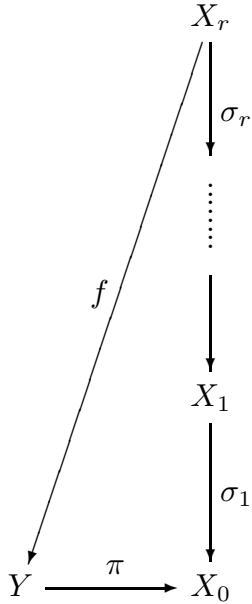
\diagram
{} &  {}   & X_{r}  \\
{}&    \ldTo(2, 6)^{f} & \dTo>{\sigma_r} \\
{}&   {} &\dDots \\
{}&  {} & \dTo \\
{}&    {} & X_1 \\
{} &  {} &  \dTo>{\sigma_1} \\
Y & \rTo^{\pi} & X_0 \\
\enddiagram
Here, each 
$X_{i+1} \overset{\sigma_{i+1}}\to\rightarrow X_{i}$
is a  blowing up
along a nonsingular center, and 
$X_r \overset{f}\to\rightarrow  Y$
is a birational morphism from the non-singular variety 
 $X_{r}$.
To see that $E$ is ruled, there is no harm in 
replacing $Y$ by $X_r$ and $E$ by its birational transform on $X_r$.
Because  $E$ is exceptional for the composition of the blowups $X_r @>>> X$,
its image on some $X_i$  must be 
an exceptional divisor for some blowup $\sigma_i$. 
But the  exceptional divisor of a
blow-up  of a non-singular variety along a non-singular center 
is a projective space bundle over the center; in particular,
 such exceptional divisors, including $E$, must be
ruled.

\bigskip
Because we are interested only in birational properties,
it is not actually  necessary to use Hironaka's deep theorem.
The idea can be adapted as follows. 
 We construct the tower
of blowups $X_i @>{\sigma_i}>> X_{i-1}$
 by blowing up the image $E_{i-1}$ of $E$ on $X_{i-1}$,
but always 
restricting to the regular (non--singular)
 loci of the $E_{i-1}$, so that each  
$X_i$ is regular.
Again, the exceptional fiber of the blow-up of a regular scheme along a regular
subscheme is a projective space bundle over the center of the blowup.
So again, we  can show that $E$ is ruled by showing that
the image of $E$ on some $X_i$ must be the exceptional
divisor for  some blowup $\sigma_i$.

This is shown by
keeping track of 
a numerical invariant that drops with each non-trivial blowup.
This numerical invariant can be taken to be the order of vanishing
along $E$
of the pull back of a local generator of $\omega_{X_i}$
to $\omega_{Y}$, provided that one can
make sense of  the sheaves $\omega_{X_i}$ and
 $\omega_{Y}$ and that one can
define a pullback map $\pi^*\omega_{X_i} @>>> \omega_Y$.
For example, when $Y$ and $X$ are of finite type over an algebraically closed 
field, there is no problem making sense of $\omega_{X}$ and
$\omega_{Y}$ and the argument is easily adapted to this case.

\medskip
In carrying out this argument, 
complications arise when the scheme $X$ (and $Y$) 
are not defined over some base field. 
When $X$ and $Y$ are both smooth over some base scheme $S$,
one can try to work with the relative canonical modules $\omega_{X/S}$.
This sometimes works (for instance, if $E$ is flat over $S$), 
but unfortunately, it breaks down precisely in the case we need it.
The trouble arises because we must work 
with  {\it regular\/} schemes that {\it may not be smooth\/}
over the  base scheme. Indeed,  the following situation is typical: the scheme 
 $X$ may be $\A^n_S$  with $S$ the 
two-point scheme $\spec{\ZZ_{p}}$.
We will blow up 
a regular subscheme $E_0$ which maps to the closed point of 
 $S$, say the closed point 
defined by $(p, x_1, \dots, x_n)$.
The resulting  blow-up scheme $X_1$ is regular, but it is
 not smooth over $S$.
In this case, it is hard to define a relative canonical module
$\omega_{X_1/S}$ that has the properties need to carry out the argument
along the lines suggested above.

 However, we will be able
to adapt the proof of Abhyankar's theorem in our case by
working with the  relative canonical modules for the
birational maps $X_i @>>> X_{0}$. 
In fact, the duals of these canonical modules,
 the so-called Jacobian ideals, are more convenient to 
work with.

\bigskip

\subsubhead{\Thm. A Digression on Relative Canonical Modules
 and Jacobian Ideals}
\endsubsubhead\label{can}
Let $Y$ be a scheme of finite type over $X$, and suppose that $Y @>>> X$
has relative dimension $d$.
We say that $Y$ is smooth over $X$ if the sheaf of  relative
K\"ahler differentials  $\Omega_{Y/X}$ is a locally free
$\Cal O_Y$  module  of rank $d$.
In this case, the relative canonical module $\omega_{Y/X}$
is defined to be the invertible sheaf  $\bigwedge^d \Omega_{Y/X}$.

When  $Y$ is normal and smooth in codimension one over $X$,
the relative canonical module  $\omega_{Y/X}$  can be defined as 
the  unique reflexive $\Cal O_Y$ module that agrees
with the above construction on the smooth locus of $Y @>>> X$.
Equivalently, $\omega_{Y/X}$ is the double dual of the $\Cal O_Y$
module  $\bigwedge^d \Omega_{Y/X}$. Although this canonical module is
 not necessarily invertible, it
still can be interpreted as the ``determinant'' of
 $\Omega_{Y/X}$ via the natural map
$$
\wedge^{d}\Omega_{Y/X}\to \omega_{Y/X} = (\wedge^{d}\Omega_{Y/X})^{**}
$$
which  is 
neither  injective nor surjective in general.

\medskip
This method of defining the canonical module fails when  $Y$ is not
smooth in codimension one over $X$ (for example, when $Y @>>> X$ is
a blowup). Let us compute the ``determinant'' of $\Omega_{Y/X}$
in a different way so that it will generalize to this case. 
Fix an embedding   $Y \overset{i}\to\hookrightarrow W$ of $Y$ in a
smooth $X$ scheme $W$, 
for instance, 
$W$ may be taken to  be an open subset of affine space over $X$.
Consider the 
 conormal complex:
 $$
\Cal I_{Y}/\Cal I_Y^2 @>d>> i^*\Omega_{W/X} @>>> \Omega_{Y/X} @>>> 0
$$
where $\Cal I_Y \subset \Cal O_W$ is the ideal sheaf of $Y$ in $W$.
If $\Cal I_Y$ is locally generated by a regular sequence, as it must be,
for instance, when both $X$ and $Y$ are regular, then  the conormal sequence
is exact also on the left:
$$
0 @>>> \Cal I_{Y}/\Cal I_Y^2 @>d>> i^*\Omega_{W/X} @>>> \Omega_{Y/X} @>>> 0.
$$ 
This  suggests a method for computing the 
``determinant" of  $\Omega_{Y/X}$ when  $X$ and $Y$ are regular. In
this case $\Cal I_{Y}/\Cal I_Y^2$ and $i^*\Omega_{W/X}$ are both
locally free, hence we can propose the definition
$$
\omega_{Y/X}:=\wedge^ni^*\Omega_{W/X}\otimes (\wedge^{n-d}\Cal
I_{Y}/\Cal I_Y^2)^{-1},
$$
where $n$ is the relative dimension of $W$ over $X$ and $d$ is the relative
dimension of $Y$ over $X$.
This
 agrees with our previous definition, but makes sense even when 
$Y @>>>X$ is not smooth in codimension one. 
The module 
$\omega_{Y/X}$ is invertible (provided $Y$ is generically smooth over $X$
so that ranks are as expected).
 We do not need to worry about
the dualizing properties of the  sheaf.

Now, in order to prove Abhyankar's theorem,
we must consider the case where  $Y @>>> X$ is a birational morphism of regular
schemes. The relative dimension is zero. 
The map of rank $n$ locally free $\Cal O_Y$-modules
$$
 \Cal I_{Y}/\Cal I_Y^2 \overset{d}\to{\hookrightarrow} i^*\Omega_{W/X}
$$
 gives rise to a map of invertible $\Cal O_Y$-modules
$$
\bigwedge^n \Cal I_{Y}/\Cal I_Y^2 \hookrightarrow \bigwedge^n i^*\Omega_{W/X}.
$$
Tensoring with 
$(\bigwedge^n i^*\Omega_{W/X})^{-1}$
we get an exact sequence
$$
0 @>>> 
\bigwedge^n \Cal I_{Y}/\Cal I_Y^2 \otimes (\bigwedge^n i^*\Omega_{W/X})^{-1}
@>>> 
\Cal O_Y @>>> \Cal  Q
@>>> 0.
$$
Here, $\omega_{Y/X} = 
\bigwedge^n i^*\Omega_{W/X}  \otimes 
(\bigwedge^n \Cal I_{Y}/\Cal I_Y^2)^{-1} $
so that its dual,
$\omega_{Y/X}^{-1} = 
\bigwedge^n \Cal I_{Y}/\Cal I_Y^2 \otimes (\bigwedge^n i^*\Omega_{W/X})^{-1}
$ is a sheaf of ideals 
 in $\Cal O_Y$.
It is   often called the
Jacobian ideal and denoted by $\Cal J_{Y/X}$.
 Note also that 
 $\Cal Q$ is some torsion $\Cal O_Y$-module supported on the non-smooth
locus of $Y @>>> X$, so that the Jacobian ideal defines the
  non-smooth locus of $Y@>>> X$.

To explain the name,
choose local coordinates
 $x_1, \dots, x_n$ for $W$ over $X$,
such that the $dx_i$ are a free basis for
 $\Omega_{W/X}$. 
 Suppose that
$\Cal I_Y$ is defined locally by the regular sequence $f_1, \dots, f_n$.
Then the map of free $\Cal O_Y$ modules
$$
 \Cal I_{Y}/\Cal I_Y^2 \overset{d}\to{\hookrightarrow} i^*\Omega_{W/X}
$$
sends the class of a generator $\bar f_i$ to 
$\sum_{i=1}^n\frac{\partial{f_i}}{{\partial x_j}} dx_j$. 
In other words, 
the map is 
defined by the Jacobian matrix 
$\left(\frac{\partial{f_i}}{\partial{x_j}}\right)$,
so that the map 
$$
\bigwedge^n \Cal I_{Y}/\Cal I_Y^2 \hookrightarrow \bigwedge^n i^*\Omega_{W/X}
$$
is defined by its determinant.  In particular, 
 the Jacobian ideal $\Cal J_{Y/X}$  is locally generated by 
this Jacobian determinant.

\medskip
The proof of Theorem \refto{ram} rests on the simple
observation that Jacobian ideals are multiplicative. 

\medskip
{\bf Exercise 22:\/}
If $Z @>\pi>> Y @>>> X$ are finite type 
birational maps of regular schemes, then 
$$
\Cal J_{Z/X} = (\Cal J_{Z/Y})(\Cal J_{Y/X}\Cal O_Z)
$$
 as ideals of
$\Cal O_Z$. 
 Here $ \Cal J_{Y/X}\Cal O_Z$  denotes the ideal of $\Cal O_Z$ 
generated by the pullbacks of generators of 
  $ \Cal J_{Y/X}$; because  $ \Cal J_{Y/X}$ is invertible, 
this is the same as  $\pi^* \Cal J_{Y/X}$.

\bigskip

 We can now give an easy proof of Abhyankar's theorem;  the idea  is 
from 
a paper of  B. ~Johnston 
\cite{Jo}.

\demo{Proof of Theorem \refto{ram} assuming $Y$ is of finite type over $X$}
Because $Y$ is normal, it is regular in codimension one. So we are 
free to replace $Y$ by an open set containing the generic point of $E$
so as to assume that $Y$ is regular.{\footnote{This is where the 
 excellence
hypothesis is used: we  must assume that the 
 regular locus is open, which holds, of course, when $Y$ is 
of finite type over a localization of a finitely generated algebra over
a field or $\ZZ$.}}

The divisor $E$ is exceptional
for the given  birational map $Y @>{\pi}>> X$,
 so its image under $\pi$ is
a subscheme of codimension at least two. Let us denote this 
image subscheme by $E_0$. Note that $E_0$ is a reduced and irreducible 
closed subscheme of the regular scheme  $X$.
 
The subscheme $E_0$ need not be regular.
 However, because it is reduced, 
 the locus of its non-regular  points is a proper
closed subscheme.
So we can replace $X$ by an open set $X_0$ in which  $E_0$ is regular.
Let $X_1 @>{\sigma_1}>> X_0$ be the blow-up of the regular scheme
$X_0$ along the regular subscheme $E_0$. 
The resulting scheme $X_1$ will be regular, and the 
exceptional fiber of the blowup will be a 
 projective space bundle over the center $E_0$ of the blowup.
In particular, the exceptional divisor is ruled.

 Let $E_1$ be the image of $E \subset Y$  
  in $X_1$ under the rational map
$Y \overset{{\sigma_1^{-1}\circ \pi}}\to\rmap X_1$.
   Of course, $E_1$ must be contained in 
the exceptional set for $\sigma_1$, but it may be strictly smaller.
If $E_1$ is codimension one in $X_1$, then $E_1$  must be this 
exceptional divisor. In this case, $E$ is ruled and the proof is complete.

 Otherwise, 
 $E_1$ has codimension larger than one in $X_1$ and we repeat the process
of replacing $X_1$ and $E_1$ by an open subset on which $E_1$ is regular
and blowing up along $E_1$. 
In this way, we construct 
a sequence of blowups $X_i @>\sigma_{i}>> X_{i-1}$.
Each $X_i$ is regular (but not necessarily smooth over any base scheme)
and each exceptional fiber is ruled and
 contains the image $E_i$  of $E$.
The 
process terminates  
(meaning $\sigma_{i+1}$ is an isomorphism)
if and only if $E_i$ is  codimension one in $X_i$.
If the process terminates, the proof is complete,
because then  $E$ is birational to  the exceptional 
divisor of some $\sigma_i$, and so $E$ must be  
ruled.  

\bigskip
\noindent
{\it{\Thm. Termination of the process.}\/} If  the process
does not terminate,
 we have a sequence
of  blowings up of regular schemes
$$
X_0 @<{\sigma_1}<< X_1 @<{\sigma_2}<< X_2 @<{\sigma_3}<< X_3 \dots 
$$
where no $\sigma_i$ is an isomorphism (we say ``$\sigma_i$ is a
non-trivial blow-up''). Such a 
 non-trivial blow-up is never smooth (nor even flat!).
Since 
the  Jacobian ideal $\Cal J_{X_i/X_{i-1}} \subset \Cal O_{X_i}$
 defines the non-smooth locus of the blowup $X_i @>>> X_{i-1}$, 
none of these Jacobian ideals can be the unit ideal.
Indeed, because  the blowup $X_i @>>> X_{i-1}$
 is not smooth along the exceptional divisor, the Jacobian ideal remains
a proper ideal after localizing along any component of an exceptional divisor. 

The rational  map $Y \overset{\pi_i}\to{\rmap}X_i$
is a morphism on some open set $Y_i$ containing the generic
point  of $E$.  In particular,
by Exercise 22, 
the morphisms 
$$
Y_i @>>> X_i @>>> X_0
$$
induce a multiplicative relation of Jacobian ideals in 
$\Cal O_{Y_i}$:
$$
 \Cal  J_{Y_i/X_0} = 
(\Cal  J_{Y_i/X_i})  (\Cal  J_{X_i/X_0}  \Cal O_{Y_i}).
$$
Localizing along $E$, we have a multiplicative relation
of {\it proper\/} ideals 
$$
 \Cal  J_{Y_i/X_0} \Cal O_{Y, E} = 
( \Cal  J_{Y_i/X_i}  \Cal O_{Y, E})(  \Cal  J_{X_i/X_0}    \Cal O_{Y, E}) 
$$
in the discrete valuation ring $\Cal O_{Y, E}$.
In particular,
the 
pullback of each Jacobian ideal
$  \Cal  J_{X_i/X_0 } $ to $  \Cal O_{Y, E}$
strictly contains the fixed ideal  
$\Cal  J_{Y_i/X_0} \Cal O_{Y, E} $. Since this  latter ideal of
the local ring $\Cal O_{Y,E}$
depends only on a small neighborhood of $E$ in $Y$, we denote it
by  $\Cal  J_{Y/X_0}$. 

Likewise, using the multiplicative property for Jacobian ideals
for the blowups $X_{i+1} @>>> X_{i} @>>> X_{0}$, we
see that  after  pulling back to $Y$ and localizing along $E$,
the ideal 
$
 \Cal  J_{X_{i}/X_0} \Cal O_{Y, E} $
strictly contains the ideal  
$ \Cal  J_{X_{i+1}/X_0} \Cal O_{Y, E}. $
We are led to a sequence of proper inclusions
$$
 \Cal  J_{Y} \subsetneq \dots \subsetneq
  \Cal  J_{X_i} \subsetneq  \Cal  J_{X_{i-1}} 
\subsetneq \dots \subsetneq  \Cal  J_{X_1} \subsetneq
 \Cal  J_{X_0}
$$
in the  discrete valuation ring 
 $  \Cal O_{Y,E }$ (the notation for ``relative to $X_0$'' and ``localize along
$E$'' has been suppressed).

This leads immediately to a contradiction: fixing a uniformizing
parameter $t$ for 
$  \Cal O_{Y,E }$,  and setting  $\Cal J_Y = (t^m)$,
it is obvious that at most $m$ ideals can be properly contained
between  $\Cal J_Y $ and
$  \Cal O_{Y,E }$. The process must terminate after at most $m$ blowups,
 and the proof is complete.
\enddemoo

{\it \Thm.  Remark.\/}
One can associate the numerical invariant given by the length of
$$
\frac{\Cal J_{X_i}}{\Cal J_{Y}}
$$ 
 to each blowup.  
The proof showed that this number is strictly decreasing 
for a non-trivial blowup.
This number can be interpreted as the ``discrepancy along $E$'' between 
differentials on $Y$ and on $X_i$.
In this sense, the proof we have given above
is very close in spirit to the proof 
we suggested in the  classical case.
The difference is that the differentials here are relative to the
scheme $X_0$, whereas 
 in the classical case, the
differentials are relative to the ground field.

\bigskip
We now 
know that  there are a host of smooth non-ruled  Fano varieties of every
dimension greater than two. In particular, there
are a host of non-rational Fano varieties. We have  a procedure for 
constructing examples. In Appendix I,  some explicit examples
over $\QQ$ are described.

\bigskip
There are many variations on this basic   method for constructing
non-ruled varieties. For example,  in \cite{K95}, Koll\'ar shows that 
a very general hypersurface of degree
$d$ in $\PP^n_{\CC}$ is not ruled whenever $d$ satisfies
$$
d \geq 2 \left\lceil\frac{n+3}{3}\right\rceil.
$$
Here, $\lceil x \rceil $ denotes the least  integer greater  than or
equal to $x$, and 
 ``very general'' means that there
is a countable union of subvarieties in the space all
  hypersurfaces in $\PP^n$ that must be avoided.
This result, however, does not produce any specific example of a non-ruled
hypersurface in $\PP^n$.  Further applications of this method
appear in \cite{K97}.

\newpage

\setcounter {Ch}  {5}
\heading  \  Appendix I: Polynomials with non-degenerate critical 
points over finite fields \\ 
by Joel Rosenberg 
\endheading
\label {joel}

In this appendix, specific examples  of  polynomials over finite fields 
with only non-degenerate critical points are recorded.
This establishes the existence of non-ruled Fano varieties over every 
field of characteristic zero.

\proclaim{Proposition}
Given  a prime $p$, and integers $n$ and $m$ with $n>0$,
$n\not\equiv{-1}\mod p$, and $mp\ge3$, let
$F\in\Bbb{F}_p[x_0,\dots,x_{n}]$ be the homogeneous polynomial of
degree $mp$
$$F(x_0,\dots,x_{n})=\sum_{i=0}^{n}x_i^{mp-1}x_{i+1},$$
where we understand subscripts to be taken mod $n+1$. Then any
dehomogenization $f$ of $F$
$$f(x_0,\dots,\hat{x_i},\dots,x_{n})=
F(x_0,\dots,x_{i-1},1,x_{i+1},\dots,x_{n})$$
will have only isolated critical points in $\bar \Bbb F_{p}$,
and all of them will be non-degenerate.
\endproclaim

\demo{Proof}
 From the cyclic symmetry of $F$, it is clear we need
only consider the dehomogenization 
$$
f(x_1,\dots,x_{n})=
F(1,x_1,\dots,x_{n}).
$$
Then we have
$$
\align
\frac{\partial F}{\partial x_i}&=x_{i-1}^{mp-1}-x_i^{mp-2}x_{i+1},\\
\frac{\partial^2 F}{\partial x_i\partial x_{i+1}}&=-x_i^{mp-2},\\
\frac{\partial^2 F}{\partial x_i^2}&=2x_i^{mp-3}x_{i+1},
\endalign
$$
and all other second partials of $F$ are zero. 
We find that any critical points of $f$ will have
$$
\align
 1-x_1^{mp-2}x_2&=0,\\
 x_1^{mp-1}-x_2^{mp-2}x_3&=0,\\
 x_2^{mp-1}-x_3^{mp-2}x_4&=0,\\
 \vdots\\
 x_{n-2}^{mp-1}-x_{n-1}^{mp-2}x_n&=0 \\
 x_{n-1}^{mp-1}-x_{n}^{mp-2}&=0 
\endalign
$$
We conclude that the critical points of $f$ are exactly those points
with
$$x_i=\zeta^{\left(\sum_{j=0}^{i-1}(1-mp)^j\right)}$$
for some $\zeta$ with
$$\zeta^{\left(\sum_{j=0}^{n}(1-mp)^j\right)}=1.$$
In particular, all critical points will be isolated and
will  have all their coordinates
nonzero.

To see that these points are non-degenerate, we write down the Hessian of $f$,
$$H=
\pmatrix
2x_1^{mp-3}x_2 & -x_1^{mp-2}    &                &        &&\\
-x_1^{mp-2}    & 2x_2^{mp-3}x_3 & -x_2^{mp-2}    &        &&\\
               & -x_2^{mp-2}    & 2x_3^{mp-3}x_4 &        &&\\
               &                &                & \ddots &&\\
&&&& 2x_{n-1}^{mp-3}x_{n} & -x_{n-1}^{mp-2} \\
&&&& -x_{n-1}^{mp-2}        & 2x_{n}^{mp-3}
\endpmatrix
$$
and compute its determinant at each of the critical points.
If we let $H_j$ be the upper left $j\times j$ submatrix of $H$, and
$h_j=\det(H_j)$, we see that we have a recursion for $\det(H)=h_{n}$:
$$
\align
h_0&=1,\\
h_1&=2x_1^{mp-3}x_2,\\
h_j&=2x_j^{mp-3}x_{j+1}h_{j-1}-x_{j-1}^{2mp-4}h_{j-2},
\text{ for }1<j\le n, 
\endalign
$$
where we understand $x_{n+1}$ to equal 1. We will show that modulo the
ideal of first partials $(\frac{\partial f}{\partial x_i})$, this
reduces to
$$h_j\equiv(j+1)\prod_{i=1}^j x_i^{mp-3}x_{i+1}.$$
Clearly this is true for $j=0$ and $j=1$. Now, inductively, for $1<j\le n$,
$$
\align
h_j&=2x_j^{mp-3}x_{j+1}h_{j-1}-x_{j-1}^{2mp-4}h_{j-2} \\
&\equiv2x_j^{mp-3}x_{j+1}h_{j-1}-x_{j-1}^{mp-3}x_j^{mp-2}x_{j+1}h_{j-2} \\
&\equiv 2jx_j^{mp-3}x_{j+1}\prod_{i=1}^{j-1} x_i^{mp-3}x_{i+1}
-(j-1)x_{j-1}^{mp-3}x_jx_j^{mp-3}x_{j+1}\prod_{i=1}^{j-2} x_i^{mp-3}x_{i+1} \\
&=(j+1)\prod_{i=1}^j x_i^{mp-3}x_{i+1},
\endalign
$$
as desired. In particular, $\det(H)=h_{n}$ is equal to $n+1$ times a
monomial, which is nonzero for any critical point of $f$. So $f$, and
similarly any dehomogenization of $F$, has no degenerate critical points.
\enddemoo

These examples show that the  methods discussed in the text
can be used to prove the existence of non-ruled smooth projective
Fano varieties over  arbitrary fields of characteristic zero. Without 
these explicit examples, our methods would {\it not\/} have 
established the existence of non-rational Fano varieties over, 
say, $\QQ$. The reason is that  the Morse Lemma fails over finite
fields, so after reducing modulo $p$, we could not guarantee any 
  polynomial with only non-degenerate critical points ({\it ie\/}
satisfying Assumption \refto{nondeg})  exists
with degrees satisfying the required
constraints.

The Proposition of this appendix shows that
 the polynomial
$$ \sum_{i=0}^{n}X_i^{2m-1}X_{i+1}
$$
 satisfies Assumption \refto{nondeg}
of Lecture 4. Furthermore, considered over $\QQ$, it is easy to check that 
that none of the critical points has critical value zero, so that the
hypersurface defined by 
$$Y^p = \sum_{i=0}^{n}X_i^{2m-1}X_{i+1}$$
in the weighted projective $\bold P$ is smooth over $\QQ$.

This  produces a host of  specific examples of non-rational
smooth projective Fano varieties of every dimension. For example, 
consider  the equation
$$
Y^2 - \sum_{i=0}^{n}X_i^{2m-1}X_{i+1} + 2G,
$$
for any even $n$ in the range $m < n+1 < 2m$, and $G$ is (weighted)
homogeneous of degree $2m$.
In the  weighted projective space $\bold P$ where $Y$ has weight $m$ and 
the $X_i$'s have weight one, this polynomial 
defines a non-rational Fano variety over $\QQ$, which is smooth for generic 
choices of $G$. In fact, 
it is smooth when $G=0$, giving a truly explicit
example.

\bigskip
\comment
\newpage

\setcounter {Ch}  {6}
\heading \  Appendix II: Proof of Manin  Theorem
 for Non-singular cubic surfaces. 
\endheading
\label {lec6}

Manin's theorem that two smooth cubic surfaces are birationally 
equivalent if and  only if they are projectively equivalent can be improved.
In this appendix, we sketch how to carry out the
proof when the surfaces may have some isolated non-smooth points,
 assuming that the cubic surfaces are {\it non-singular.\/} 

A variety is non-singular if all its local rings are regular
local rings.  This is not equivalent to smoothness in general. Indeed,
a  variety 
$X$ defined over $k$
 is smooth  over $k$ if and only if it 
is it is geometrically regular; that is, that the 
variety $X \times_{\spec k} \spec {\bar k} $ should be non-singular.
Of course, when $X$ is defined over an algebraically closed field, or even
a perfect field, there is no difference between the 
concepts of  non-singular and smooth.

If a cubic surface has only isolated non-smooth points, it can be said to 
be {\it geometrically normal:\/} after base change,
 $X \times_{\spec k} \spec {\bar k}$ is regular in codimension one.

\proclaim{Theorem}
Any two birationally equivalent 
 non-singular  cubic surfaces with only isolated non-smooth points
of Picard one are actually projectively equivalent.
\endproclaim

\demo{Proof}
The proof begins in precisely the same manner as the proof in the smooth
case presented in Lecture 3.  One has a  birational map 
$$
S \overset{\phi_{\gamma}}\to\rmap S'
$$ 
which may have some base points $P_1, \dots, P_r$, and one
computes numerical formulas for their multiplicities. {\it Assuming
the base points occur only at smooth points \/} of $S$, these numerical
computations are exactly the same as in the case where $S$ is smooth.
Fortunately, the base points of a birational map between non-singular 
complete surfaces occur only at smooth points,  which can be 
deduced easily from Zariski's
theorem on factoring birational maps by blow-ps and blow-downs 
\cite{K97, Cor 5.3.3}.

\endcomment
\newpage

\setcounter {Ch}  {7}
\heading   Solutions  for  Exercises
\endheading
\label {lec7}

\subhead{Exercise 1}\endsubhead
This was originally proved by Nishimura in 1955 \cite{Ni}.
The following proof is due to Endre Szab\'o.

Use induction on the dimension of $Y$. If $Y$ has dimension one, 
rational maps are morphisms defined everywhere, and the result is obvious.
If $Y$ is a smooth variety with a $k$-point $P$, blow up $P$ to get 
a variety $\tilde Y$. The blowup map $\tilde Y @>>> Y$ 
is defined over $k$, and the exceptional fiber, being isomorphic
to a projective space $\PP$, has lots of $k$-points. Any  rational map 
$Y \overset{\phi}\to\rmap Y'$  defined over $k$ determines a rational
map 
$\tilde Y \overset{\tilde{\phi}}\to\rmap  Y'$. Because 
$\tilde Y$ is smooth and $Y'$ is projective, 
the locus of indeterminacy has codimension at least two. 
This means that $\tilde\phi$ restricts to a rational map of the 
exceptional fiber 
$\PP$. Because this variety has smaller dimension, 
we are done by induction.

\smallskip
If $Y$ is not smooth, Nishimura's Lemma can fail. 
Indeed, let $Y$ be the projective closure of the affine
cone over a smooth projective variety $X$ with no
$k$-points.  Then $Y$ has exactly one $k$ point, the vertex
of the cone.  
Blowing up the vertex, we achieve
a smooth projective variety  $\tilde Y$ with no $k$-points, since the 
exceptional fiber is $k$-isomorphic to  $X$.
The rational  map $Y \rmap \tilde Y$ 
gives the counterexample to Nishimura's Lemma in the case where the source is
not smooth.

\subhead{Exercise 2}\endsubhead
(1) Choose coordinates so that  the disjoint $n$ planes $L_1$ and $L_2$ are 
given by $\{X_0 = X_1 = \dots = X_n = 0\}$ and
 by $\{X_{n+1} = X_{n+2} = \dots = X_{2n+1} = 0\}$
respectively.

A cubic given by an equation of the form 
$$
\sum_{i \leq n; j >n} a_{ijk} X_i X_j X_k
$$
obviously contains both planes. 
The generic member in this 
 linear system of cubics in $\PP^{2n+1}$
 is smooth, 
 because it admits  the following smooth special member:
$$
\sum_{i=0}^n(X_i^2X_{i+n+1} + X_iX_{i+n+1}^2).
$$
This is easily checked by the Jacobian 
criterion (assuming the characteristic is not 3).

\smallskip
(2)
To count the dimension of the linear system of cubics containing
a fixed pair of disjoint planes $L_1$ and $L_2$, we count
the number of monomials of degree 3 in $2n+2$ variables minus the number
of those monomials involving only variables 
generating the ideal of $L_1$ or of $L_2$.
The total is
$$
\binom{2n+1+3}{3} - \binom{n+3}{3} - \binom{n+3}{3} = (n+1)^2(n+2),
$$
so the dimension of the linear system is $(n+1)^2(n+2) - 1$.

Finally, to find the dimension of the space of all cubics containing
{\it any\/} pair of disjoint planes, we need to add the dimension
of the space of such pairs of planes. As a generic pair
of planes are disjoint,  this is twice the dimension of
the Grassmannian of $n$ planes in $\PP^{2n+1}$, or $2(n+1)^2$.
So the dimension of the space of all smooth 
 $2n$ dimensional cubic hypersurfaces
containing a pair of disjoint $n$-planes is 
$(n+1)^2(n+4) - 1$.

\smallskip
(3)
A cubic surface containing a linear subspace of dimension greater than $n$
is never smooth. Indeed, choosing coordinates so that 
$X$ contains the space defined by $ \{X_0 = X_1 = \dots = X_{n-1} = 0\}$,
$$
X = \left\{ \sum_{i=0}^{n-1} X_i f_i = 0 \right\}
$$
where $f$ are degree two. Now $X$ can not be smooth 
along the locus of points where
 $ \{X_0 = X_1 = \dots = X_{n-1} = 0\}$ and $\{f_0 = f_1 = \dots = f_{n-1} = 
0\}$, because the Zariski tangent space at these points 
is $2n+1$ dimensional.
 But because this locus is defined by only $2n$
equations, it must have non-empty intersection with $X$.

\subhead{Exercise 3}\endsubhead
Recall that a linear system of plane cubics with up to seven 
assigned base points  (including possibly one infinitely near another),
 no four collinear and no seven on a conic, 
has no unassigned base points \cite{H, p399}.
Also recall that a linear system on a smooth surface
 is very ample if and only if 
imposing two more base points (including one infinitely near another)
causes the dimension to drop by exactly two \cite{H p396}.

Consider four general points $P_1, \dots, P_4$ in $\PP^2$ and
let $ \beta = |3H - P_1 - P_2 - P_3 - P_4|$ be the linear system 
of cubics in $\PP^2$ passing through these points. Using the criterion above,
we see that (the pull back of) this
linear system, $\beta = 
 |3H - E_1 - E_2 - E_3 - E_4|$  to the blowup of $\PP^2$
at the four points 
is very ample.
 Using this linear system 
 embed the blowup as a surface  $S$ in $\PP^5$.

{\it Claim:\/}  
A generic projection of $S$ to $\PP^4$ is a surface 
$S'$ with exactly one singular point. 

First, for two general points $P$ and $Q$ on $\PP^2$,
 consider the following two linear subsystems of $\beta$ on $S$.
Fixing defining equations $s_0, \dots, s_5$ for generators of $\beta$, 
consider the linear subsystems whose defining equations satisfy:
$$
\align
\gamma := & \, 
\left\{ s \in \beta   \,\,\, | \,\, \frac{s(P)}{s_0(P)} =
 \frac{s(Q)}{s_0(Q)} \right\}\\
\alpha := & \,
\left\{ s \in \beta  \,\,\, | \,\, {s(P)} = {s(Q)}  = 0 \right\}.\\
\endalign 
$$
Note that $\alpha \subset \gamma \subset \beta$, and the dimensions
 drop by exactly one with each successive condition imposed.
The linear system $\gamma$ determines a projection $\pi$ of $\PP^5$ to 
$\PP^4$, sending $S$ to, say,  $S'$. By definition of $\gamma$, we have
$[s_0(P): s_1(P) : \dots : s_5(P)]  = [s_0(Q): s_1(Q) : \dots : s_5(Q)]$,
so that $\pi$ sends $P$ and $Q$ to the same point of $S'$. 

If $\pi$ collapses some other point  $P'$  (possibly infinitely near $P$ or
 $Q$)
to $\pi(P) = \pi(Q)$,
then we have that whenever $s(P) = s(Q) = 0$ for some $s \in \beta$,
the vanishing $s(P')$ is forced as well. This makes $P'$ an unassigned
base point of $\alpha$, contradicting the genericity assumption.
Likewise, if $\pi$ collapses  two other points $P'$ and $Q'$ 
($Q'$ may be infinitely near $P'$) to a single point
of $S'$, then  the linear system $|3H - P_1 - P_2 - P_3 - P_4 - P - Q - P'|$
has an unassigned base point $Q'$, again a contradiction, since six general
points and the one  special point $P'$ impose no extra conditions.

Thus any projection   of $S \subset \PP^5$ to $\PP^4$  that collapses 
two general points of $S$ to a single point  of $S'$
can collapse {\it only\/} these two points
to a single point. 
The argument will be complete once we have shown that a general projection 
 $\PP^5 \rmap \PP^4$ cannot be one-to-one on $S$.

Consider the incidence correspondence
$$
\Gamma = \{(P, Q, x) \, | \, P, Q, x {\text {  collinear}} \} \, \subset 
S \times S \times \PP^5.
$$
Through any two distinct  points of $S$,  there is a unique line in $\PP^5$,
so the projection $\Gamma \rightarrow S \times S$ is surjective,
 and its fibers
are all one-dimensional.
 It follows that $\Gamma$ is irreducible and of dimension $5$.

Consider the other projection $\Gamma \rightarrow \PP^5$. 
We know that if $P$ and $Q$ are collapsed to the same point of $S'$ via $\pi$,
then these are the only two points collapsed under $\pi$. 
This implies that the fiber
over any point in the image of $\Gamma  @>>> \PP^5$ is simply the triple
$(P, Q, x)$, so the fibers are zero-dimensional.
From this we conclude that $\Gamma  @>>> \PP^5$ is surjective.
This means a generic projection from any point in $\PP^5$ 
cannot be one-to-one on $S$. 

This implies that a generic projection of $S \subset \PP^5$ to 
a hyperplane in $\PP^5$ collapses precisely two points of $S$ to a
single point $S'$ in the image (which is therefore a singular point of $S'$).
This  completes the proof.

\smallskip
The surface we described is called a del Pezzo surface in $\PP^5$. 
The reader familiar with rational quartic scrolls in $\PP^5$
should be able to prove that these surfaces also have the property that 
a generic projection to $\PP^4$ produces exactly one double point.
In 1901, Severi claimed that these are the only  two examples
of such surfaces with a single ``apparent double point,'' as he called them
\cite{Sev, p 44}. Unfortunately, there was a gap in his argument. This gap
is being considered in the developing thesis of Mariagrazia Violo
(under the direction of Edoardo Sernesi)  which also  includes  further 
generalizations of these ideas \cite{V}.

\smallskip

\subhead{Exercise 4}\endsubhead
Use the  exact sequence 
$$
0 \rightarrow \Omega_{\PP^n} \rightarrow
\Cal O_{\PP^n}(-1)^{\oplus (n+1)} \rightarrow 
\Cal O_{\PP^n} \rightarrow 0
$$
(see \cite{H, p176} for the derivation of this sequence).
Since $\Omega_{\PP^n} \subset 
\Cal O_{\PP^n}(-1)^{\oplus (n+1)}$, it follows that 
 $(\Omega_{\PP^n})^{\otimes m} \subset 
(\Cal O_{\PP^n}(-1)^{\oplus (n+1)})^{\otimes m} \cong
\Cal O_{\PP^n}(-m)^{\oplus (n+1)^m}$.
Obviously now
$\Omega_{\PP^n}^{\otimes m}$ has no non-zero global sections,
 since $\Cal O_{\PP^n}(-m)$ 
has none.

\subhead{Exercise 5}\endsubhead
\,(1)\,\,
To prove that the plurigenera of a separably unirational variety $X$
vanish, reduce to the case of $\PP^n$ exactly as in the proof of Theorem
\refto{diffvan}. For $\PP^n$, 
we have $\bigwedge^n\Omega_{\PP^n} \cong \Cal O_{\PP^n}(-n-1)$ from the
above sequence, whence the vanishing of  $H^0(\Cal O_{\PP^n}(-mn -m))$
is immediate. This proves corollary 1.12.

\smallskip
(2) \,\, 
By the adjunction formula, the canonical class of a hypersurface $X$ in 
$\PP^n$ is $K_X = (K_{\PP^n} + X )|_X$,
so $\Cal O_X(K_X) = \Cal O_X((-n-1+d)H)$ where $X$ is degree $d$. 
So for 
$d > n$, the global sections of
all powers are non-zero and the plurigenera do not vanish.

\subhead{Exercise 6}\endsubhead \,(1)\,\,
A variety $X/k$ of dimension $d$ is unirational if and only if its
function field has an algebraic extension that is 
purely transcendental over $k$.
If $X' @>>> X$ is purely inseparable, then by definition,
the functional fields have the following relationship
$$
\{k(X')\}^{p^e} \subset k(X) \subset k(X')
$$
where $p^e$ is some power of the characteristic $p$.
If $X$ is unirational, then $k(X)$ has an algebraic extension
$k(t_1, \dots, t_d)$,
and hence $k(X') \subset k^{1/p^e}(t_1^{1/{p^e}}, \dots, t_d^{1/p^e})$.
Since $k$ is perfect, $k = k^{1/{p^e}}$,
 and $k(X')$ is a subfield
of a purely transcendental extension of $k$.
Thus $X'$ is unirational over $k$.

\smallskip
(2)\,\,
To construct a unirational variety of arbitrary degree, fix 
any polynomial $f$ of degree $p^e$ in $X_0, X_1, \dots, X_n$,
the homogeneous coordinates of $\PP^n$. Let $Y$ be an indeterminate, and let 
 $g$ be the polynomial $Y^{p^e} - f(X_0, X_1, \dots, X_n)$.
Then $g$
 defines a hypersurface in $\PP^{n+1}$, of degree $p^e$,
which is a purely inseparable cover of $\PP^n$. 
The previous computation shows that this
hypersurface is unirational. Indeed, its function field
is isomorphic to  $k(x_1, \dots, x_n)(f^{1/p^e})$, a subfield
of the purely transcendental extension $k(x_1 ^{1/p^e}, \dots, x_n^{1/p^e})$.
  
\subhead{Exercise 7}\endsubhead
 This is sometimes called Tsen's theorem \cite{Ts}.  The case $d = n = 2$ 
is due to Max Noether in 1871 \cite{No}.
We seek solutions $x_i = \sum_{j = 0}^{m} a_{ij} t^j$, where the $a_{ij}$ are 
unknown elements of $\CC$, to the degree $d$ polynomial $F(X_0, \dots X_n)$. 
Plugging in $X_i = x_i$, and gathering up all terms $t^r$, we see that
the coefficient of $t^{r}$ is a  polynomial in the unknowns 
$a_{ij}$.
We have a solution if and only if we can choose the $a_{ij}$ so that
 the coefficient of each  $t^r$  is zero.
Note that a  collection of complex numbers $a_{ij}$ is a solution
if and only if $\lambda a_{ij}$ is solution, where $\lambda$ is a 
non-zero scalar,
so the solutions naturally live in a  projective space over $\CC$.

Without loss of generality, we can assume that the coefficients of $F$ are
polynomials in $t$, say of degree less than $c$.
In this case, the highest occurring power of $t$ in
the expansion  
of  $F(x_0, \dots, x_n)$ is at most $dm + c$. Thus to solve for
the $a_{ij}$ is to solve a system of $dm + c$ 
equations in $m(n+1)$ unknowns. Since $n \geq d$, there are more
unknowns than equations when 
$m\gg 0$, so there are solutions for the $a_{ij}$ in projective space over
$\CC$.

The argument for $\CC(t, s)$ is similar.
Of course the same argument works with any algebraically closed field in
 place of
  $\CC$.
If $d > n$, the degree $d$ hypersurface in $\PP^n$ may have no $\CC(t)$ points.
For example, there are no $\CC(t)$ points of the variety defined by
 $\sum_{i = 0}^{d-1}t^iX_i^d = 0$ in $\PP^{d-1}$. 
By  projectivizing the affine cone over this
example, we get an example  of a degree $d$ hypersurface in $\PP^d$
that has exactly one
 $\CC(t)$-point,
 $[0: 0 : \dots : 0 : 1]$.

\subhead{Exercise 8}\endsubhead
Think of  $X_{a,2} \subset \A^1 \times \PP^n$ 
as defined by an equation $\sum_{ij} a_{ij}(t) X_i X_j$
where the $a_{ij} \in \CC(t)$ have 
  degree  less than or equal to $a$. 
This is a quadric in  $\PP^n_{\CC(t)}$ 
By exercise 7, we know this quadric has a $\CC(t)$-rational point
provided $n \geq 2$. This makes the quadric rational over $\CC(t)$,
whence its function field is isomorphic to $\CC(t)(x_1, \dots, x_n) 
\cong \CC(t, x_1, \dots, x_n)$.
This proves that $X_{a,2}$ is a  rational $\CC$-variety. 

Geometrically, the point is that  
 the projection  $X_{a, 2} @>>> \PP^1$ has a section, making the
family trivial on an open set.

The proof of (2) is similar.

\subhead{Exercise 9}\endsubhead
(1)
The variety $Y$ of $m \times n$ matrices of  rank at most $t$ is
defined by the $t+1$-minors of an $m \times n$   matrix of indeterminates.
It is easy to check that its dimension is $mn - (m-t)(n-t)$.
Define the rational map
$$
Y \rmap  \A^{mn - (m-t)(n-t)}
$$
$$
\bold{\lambda} \mapsto \{( \dots, \lambda_{ij},  \dots ) |
\,\, i {\text{ or }} j \leq t \},
$$
sending a matrix $\bold{\lambda}$ to the indicated string of its entries.
This map is a birational equivalence 
because,  whenever  the upper left hand $t$-minor $\Delta$ of the
$m \times n$ matrix  $\bold{\lambda}$
 is non-zero, we  solve uniquely
for each $\lambda_{ij}$ with both $i$ and $j$ greater than $t$. 
Indeed, since all $(t+1)$-minors vanish, 
we can use the Laplace expansion  to express 
any such $ \lambda_{ij}$  with $i, j > t$  as a polynomial in
the $\lambda$'s from the first $t$ columns and rows with denominators $\Delta$.
This proves that $Y$ is a rational variety over any field.

The singular locus of $Y$ is
the subvariety of matrices of rank strictly less than $t$.
Indeed, if an $m \times n$ matrix has rank $t$, than some 
$(t-1)$-minor is non-vanishing, so using the analogous map described above,
we can map an open subset of $Y$ containing this matrix
 isomorphically to an open subset of affine 
space. Thus every full rank $t$ matrix in $Y$ is a smooth point.
 Conversely, if some matrix
has rank less than $t$, all $t$-minors vanish, and considering the
Laplace expansion of the  $t+1$-minors defining $Y$,
we see easily that the Jacobian matrix is zero in this case. This says
that  the
tangent space at such a point has dimension $mn$ and the point is 
a singular point of $Y$.

\smallskip
(2)
 Let $X$ be the subvariety of $\PP^n$ defined by the vanishing of the
determinant of the  $n \times n$ matrix
 $\bold L$ of general linear forms in $n+1$
 variables. For each $n+1$ tuple $x = (x_0,  x_1,  \dots, x_n)$, 
consider $\bold{L}(x)$ as a linear map $k^n @>>> k^n$.
This defines a 
rational map 
$$
X \subset \PP^n \rmap  \PP^{n-1}
$$
$$
x = [x_0 : x_1 : \dots : x_n] \mapsto \{{\text{kernel of the matrix }}
 \bold {L}(x)\}.
$$
The genericity assumption guarantees that the
 matrix $\bold{L}$ has rank exactly $n-1$ generically on $X$. Thus, for
a generic $x \in X$, 
the kernel of the matrix $\bold{L(x)}$ is a  one dimensional subspace of
$k^n$, and  so determines
a well defined point in $\PP^{n-1}$. 

It is easy to check that this map is birational:
the genericity hypothesis on the linear forms guarantees that for 
distinct general elements $x, y$ in $X$, the matrices
$\bold{L}(x)$ and $\bold{L}(y)$ 
have distinct null spaces. 

As above,
the singular locus of $X$ is defined by the vanishing of the
 $(n-1) \times (n-1)$
subdeterminants of the matrix of linear forms.

\bigskip
(3) We now show that  every smooth cubic surface is determinantal. 
The earliest proof of this fact appears to be in 
an 1866 paper of  Clebsch, who credits Schr\"oter \cite{Cl}.
We give here two different proofs.{\footnote{Both are classical;
I am grateful to  I. Dolgachev for suggesting the 
first, which was worked out together
with J. Keum, R. Lazarsfeld, and C. Werner, and to 
 T. Geramita for suggesting the second, which 
can be found in \cite{Ger}.}}

\smallskip
{\it First Proof:\/}
This proof uses the configuration of the 
twenty seven lines on the cubic surface $S$.  We claim that there
are nine lines on the surface that can be represented in two different
ways as a union of three hyperplane sections. That is, there are
six different
 linear functionals $l_1, l_2, l_3, m_1, m_2, m_3$ on $\PP^3$
such that the hyperplane sections of $S$ determined by each is a union of
three distinct lines, and the nine lines obtained as 
hyperplane sections with the $l_i$'s
are the same nine lines obtained from the  the $m_i$'s. 
Assuming this for a moment, 
the 
cubics $l_1l_2l_3$ and $m_1m_2m_3$ both define the same subscheme  of $S$,
which means that up to scalar, these cubics agree on $S$. 
In other words, the cubic
$$
l_1l_2l_3 - \lambda m_1m_2m_3 
$$
is in the ideal generated by the cubic equation defining $S$, and hence
it must generate it. 
On th e other hand, the cubic
$l_1l_2l_3 - \lambda m_1m_2m_3$ is obviously the
determinant of the matrix
$$\matrix
l_1 & m_1 & 0 \\
0 & l_2 &  m_2\\
-\lambda m_3 & 0 & l_3
\endmatrix.
$$

The proof will complete upon establishing the existence
of the special configuration of lines. 
First recall that $S$ is the blow-up of six points
$P_1, P_2, \dots, P_6$  in $\PP^2$, no
three on a line and no five on a conic. We embed $S$ in $\PP^3$ using the
linear system of plane cubics through these six points. 
The twenty seven lines on $S \subset \PP^3$ are obtained as follows:
\roster
\item 
for each pair of two points $P_i$ and $P_j$, the birational transform of 
the line $\overline{P_iP_j}$  in $\PP^2$ joining them;
\item
for each point $P_i$, the birational transform of the conic $Q_i$ through 
the remaining five points;
\item
for each point $P_i$, the fiber $E_i$ over $P_i$.
\endroster
For any pair of indices $i, j$,
 the three lines $\overline{P_iP_j}$,
 $E_i$, and $Q_j$ form a (possibly degenerate)
triangle on $S$; Indeed, 
thinking of the hyperplane sections of $S$ as cubics in the plane through the
six points, this triangle is  the hyperplane section given by
the cubic obtained as the union of 
$\overline{P_iP_j}$ and $Q_j$.  
Now it is easy to find such a configuration. For instance,
the nine lines
$$
\{E_1, Q_2, \overline{P_1P_2}\} \cup
 \{E_2, Q_3, \overline{P_2P_3}\} \cup \{E_3, Q_1, \overline{P_1P_3}\}
$$
are the same as the nine  lines
$$
\{Q_1, E_2, \overline{P_1P_2}\} 
\cup \{Q_2, E_3, \overline{P_2P_3}\} \cup \{Q_3, 
 E_1,\overline{P_1P_3}\},
$$
with the groupings indicating the two different configurations of triangles.

\smallskip

{\it Second Proof:\/} This proof is more algebraic.
Let $\gamma$ be the linear system of plane cubics through six
the six points
used to embed the blowup up of $\PP^2$ at these six points as the cubic surface
in $\PP^3$.  A set of
 defining equations  $\{s_0, s_2, \dots, s_3\}$
 for a basis
of $\gamma$ is necessarily the homogeneous ideal $I$
 of functions vanishing at the
six points.
The blowup of the six points  is by definition the (closure of the) graph of 
 rational map $\PP^2 \overset{\phi_{\gamma}}\to\rmap \PP^3$
given by $x \mapsto [s_0(x): s_1(x): s_2(x): s_3(x)]$.

The
homogeneous coordinate
ring $R/I$
 of the six points 
has a resolution 
$$
0 @>>> R^3 @>A>> R^4 @>>> R @>>> R/I @>>> 0
$$ 
where $R = k[X_0, X_1, X_2]$ is the homogeneous 
coordinate ring of $\PP^2$, 
and $A$ is a 
 $4 \times 3$ matrix  whose 
 3-minors are precisely the generators $s_i$ for $I$. 
That defining ideals of
 codimension two subvarieties  of projective space have
resolutions of this form was proved by Hilbert in 1890 \cite{Hil};
nowadays we recognize it as a special case of 
the well-known Hilbert-Burch theorem; see \cite{E, p502}.

Let $Y_0, Y_1, Y_2, Y_3$ be homogeneous coordinates for $\PP^3$.
The entries of  the $3 \times 1$ matrix
$$
 A^{tr} {\bmatrix Y_0 \\ Y_1 \\ Y_2 \\ Y_3 \endbmatrix} = {\bmatrix H_1 \\ H_1
 \\ H_2 \endbmatrix}
$$
are bihomogeneous of degree $(1, 1)$
 in each set of variables $X_i$ and $Y_i$, and the
resulting closed subvariety of $\PP^2 \times \PP^3$ is the graph 
of the rational map $\PP^2 \rmap \PP^3$ defined by $\gamma$,
as one readily checks using  Cramer's rule for solving linear systems
of equations. 

The cubic surface $S$ is
the  projection of this graph to $\PP^3$.
Factoring  above matrix equation differently:
$$
 {\bmatrix H_1 \\ H_1
 \\ H_2 \endbmatrix} = B {\bmatrix X_0 \\ X_1 \\ X_2 \endbmatrix}
$$ 
where $B$ is a $3 \times 3$ matrix  of linear polynomials in the $Y_i$.
The determinant of $B$ is a degree three polynomial in the $Y_i$;
this cubic obviously vanishes on the projection  to $\PP^3$  of the 
subvariety of 
$\PP^2 \times \PP^3$
defined by the
$H_i$, and hence the determinant of $B$ must define the original
cubic
surface in $\PP^3$.

\subhead{Exercise 10}\endsubhead
We first prove an even stronger result for lower degree
hypersurfaces: 
 if $d < n$, then there is a {\it line\/} through every point on 
$X$. 

Without loss of generality, the point may be assumed to be the
affine origin $P = [1: 0: \dots : 0 : 0]$. 
Since $X$ passes through $P$, it has affine equation
$$
f_1(x_1, \dots, x_n) + f_2(x_1, \dots, x_n) + \dots + f_d(x_1, \dots, x_n)
$$
where each $f_i$ is homogeneous of degree $i$. 

Any line through $P$ has the form 
$$
(a_1t, a_2t,\dots, a_nt); \,\,\, t \,{\text{ ranging through }} k  
$$
where $(a_1, \dots, a_n)$ is a non-zero point on the affine plane.
To find such a line on $X$, we need 
$$
f_1(a_1t, \dots, a_nt) + f_2(a_1t, \dots, a_nt) + \dots + 
f_d(a_1t, \dots, a_nt) = 0
$$
for all $t$.
Because each polynomial $f_i(a_1t, \dots, a_nt)$
 is homogeneous of degree $i$ in 
$t$, 
we need $(a_1, \dots, a_n)$ such that  each $f_i(a_1, \dots, a_n) = 0$. 
Because $d < n$, the equations  $\{f_i = 0\}_{i=1}^d$ 
have a solution in $\PP^{n-1}$,
and we have found a line on $X$.

\medskip
When $d = n$, there may not be a line on $X$. 
Instead, we prove the following:
{\it Through every  point on a degree $n$ hypersurface in $\PP^n$, there
passes a plane conic.\/} 

Let $\Cal C$ be the variety of plane conics in $\PP^n$ passing through $P$.
The linear
system of plane conics forms a projective space of dimension 5,
and those passing through $P$  is a hyperplane in this space.
So 
the variety of
conics through $P$ is a $\PP^4$-bundle over the Grassmannian of
planes in $\PP^n$ through $P$. 
This Grassmannian has dimension $2(n-2)$, so the dimension
of the variety $\Cal C$ of conics is $2n$. 

The hypersurfaces of degree $n$ in $\PP^n$ passing through $P$
naturally form  a hyperplane $\Cal X$ in the  $\binom{2n}{n}-1$-dimensional
 projective space of all degree $n$ hypersurfaces in $\PP^n$. 
Consider 
the incidence correspondence
$$
\Gamma = \{(X, Q) \,|, Q \subset X\} \subset \Cal X \times \Cal C
$$
together with the two projections 
$\Gamma @>{\pi}>> \Cal X$ and $\Gamma @>{q}>> \Cal C$.

The elements in the 
 fiber of $\pi$  over a hypersurface $X \in \Cal X$ can be identified with 
the conics on $X$ through  $P$. 
In order to show that through every point on a degree $n$ (or less)
hypersurface in $\PP^n$  there passes a conic,  we need to show
that the projection $\Gamma @>{\pi}>> \Cal X$ is surjective.

We compute the dimension of $\Gamma$ using the other projection
 $\Gamma @>{q}>> \Cal C$.  Fix a conic $Q$ through $P$.
We need to compute the dimension of $\pi^{-1}(Q)$, the space
of degree $n$ hypersurfaces containing $Q$. Choose coordinates
so that $Q$ is given by $x_3 = x_4 = \dots = x_n = 0$ and
a  homogeneous degree 2 polynomial $g(x_0, x_1, x_2)$.
The hypersurfaces of degree $m$ containing $Q$ can be written 
uniquely in the form:
$$
\align
x_n h_n(x_0, \dots, x_n) &+ x_{n-1} h_{n-1}(x_0, \dots, x_{n-1}) + 
x_{n-2}h_{n-2}(x_0, \dots, x_{n-2}) +  \dots \\
& \dots + x_3h_3(x_0, x_1, x_2, x_3) +
g \cdot  h(x_0, x_1, x_2),
\endalign
$$
where the $h_i$ are homogeneous of degree $m-1$ and $h$ is a homogeneous
of degree $m-2$. 
When $m = n$, 
the space of hypersurfaces of this form is
of dimension $\binom{2n}{n} -2. $ 

Now, because $\Gamma$ and $\Cal X$ have the same dimension, the
projection map $\Gamma @>{\pi}>> \Cal X$ is surjective if the fiber
over some  point in the image is of dimension zero. 
So the proof is complete upon  exhibiting any  particular
 hypersurface of degree $n$ containing only finitely
many conics through a point  $P$. 
We leave it to the reader to verify that 
the hypersurface defined by $X_0^n - X_1X_2 \dots X_n$ contains
only finitely many conics through the point $[1: 1: \dots: 1 :1]$.  

\subhead{Exercise 11}\endsubhead
This  real cubic surface and its properties were considered in a
1962 paper
of Swinnerton-Dyer \cite{S-D}.

To see that the surface has two real components, 
consider the affine chart where $t = 1$. Setting
 $v^2 = x^2 + y^2$, we see that the surface is a 
 surface of revolution  for the elliptic curve
$$
v^2  = (4z-7) (z^2-2).
$$
Because the function $f(z) = (4z - 7) (z^2 - 2)$ has
three distinct  real roots, we know from Example \refto{nonrat}
that the curve, and hence the surface of revolution, has two disjoint
real components.  The two real components correspond to $z \geq 7/4$
and $|z| \leq \sqrt 2$ respectively. 

\smallskip
On the real component where $z \geq 7/4$, the $\QQ$-points are dense.
Indeed, this component contains $[x:y:z:t] = [1:1:2:1]$. The tangent plane
to the surface at this point intersects with the 
surface to produce an irreducible singular cubic on $S$. This curve is
rational over  $\QQ$, and its $\QQ$-points 
are dense among its $\RR$-points. In 
particular, the $\QQ$-values for $z$ are dense among all real values for 
$z \geq 7/4$. For each  of these 
fixed $\QQ$-values $z_0$, the plane $z = z_0$ intersects the surface $S$
in the circle $x^2 + y^2 = f(z_0)$.
This conic is
$\QQ$-rational 
and its  $\QQ$-points
 are dense among its $\RR$-points.
Thus the $\QQ$-points of $S$ are dense on the manifold component
where $z/t \geq 7/4$.

\smallskip
There are no $\QQ$-points on the component where $|z/t| \leq \sqrt 2$.
To see this, suppose that 
  $[x: y: z: t]$  is such a  $\QQ$-point,  where without loss of
generality, 
 $t$ and $z$ are assumed relatively prime integers, with $t > 0$. 
So 
$$
t(7t - 4z) (2t^2 - z^2) = (tx)^2 + (ty)^2
$$
is an integer which is the sum of two rational squares.
Thus any prime $p$ congruent to 3 modulo 4
that divides $t(7t - 4z) (2t^2 - z^2)$  must divide
it an even number of times.

Because $|\frac{z}{t}| \leq \sqrt 2$, each of the integer factors
$$
t, \,\,
(7t - 4z), \,\,  (2t^2 - z^2) 
$$
is {\it positive}. We claim that none is
 congruent to 3 modulo 4. Indeed, no prime $p$ congruent to 
3 modulo 4 can divide any one of these factors to an odd power.
For if some such $p$ does, then it must divide precisely two of the factors
an odd number of times. But because $t$ and $z$ are relatively prime, 
it follows that $t$ and $2t^2 - z^2$ are relatively prime, and the
only possible common prime factor of $t$ and $(7t - 4z)$ is 2. 
Furthermore,  if $p$ divides both $(7t - 4z)$ and  $(2t^2 - z^2)$,
then $p$ divides $ (8t + 7z)(7t - 4z)  -  28(2t^2 - z^2) = 17tz$. 
Since such $p$  divides neither $z$ nor $t$, the only possibility 
is $p = 17$, which is not congruent to 3 modulo 4.

Now if $t$ is even, then $z$ must be odd, but this would force 
$ (2t^2 - z^2) $ to be congruent to 3 modulo 4.
On the other hand, if  $t$ is odd, then it must be congruent to 1 modulo 4,
but this forces $(7t - 4z)$  to be congruent to 3 modulo 4. This
contradiction implies that there is no $\QQ$-rational point on the 
component of the surface where $|z/t| \leq \sqrt 2$.

\subhead{Exercise 12}\endsubhead
(1).
Let $T_i = \sum_{L_j \in O(i) } L_j$ where $O(1), \dots, O(r)$ are the orbits
of $G$ on the twenty seven lines of $S/{\bar k}$.
Suppose that $C = \sum_{j = 1}^{27} m_j L_j$ is a $G$-invariant effective 
curve, where
without loss of generality  $G$  is assumed finite.
Grouping the $L_i$ into orbits, we write $C = \sum_{i=1}^r C_i$ where
$C_i = \sum_{L_j \in O(i)} m_{ij}L_j$.
Consider the curve class
$$
|G|C = 
\sum_{g \in G} gC = \sum_{i=1}^r (\sum_{g \in G} gC_i).
$$
Each term $\sum_{g \in G} gC_i $
has each $L_j$ appearing with exactly
the same coefficient, namely  
 $ c = \frac{|G|}{|O(i)|}\sum_{j \in O(i)} m_{ij}$.
This implies that
$$
C = \frac{1}{|G|} \sum_{g\in G} gC  = \sum_{i=1}^r cT_i.
$$
This proves that 
the orbit sums $T_i$
generate the cone of curves for $S/k$.

\smallskip
(2). Let $C$ be an effective curve such that $C^2 > 0$.
 To show that $C$ is in the interior of the cone of curves, it is
sufficient to show that for any divisor $D$, the divisor $(C + \epsilon D)$
is also effective, for sufficiently small positive $\epsilon$. 

Note that $(C + \epsilon D)^2 = C^2 + 2\epsilon C\cdot D + \epsilon^2 D^2$,
so this self intersection number is positive for $\epsilon $
 sufficiently small.
Also if $H$ is any ample divisor,
then $H \cdot (C + \epsilon D) = H \cdot C + \epsilon (H \cdot D)$ 
is positive for small enough $\epsilon,$  because $H \cdot C > 0$.

Because $(C + \epsilon D)$ has positive self intersection and positive
 intersection with all ample divisors, it follows that $C+ \epsilon D$ is
effective, as is easily seen by applying Riemann-Roch to the divisors
$n(C + \epsilon D)$ for $n \gg 0$ (see \cite{H, p363}). This completes the
proof. 

\subhead{Exercise 13}\endsubhead
Consider a line $L$  in the ambient three-space, together with the
 smooth  cubic surface
$S$ defined by $u^3 = f(x, y)$. If the line
 $L$ lies on
$S$, then $L$ projects to a line triply tangent to the smooth
plane cubic
curve defined by $f(x, y) = 0$ in the $xy$-plane. 
Let $L'$ denote this projection, and suppose it has equation
$y = mx+b$. 
The line $L'$ is triply tangent if
and only if  $f(x, mx +b)$ is a cube of a linear form, say $(cx + d)^3$.
Whenever we have this perfect cube, there are three lines on $S$ 
projecting to $L'$. These three lines
have parametric equations 
$$
(x, \,\, mx+b, \,\,\omega(cx +d))
$$
where $\omega$ is one of the three cube roots of unity.

Because the plane cubic $\{ f(x, y) = 0\}$ has nine points of triple
tangency, all twenty seven lines in $S$ are constructed in this way.

\smallskip
We work out explicitly 
the lines on the Fermat surface given in homogeneous coordinates 
by 
$$a_0X_0^3  + a_1X_1^3 + a_2X_2^3 + a_3X_3^3.$$ 
Factoring the first two 
two terms $ a_0X_0^3 + a_1X_1^3$ completely
into distinct homogeneous linear polynomials $l_1l_2l_3$, and
likewise factoring $ a_2X_2^3 + a_3X_3^3 = m_1 m_2 m_3$, 
the  Fermat cubic  has the form
$$
l_1 l_2 l_3 + m_1 m_2 m_3 =  0.
$$
The linear factors are distinct because the surface is smooth.
This produces nine lines on the surface, defined by 
the nine different pairs of planes containing it:
$$
\{ l_i =  m_ j = 0 \}.
$$
By considering factorizations of the other two groupings of the terms
$(a_0X_0^3 + a_2X_2^3) + (a_1X_1^3 + a_3X^3_3)$ and 
$(a_0X_0^3 + a_3X_3^3) + (a_1X_1^3 + a_2X_2^3)$  we can produce the
 remaining eighteen
lines on the surface in similar fashion.

\smallskip 
The case where $u^2 = f(x, y)$ is similar. 
However, there are three distinct lines on the surface in the plane at 
infinity (meeting in an Eckardt point). The remaining 24 lines on $S$ 
project  to twelve lines in the plane
 $\{ u = 0\}$ tangent to the smooth curve defined by $\{f = 0\}$. 
These twelve lines can be found by solving for $m $ and $b$ 
such that $f(x, mx+ b)$ is a perfect square.

\medskip
We show that the cubic surface  
defined by $x_1^3 + x_2^3 + x_3^3 =  a$ (where $a$ is not a cube) is 
not rational over $\QQ$. 
By Segre's theorem, it suffices to show that Picard number is one, 
and by Theorem \refto{gal}, it is enough to show that 
no Galois orbit consists of disjoint lines on the surface. 

The  computation above indicates  that all lines are 
defined over the splitting field of $t^3 - a$ over $\QQ$. This splitting
field is degree six over $\QQ$, and 
is generated by $\omega$, a primitive third root of unity and a real third root
$\beta$ of $a$. The Galois group is the full group
$S_3$ of all permutations of the 
 roots $\beta$,  $\beta\omega$,
and  $\beta \omega^2$ of \, $t^3 - a$. 

Factoring the equation for the cubic 
$$
(x_1^3 + x_2^3) + (x_3^3 + ax_0^3) = 
(x_1 + x_2)(x_1 + \omega x_2)(x_1 + \omega^2 x_2) + 
(x_3 + \beta x_0)(x_3 + \beta \omega x_0)(x_3 + \beta \omega^2 x_0) 
$$ 
we consider the lines  as described above.
 
Any  line defined by 
$$
\{x_1 + x_2 = 0; \,\,\, x_3 + \beta \omega^{i} x_0 = 0\/\}
$$
for some $i = 0, 1, 2$ contains all other lines of this type in its orbit,
since the Galois group acts transitively on the $\beta \omega^i$. 
This orbit  consists of three line in the plane $\{x_1 + x_2 = 0\}$;
 in particular, the lines of this orbit are not disjoint.

Now consider the orbit of a  line defined by 
$$
\{x_1 +  \omega x_2 = 0; \,\,\, x_3 + \beta \omega^{i} x_0 = 0\/\}
$$
for some $i = 0, 1, 2$. The cyclic permutation 
$\beta \mapsto \beta \omega
\mapsto \beta \omega^2 \mapsto \beta$ fixes $\omega$. So the orbit of this
line contains three lines in the plane $\{x_1 +  \omega x_2 = 0\}$, 
and hence can not consist of disjoint lines. Furthermore, the permutation
 interchanging $\beta \omega$ and $\beta \omega^2$ 
sends these lines to lines of
the form 
$$
\{x_1 +  \omega^2 x_2 = 0; \,\,\, x_3 + \beta \omega^{i} x_0 = 0\/\}.
$$
The same cyclic 
permutation
now takes this line to all others of this form, 
that is, to all others in the plane $\{x_1 + \omega^2 x_2 = 0\}$. 
So the six lines in these two planes constitute another orbit.

Considering the other two  groupings of lines, we find that there
are two more orbits consisting of three lines in the same plane, 
and two more orbits consisting of six lines in two planes.
None of these six orbits consists of disjoint lines, 
so we conclude that the Picard number of the surface is one.

\medskip
\subhead{Exercise 14}\endsubhead
Note that $E^2 = -1$, and that,
 since $P$ has multiplicity $m$ on $C$, $C' \cdot E = m$. Thus
 $C^2 = (C'+mE)(C'+mE) = (C')^2 + 2mC'\cdot E + m^2 E^2 = (C')^2 + m^2$.
For the other equality, first verify that $K_{S'} \cdot E = -1$ using the
adjunction formula   $\deg K_E = (K_{S'} + E)\cdot E$.
Then compute that 
 $C \cdot K_S = (C' + mE)(K_{S'} - E) = C'\cdot K_{S'}  - C' \cdot E + 
m E \cdot K_{S'}- m E^2 = C' \cdot K_{S'} - m$.

\medskip

\subhead{Exercise 15}\endsubhead
The map $S' @>{q}>> \PP^2 = \PP(T_P\PP^3)$  can be described as follows:
for $Q \in S' $ but not in the exceptional fiber, thinking of $Q$ as
a point in $S$,   $q(Q)$ is the
line $L$ through $P$ and $Q$; for $Q$ in the exceptional fiber,
thinking of $Q$ as a direction at $P$, $q(Q)$ is
the line through $P$ in the direction of $Q$.

Clearly the fiber of $q$ over a point $L \in \PP^2$ 
consists of the two points
$Q_1$ and $Q_2$ that, together with $P$, make up the intersection 
$L \cap S$. Ramification occurs precisely when $Q_1 = Q_2$.
To find the equation of this ramification locus, choose coordinates so that
$P = [0:0:0:1]$ is the origin in an affine patch where the surface $S$ is
 given  by an equation of the form 
$$
 f_1(x, y, z) +  f_2(x, y, z) +  f_3(x, y, z) 
$$
with $f_i$ homogeneous of degree $i$.

A line $L$ through $P$ is given by parametric equations
$(at, bt, ct)$  corresponding to a   point
$[a:b: c]$ in $\PP^2$. 
The intersection points of this line with $S$ are given by
the solutions of the  equation
$$
tf_1(a, b, c) +  t^2f_2(a, b, c) +  t^3f_3(a, b, c) = 0
$$
The two solutions  (other than $t = 0$) 
define the fiber over $L$
under the map $q$. Ramification occurs when the two solutions 
are identical,  so 
is given by the discriminant.
 Thus the ramification locus in $\PP^2$ 
 is the  quartic  defined by the homogeneous equation 
$
f_2^2 - 4f_1f_3.
$

To see that this ramification locus is smooth,
note that 
because 
$S'$ is a degree two cover of $\PP^2$, 
locally $S'$ is defined by a quadratic polynomial of the
form $u^2 - g(s, t)$, where $s, t$ are local coordinates on $\PP^2$.
The ramification locus is locally
 defined by $g=0$. On the other hand, since the polynomial
$u^2 - g(s, t)$ defines a smooth variety (namely $S'$), the Jacobian 
criterion implies that  
 $g(s,t)$ also 
 defines a smooth variety in $\PP^2$.

\subhead{Exercise 16}\endsubhead
Consider the linear system $\gamma \subset |2H|$ of quadric sections
 on the cubic surface $S$ 
 passing through both $P_1$ and $P_2$ 
with multiplicity at least two.  If $P_2$ is infinitely
near $P_1$, ``passing through $P_2$" should be interpreted
as ``in the direction of $P_2$." 
This linear system contains the symmetric square of the pencil
of hyperplanes containing through $P_1$ and $P_2$;
together with the divisor $T_{P_1}S\cap S + T_{P_2}S \cap S$
(interpreted as  $2T_{P_1}S\cap S$ if $P_2$ is infinitely near) these divisors
generate all of $\gamma$.

The only base points of $\gamma$ are $P_1$ and $P_2$. Indeed, since the line
$\overline{P_1P_2}$ does not lie on $S$, the linear system of hyperplanes 
through $P_1$ and $P_2$ has 
 a unique third base point $Q$
where $\overline{P_1P_2}$ intersects $S$. 
But since $Q$ can not lie in $T_{P_1}S\cap S + T_{P_2}S \cap S$, 
we see that $\gamma$ has exactly two base points, $P_1$ and $P_2$.

The dimension of $\gamma$ is three. Indeed, if $s_1$ and $s_2$ are defining
equations for the linear system of hyperplanes through $\overline{P_1P_2}$,
then defining equations for the generators of $\gamma$ are
$$
s_1^2, \,\, s_1s_2, \,\, s_2^2, \,\, q
$$
where $q$ is a defining equation for $T_{P_1}S + T_{P_2}S$.

It is now easy to check that the image of the rational
map, which is defined over $k$, 
$$
S \rmap \PP^3
$$
is a (singular) quadric surface in $\PP^3$, and that the map is two-to-one
everywhere it is defined.

This two-to-one map allows us to define an involution of $S$: 
we interchange points of $S$ which map to the same point under the map
given by $\gamma$. It is clearly defined over the ground field $k$.

\subhead{Exercise 17}\endsubhead
Regardless of the characteristic of the ground field, 
the non-smooth locus of an affine
  hypersurface defined by $G$ is the locus defined
by $G$ and all its partial derivatives. 

In particular,
the non-smooth locus of the hypersurface defined by $y^p - f$ is
the closed set defined by the ideal generated by $y^p - f$, $py^{p-1}$
and the partial derivatives of $f$.
In characteristic zero, therefore,
any non-smooth point will have $y$ coordinate zero. 
So a non-smooth 
 point has the form $(y, x_1, \dots, x_n) =
 (0, \lambda_{1}, \dots, \lambda_{n})$, where all the partial derivatives
of $f(x_1, \dots, x_n)$ vanish at $(\lambda_{1}, \dots, \lambda_{n})$.
Since $y^p - f = 0$,  it must be that 
 $(\lambda_{1}, \dots, \lambda_{n})$ is a critical point of critical value
zero. 
Equivalently, 
 $(\lambda_{1}, \dots, \lambda_{n})$ is  a non-smooth
point of  the hypersurface in $n$ space defined by $f$. 
But a sufficiently general polynomial $f$ defines a 
smooth  hypersurface, so in characteristic zero  we expect 
a general hypersurface of the form $y^p - f$ to be smooth.

In characteristic $p$, the derivative with respect to $y$ vanishes,
so the non-smooth locus is defined by the ideal
$(y^p - f, \frac{\partial f}{\partial x_1}, \dots
 \frac{\partial f}{\partial x_n})$. 
Every critical point $(\lambda_1, \dots, \lambda_n)$
 of $f$ determines exactly one non-smooth point,
by setting $y_i = f(\lambda)^{1/p}$. 
(Of course, $f$ could fail to have critical points at all, as in the example
 $f = x_1 + x_2^{mp}$, but this means simply that the critical points are 
hiding at infinity.)
For a general $f$, 
the expected dimension of
the locus where all the $\frac{\partial f}{\partial x_i}$ vanish 
is zero. Thus a general hypersurface of the form $y^p - f$ in
characteristic $p$ has only isolated non-smooth points.

\subhead{Exercise 18}\endsubhead
The hypersurface in $\PP^{n+1}$ is given by a homogeneous polynomial
$$G = y^pt^{mp-p} - (F_0t^{mp} + F_1t^{mp-1} \dots + F_{mp-1}t + F_{mp}),
$$
where each $F_i$ is a homogeneous polynomial in $x_1, \dots, x_n$.
The singular locus is defined by the homogeneous ideal
$$
(G,  \frac{\partial G}{\partial t}, 
 \frac{\partial G}{\partial y},
 \frac{\partial G}{\partial x_1}, \dots,  
 \frac{\partial G}{\partial x_n}).
$$
Note that each derivative above is contained in the ideal
$(t, x_1, \dots, x_n)$. So the singular locus contains the
point $[y: x_1 : \dots : x_n : t] = [1 : 0 : \dots :0]$, regardless of
the characteristic.

\subhead{Exercise 19}\endsubhead
Consider a connection $\LL @>{\nabla}>> \LL \otimes \Omega_X$
on $X$. On an open set $U$ where $\LL$ is trivial, fix
an isomorphism $\Cal O_X(U) \cong  \LL(U)$, with $g \in \LL$ 
corresponding to 1 in $\Cal O_X$. 
Set $\nabla(g)  = g \otimes \eta$, for some one-form $\eta \in \Omega_X(U)$.
On $U$, we have 
$$
\LL @>{\nabla}>> \LL \otimes \Omega_X 
$$
$$
fg \mapsto   f\nabla(g) = g\otimes df =   g \otimes (df + f\eta)
$$
So we can think of $\nabla$ as a gadget that associates to the local
section $f$, the one-form $df + f\eta$.

Now, we can 
``differentiate a section of $\LL$  in any tangent
direction.'' Indeed, a tangent direction is interpreted as section of
the sheaf 
of derivations of $\Cal O_X$ to $\Cal O_X$.
Since  Der($\Cal O_X, \Cal O_X$) = Hom(${\Omega_X}$, ${\Cal O_X}$),  
 each  derivation $\theta$  produces a homomorphism
 $\Omega_X @>{\theta}>> \Cal O_X$. Its value on the  
 one form $df + f \eta$  can be considered the derivative
of $f$ in the direction of $\theta$. 
So the local section $fg$ of $\LL$ is sent to  the section 
$\theta(df + f\eta) \, g$ of $\LL$.

\smallskip
 (2). For any line bundle $\LL$, we can naively try to differentiate
local sections as follows. Over an open set $U$ where $\LL$ is trivial 
with fixed
generator $g$, if we are 
 given a section  $\theta$ of the sheaf of derivations, 
we  try sending  each section $f \cdot g$ of $\Cal O_X  \cdot g$ to 
 $(\theta  f) \cdot g \in \Cal O_X  \cdot g \cong \LL$. 
In general, of course, this does not lead to a globally  well defined
method for differentiating sections of $\LL$, 
because patching fails. Indeed, let $g_1$ and $g_2$ 
be local generators for $\LL$, related by the transition function 
$ g_1 = \phi g_2$. If $s$ is a local section of $\LL$, then
writing $s = fg_1  = (\phi f) g_2$, 
we see that $\theta(s)$ is well defined
if and only if $\theta(\phi f)  = \phi\theta(f)$. Using the Leibnitz rule for
derivations, we see that this 
is equivalent to $\theta(\phi) = 0$.
Thus, this naive approach to differentiating sections
gives a well defined  global connection on $\LL$ if and only if
$\LL$ admits transition functions that are killed by all derivations.
Because derivations annihilate any function that is a $p^{th}$ power, 
it follows that any 
line bundle $\Cal L$ that is a $p^{th}$ power of another line bundle $\Cal M$
admits this  natural connection, as the transition functions for 
$\Cal L$ can be taken to be $p^{th}$ powers of transition functions for 
$\MM$. 

\smallskip
Let $\LL = \Cal O_{\PP^n}(mp)$ and fix a global section $f$ of $\LL$.
A convenient choice of local trivialization
for $\LL$ is to let $U_i$ be the set where the homogeneous coordinate $x_i$
does not vanish. On $U_i$, we can consider $x_i^{mp}$ to be a local generator.
Thinking of $f$ as a homogeneous polynomial of degree $mp$ in the 
homogeneous coordinates for $\PP^n$, it has representation 
$(f/x_i^{mp}) x_i^{mp}$ on $U_i$.

\subhead{Exercise 20}\endsubhead
The ground field may be assumed algebraically closed.
We check the case where $n$ is even; the case where $n$ is odd is similar. 
At an isolated non-smooth point, after making suitable
linear changes of coordinates, we can assume the equation has the form
$$
y^p - x_1x_2 - \dots - x_{n-1}x_{n} - f_3(x_1, \dots, x_{n}).
$$
Blowing up the ideal generated by $y, x_1, \dots, x_{n}$,
the resulting scheme is covered by  affine patches where $x_i \neq 0$.
Consider one of these, say where $x_1 \neq 0$. There are  local coordinates
$y', x_1, x_2', x_3', \dots, x_n'$ where $x_1y' = y$, and $x_1x_i' = x_i$ for 
$i >1$. 
The blown up hypersurface is defined
by 
$$
{x_1}^{p-2}{y'}^p - x_2' - x_3'x_4' -  \dots x'_{n-1}x'_{n} - x_1
f'_1(x_1, x_2' \dots, x_{n}')
$$ 
where $f_1'$ has order at least one in $(x_1, x_2', \dots, x_{2n}')$.
This hypersurface is easily verified to be smooth, using the Jacobian 
criterion. Alternatively, it is sufficient to check that the exceptional 
divisor, namely the divisor defined by $x_1 = 0$ on this
hypersurface,  is smooth. This is obvious, since its equation is 
$x_2' -  x_3'x_4' - \dots - x'_{n-1}x'_{n}$ (or 
${y'}^p - x_2' - x_3'x_4' - \dots - x'_{n-1}x'_{n} $ when  $p = 2$). 

\subhead{Exercise 21}\endsubhead
Let $x_1, \dots, x_n$ be local coordinates around $P$.
In these coordinates,
polynomials with a critical point at $P$ all have the form
$$
f(x_1, \dots, x_n) = a + \sum_{i\leq j} a_{ij} x_{i}x_j + 
{\text {  higher order terms}}
$$
where $a_{ij} \in k$ (by considering a Taylor series expansion, for instance).
The Hessian of $f$ is the symmetric matrix
$$
A =\pmatrix
2a_{11} & a_{12} & \dots & a_{1n} \\
a_{12} & 2a_{22} & \dots & a_{2n} \\
\vdots & \vdots  & \ddots & \vdots  \\
a_{1n} & a_{2n} & \dots & 2a_{nn} \\
\endpmatrix
$$
and the invertibility of $A$ 
is equivalent to the non-vanishing of the determinant
of this symmetric matrix.
 (It is easy to verify that this is also  equivalent to 
the condition that the 
 $\frac{\partial f}{\partial x_i}$'s generate
the maximal ideal $(x_1, \dots, x_n)$ of $P$.) 

In any characteristic other than two, the determinant of a symmetric
matrix is a non-zero polynomial in the entries. Therefore,  on a Zariski
 open subset of the finite dimensional vector
space of quadratic polynomials in the $x_i$, the coefficient matrix  
 $a_{ij}$ has non-zero
 determinant.
Thus a 
``sufficiently general'' polynomial over an infinite field has only
non-degenerate critical points, assuming the
 characteristic is not two. 

In characteristic two, however, a symmetric matrix is also an
alternating matrix. In $n$ is odd, 
 it  always has determinant zero, so all critical points 
of $f$ are degenerate in this case.  If $n$ is even, however, 
the determinant of a symmetric $n \times n$ matrix 
 is a  non-zero polynomial, and again we conclude that a
generic $f$ in an even number of variable
 has non-degenerate critical points.  See \cite{J. p 332-335} for these
basic facts on alternating forms, convince  yourself by looking at the cases 
$n \leq 4$.  

\medskip
{\bf Exercise 22.\/}
The question is local, so assume $Z @>>>Y @>>> X$ 
 are maps of   affine schemes corresponding to the
ring maps 
$A @>>> B @>>> C$.  By our finite type assumption,
we know $B$ is a finitely generated $A$-algebra with generators, say 
$x_1, \dots, x_n$ and relations, say $f_1, \dots, f_n$. We can assume the
same number of generators and relations because both $A$ and $B$ are
regular of the same dimension.
Likewise,
 $C$ is a finitely generated $B$ algebra with generators
$y_1, \dots, y_m$ and relations,
say $g_1, \dots, g_m$.
 Thus  $x_1, \dots, x_n, y_1, \dots, y_m$ are $A$ algebra
generators for $C$, with relations 
 $f_1, \dots, f_n, g_1, \dots, g_m$.

In this case, the Jacobian ideal for $B$ over $A$ 
is the principal 
ideal  given by the 
 determinant of the $n \times n$ Jacobian matrix 
$$
(\frac{\partial f_i}{\partial x_j})
$$
and the Jacobian ideal for $C$ over $B$ is given by the determinant
of the $m \times m$ matrix
$$
(\frac{\partial g_i}{\partial y_j}).
$$
Because the Jacobian ideal of $C$ over $A$ is the determinant of the
$(m+n) \times (m+n)$ matrix 
$$
\pmatrix
\frac{\partial f_i}{\partial x_j}  & 0 \\
\frac{\partial g_i}{\partial x_j}  & \frac{\partial g_i}{\partial y_j}  \\
\endpmatrix
$$
(in block form), the result is immediate.

\enddocument